\documentclass{uqthesis}

\usepackage{epsf,latexsym,graphics,amsmath,amssymb,epsfig}
\usepackage{cite}

\newtheorem{theorem}{Theorem}[section]
\newtheorem{lemma}{Lemma}[section]

\newtheorem{definition}{Definition}[section]

\input epsf
\begin{document}

\bibliographystyle{plain}

\pagestyle{empty}
\vspace*{1.5cm}
\begin{center}
{\huge \bf Solutions to the Yang-Baxter Equation and Casimir Invariants for the Quantised Orthosymplectic Superalgebra}\\[3cm]

{\LARGE Karen Dancer \\[.5cm]}{\Large
        B.Sc. (Hons)\\[5cm]}
{\Large Centre for Mathematical Physics\\
        School of Physical Sciences\\
	The University of Queensland\\[2cm]
 A thesis submitted for the degree of Doctor of Philosophy \\[.5cm]
August, 2004}
\end{center}

\renewcommand{\baselinestretch}{1.5}
\normalsize


\newpage

\vspace*{1cm}
\begin{center}
{\bf Statement of Originality}
\end{center}

\noindent I declare that, to the best of my knowledge and belief, the work contained in this thesis is Karen Dancer's own work, except as acknowledged in the text.  Furthermore, this material has not been submitted, either in whole or in part, for a degree at this or any other university.

\vspace*{2cm}
\noindent Karen Dancer \hspace{7cm} Mark Gould


\newpage
\pagestyle{plain}
\pagenumbering{roman}

\vspace*{1cm}

\begin{center}
  {\LARGE \bf Acknowledgements}
\end{center}
\noindent First and foremost, my thanks go to my supervisors Mark Gould and Jon Links for their mathematical insight and assistance, their approachability, and for their neverfailing faith in me.  I am also very appreciative of the help and encouragement given to me by Maithili Mehta and the other mathematical physicists at UQ.  Lastly I wish to thank my friends for giving me joy and keeping me sane, and my family for their constant love and support.

\newpage

\vspace*{1cm}
\begin{center}
{\bf Abstract}
\end{center}

\noindent  For the last fifteen years quantum superalgebras have been used to model supersymmetric quantum systems.  A class of quasi-triangular Hopf superalgebras, they each contain a universal $R$-matrix, which automatically satisfies the Yang--Baxter equation.  Applying the vector representation to the left-hand side of a universal $R$-matrix gives a Lax operator.  These are of significant interest in mathematical physics as they provide solutions to the Yang--Baxter equation in an arbitrary representation, which give rise to integrable models.

\

\noindent In this thesis a Lax operator is constructed for the quantised orthosymplectic superalgebra $U_q[osp(m|n)]$ for all $m > 2, n \geq 2$ where $n$ is even.  This can then be used to find a solution to the Yang--Baxter equation in an arbitrary representation of $U_q[osp(m|n)]$, with the example of the vector representation given in detail.

\

\noindent  In studying the integrable models arising from solutions to the Yang--Baxter equation, it is desirable to understand the representation theory of the superalgebra.  Finding the Casimir invariants of the system and exploring their behaviour helps in this understanding.  In this thesis the Lax operator is used to construct an infinite family of Casimir invariants of $U_q[osp(m|n)]$ and to calculate their eigenvalues in an arbitrary irreducible representation.


\newpage
\tableofcontents
\newpage

\pagenumbering{arabic}
\setcounter{page}{1}

\chapter{Introduction}
        \noindent Generally speaking, mathematicians don't get arrested for espionage.  But then, Norwegian Sophus Lie (1842-1899) was often unusual.  A brilliant mathematician, he started a new field of study by introducing what were later named Lie algebras, of which the superalgebras used in this thesis are a generalisation.  Unfortunately, his mathematical insight reportedly exceeded his communication skills.  Perhaps this is why, leaving France during the Franco-Prussian war, he was arrested as a German spy, his mathematics notes believed to be top-secret coded documents!  Fortunately a French mathematician vouched for the innocence of both Lie and his notes, and Lie was safely released from prison \cite{WebLie}.

\

\noindent  The study of Lie algebras has advanced much since then, in part because of their interest to physicists.  Applications were known as early as the 1920s, with one of the earliest being the description of the electron configuration of atoms \cite{Weyl}.  While very useful in modelling non-commutative systems, Lie algebras have some unfortunate limitations.  In particular, during the drive for unified physical theories a model was sought for systems involving both bosons and fermions.  Lie algebras are not a viable option as some of the operators in such systems obey anti-commutation relations.

\

\noindent  The answer to this problem was to use a Lie superalgebra, originally known as a $\mathbb{Z}_2$-graded Lie algebra.  In this generalisation of a Lie algebra the operation is sometimes commutative and sometimes anti-commutative, depending on the grading of the operators involved.  The usual Serre relations for a Lie algebra \cite{Serre} are also altered, with many superalgebras containing higher order relations known as the extra Serre relations \cite{Yamane}.  Superalgebras were being examined as early as 1955 \cite{Frolicher, Nijenhuis}, and their involvement in the deformation of algebraic structures was investigated in the 1960s \cite{Gerstenhaber}, but they didn't become a prominent area of research until the 1970s \cite{Corwin}. This was when their relevance to quantum physics in the context of supersymmetries was recognised, with their application being to systems containing both bosonic and fermionic particles.  With a more complicated root system and representation theory than their non-graded counterparts, they presented quite a challenge to mathematicians and physicists. 

\

\noindent  Nonetheless, with four groups on three different continents competing in the exploration of superalgebras, progress was bound to be swift.  One of the important early problems was to classify all the finite-dimensional Lie superalgebras.  Unsurprisingly the honours went to Victor Kac, who completed the classification in 1977 \cite{Kac}.  Four infinite families of non-exceptional superalgebras were found, known as the $A, B, C$ and $D$ series (or type) superalgebras.  This thesis concentrates on solving problems for the $B$ and $D$ series, some of which have already been answered for the $A$ and $C$ series.  
 
\

\noindent The 1970s also saw the investigation of the enveloping algebras of Lie algebras \cite{Dixmier}.  Interesting in themselves, these polynomial algebras can also be ``$q$-deformed'' to produce more generalised algebras dependent on a complex parameter $q$ \cite{Drinfeld, Jimbo}. Several groups then extended the concept to superalgebras \cite{BGZ, Chaichian, Deguchi, Delius, Kulish}, with the results being referred to as either quantum supergroups or, more correctly, quantum superalgebras.  These form a class of quasi-triangular Hopf superalgebras, which implies they each admit a universal $R$-matrix, making them systems of significant interest.

\

\noindent The Yang--Baxter equation originally arose in McGuire's and Yang's studies of the many-body problem in one-dimension with repulsive delta-function interactions \cite{McGuire,Yang} and Baxter's solution of the eight-vertex model from statistical mechanics \cite{Baxter72}.  It has since appeared in the study of other exactly solvable lattice models \cite{Baxter}, knot theory \cite{Witten, Wadati} and the quantum inverse scattering method \cite{KS}, with a mathematical examination given in \cite{Jimbo89}.  By finding solutions to the Yang--Baxter equation in the affine extensions of quantum superalgebras and studying the representation theory we can construct new supersymmetric integrable models, which have a variety of physical applications.  

\

\noindent  One such application is in knot theory, where each representation gives rise to a link invariant \cite{LGZ,Zhang95}.  Constructing solutions to the Yang-Baxter equation is an essential step towards evaluating the invariants.  Another application is in strongly correlated electron systems.  As electrons are fermions, such systems are often supersymmetric.  Thus it is unsurprising that quantum superalgebras provide a suitable framework in which to work \cite{Arnaudon, Ge, GLZT, MR, GHLZ}.  One of the simplest examples is the $q$-deformed $t-J$ model \cite{Foerster, Gonzalez}, which describes a doped antiferromagnet, in which at each site of a one-dimensional lattice the occupancy of two electrons in different spin states is forbidden as a result of the on site Coulomb interaction.  For a certain choice of couplings this model is invariant with respect to the superalgebra $U_q[gl(2|1)]$, and the Hamiltonian can be derived through the quantum inverse scattering method.  Having more information about the higher order quantum superalgebras will assist in the study of more complex models.

\

\noindent Many of the applications of the Yang-Baxter equation arise in the spectral parameter dependent case.  Such solutions are associated with representations of affine quantum superalgebras; the representations of the $R$-matrices in these cases automatically satisfy the Yang-Baxter equation.  However even in the non-affine case the theory of quantum superalgebras is largely undeveloped.  In this thesis the Lax operator, which is the universal $R$-matrix with the vector representation acting on the first component, is constructed for the $B$ and $D$ type quantum superalgebras.  Previously the $R$-matrix with the vector representation acting on both components has been constructed \cite{Martins, Mehta}, but not the Lax operator.  In principle this could be calculated from the results of Khoroshkin and Tolstoy \cite{Khoroshkin}, but that would be difficult technically.

\

\noindent  When studying the representation theory of classical Lie algebras, understanding the central elements known as Casimir invariants proved very useful \cite{Perelomov, Popov, Nwachuku}.  Similarly, knowledge about the Casimir invariants of the quantum superalgebras will assist in the study of the integrable models. Thus we wish to find the Casimir invariants of the superalgebra, and also to calculate their eigenvalues in an arbitrary representation.  This has been done for the non-exceptional classical superalgebras \cite{Bincer,GS, Jarvis, Scheunert83}, but only for the A and C series quantum superalgebras \cite{GLZ, LinksZhang}.  In this thesis these results are extended to cover the quantised B and D type superalgebras.

\

\noindent  Chapter 2 provides an introduction to the mathematics used in the thesis.  It begins by setting up the classical orthosymplectic superalgebra, including the root system chosen, the generating elements, and their defining relations.  A $q$-deformation is then performed on the enveloping algebra to produce the quantised orthosymplectic superalgebra, which includes both the $B$ and $D$ series quantum superalgebras.  A brief introduction to the Yang--Baxter equation and universal $R$-matrices concludes this chapter.

\

\noindent  In the following chapter one of the properties of universal $R$-matrices is examined in the context of an arbitrary representation, leading to a set of simple generators and defining relations which uniquely determine a solution.  In Chapter 4 the other relevant $R$-matrix properties are checked, confirming that the solution is indeed a Lax operator.  This is in turn used to construct another, related Lax operator known as its opposite. The defining relations are also examined more closely to confirm they incorporate not only the standard, but also the higher order, $q$-Serre relations.

\

\noindent  An example of how to use the Lax operator to construct a solution to the Yang--Baxter equation in a particular representation is included as Chapter 5.  Although this is done only for the vector representation, exactly the same method can be used for any other representation.  The result agrees with a previously constructed $R$-matrix for the vector representation \cite{Mehta}.

\

\noindent  Finally, the Lax operator is used to construct Casimir invariants for the quantised orthosymplectic superalgebra.  This follows the method used in \cite{Bincer} and \cite{Scheunert83} for various classical superalgebras, which was adapted in \cite{LinksZhang} to cover the quantum superalgebra $U_q[gl(m|n)]$.  The calculations are more complex than in those cases, however, both because they include $q$-factors and because orthosymplectic superalgebras possess a more complicated root system than general linear superalgebras.

\chapter{The Construction of $U_q[osp(m|n)]$} \label{osp}
	\noindent To construct the quantised orthosymplectic superalgebra $U_q[osp(m|n)]$ we closely follow the method used in \cite{GouldZhang99} and \cite{GouldZhang00}. We begin by developing $osp(m|n)$ as a graded subalgebra of $gl(m|n)$.  The enveloping algebra of $osp(m|n)$ is then deformed to yield $U_q[osp(m|n)]$, which reduces to the original enveloping superalgebra as $q \rightarrow 1$.

\section{The Construction of $osp(m|n)$}

\normalsize
\noindent  We start with the standard generators $e^a_b$ of $gl(m|n)$, the $(m+n) \times (m+n)$-dimensional general linear superalgebra, whose even part is given by $gl(m) \oplus gl(n)$.  Now the commutator for a $\mathbb{Z}_2$-graded algebra satisfies the relation

\begin{equation*}
[A,B] =  - (-1)^{[A][B]}[B,A],
\end{equation*}

\noindent where $A,B$ are homogeneous operators and $[A] \in \mathbb{Z}_2$ is the grading of $A$.  In particular, the generators of $gl(m|n)$ satisfy the graded commutation relations 

\begin{equation*}
[e^{a}_{b},e^{c}_{d}] = \delta^{c}_{b} e^{a}_{d} - (-1)^{([a]+[b])([c]+[d])} \delta^{a}_{d} e^{c}_{b}
\end{equation*}

\noindent where 

\begin{equation*}
[a] =  
   \begin{cases}
     0, \qquad \; a=i,     & 1 \leq i \leq m, \\
     1, \qquad \; a = \mu, & 1 \leq \mu \leq n.
   \end {cases}
\end{equation*}

\noindent   Throughout the thesis we use Greek indices $\mu, \nu$ etc. \hspace{-4mm} to denote odd objects and Latin letters $i, j$ etc. \hspace{-4mm} for even indices.  If the grading is unknown, the usual $a, b, c$ etc. \hspace{-4mm} are used.  Which convention applies will be clear from the context.   We will only ever consider the homogeneous elements, but all results can be extended to the inhomogeneous elements by linearity.

\

\noindent   The orthosymplectic superalgebra $osp(m|n)$ is a subsuperalgebra of $gl(m|n)$ with even part equal to $o(m) \oplus sp(n)$, where $o(m)$ is the orthogonal Lie algebra of rank $m-2$ and $sp(n)$ is the symplectic Lie algebra of rank $n-1$.  The latter only exists if $n$ is even, so we set $n=2k$.  We also set $l = \lfloor \frac{m}{2} \rfloor$, so $m=2l$ or $m=2l+1$. 

\

\noindent  To construct $osp(m|n)$ we require an even non-degenerate supersymmetric metric $g_{ab}$.  Any can be used, but for simplicity's sake we choose $g_{ab} = \xi_{a} \delta^a_{\overline{b}}$, with inverse metric $g^{ba} = \xi_{b} \delta^a_{\overline{b}}$. Here

\begin{equation*}
\overline{a}=
  \begin{cases}
     m + 1 - a, & [a]=0, \\
     n + 1 - a, & [a]=1,
  \end {cases}
\qquad  \text{and} \quad
\xi_{a} = 
  \begin{cases}
     1, & [a] = 0, \\
     (-1)^a, & [a]=1.
  \end{cases}
\end{equation*}

\

\noindent Then the operators

\begin{equation*}
\sigma_{ab} = g_{ac} e^c_b - (-1)^{[a][b]} g_{bc} e^c_a = -(-1)^{[a][b]} \sigma_{ba}
\end{equation*}

\noindent generate the orthosymplectic superalgebra $osp(m|n)$.  These satisfy the commutation relations

\begin{multline*}
[\sigma_{ab}, \sigma_{cd}] = g_{cb} \sigma_{ad} - (-1)^{([a]+[b])([c]+[d])} 
  g_{ad} \sigma_{cb} \\
- (-1)^{[c][d]} \bigl( g_{db} \sigma_{ac} - (-1)^{([a]+[b])([c]+[d])} g_{ac} 
  \sigma_{db} \bigr).
\end{multline*}

\noindent This $\mathbb{Z}_2$-graded subalgebra actually arises naturally from considering the automorphism $\omega$ of $gl(m|n=2k)$ given by:

\begin{equation*}
\omega (e^a_b) = - (-1)^{[a]([a]+[b])} \xi_a \xi_b e^{\overline{b}}_{\overline{a}}.
\end{equation*}

\noindent This is clearly of degree 2, with eigenvalues $\pm 1$, so it gives a decomposition of $gl(m|n)$:

\begin{equation*}
gl(m|n) = \mathcal{S} \oplus \mathcal{T}, \text{ with } [\mathcal{S},\mathcal{S}] \subset \mathcal{S},\; [\mathcal{T},\mathcal{T}] \subset \mathcal{S} \text{ and } [\mathcal{S},\mathcal{T}] \subset \mathcal{T},
\end{equation*}

\noindent where

\begin{alignat*}{2}
&\omega (x) = x && \forall x \in \mathcal{S}, \\
&\omega (x) = -x &\qquad & \forall x \in \mathcal{T}.
\end{alignat*}

\noindent Here $\mathcal{T}$ is generated by operators

\begin{equation*}
T_{ab} = g_{ac} e^c_b + (-1)^{[a][b]} g_{bc} e^c_a = (-1)^{[a][b]} T_{ba},
\end{equation*}

\noindent while the fixed-point $\mathbb{Z}_2$-graded subalgebra $\mathcal{S}$ is generated by

\begin{equation*}
\sigma_{ab} = g_{ac} e^c_b - (-1)^{[a][b]} g_{bc} e^c_a = -(-1)^{[a][b]} \sigma_{ba},
\end{equation*}

\noindent so is simply the orthosymplectic superalgebra $osp(m|n)$.  As a more convenient basis for $osp(m|n)$ we choose the set of Cartan-Weyl generators, given by:

\begin{align}
\sigma^a_b &= g^{ac} \sigma_{cb} \notag \\
&= e^a_b - (-1)^{[a]([a]+[b])} \xi_a \xi_b e^{\overline{b}}_{\overline{a}}. \label{CW}
\end{align}

\noindent Then the Cartan subalgebra $H$ is generated by the diagonal operators

\begin{equation*}
\sigma^a_a = e^a_a - e^{\overline{a}}_{\overline{a}},
\end{equation*}

\noindent which satisfy $$[\sigma^a_a, \sigma^b_b]= 0, \qquad \forall a,b.$$

\

\noindent As a weight system, we take the set $\{ \varepsilon_i,\; 1 \leq i \leq m \} \cup \{ \delta_\mu,\; 1 \leq \mu \leq n \}$, where \mbox{$\varepsilon_{\overline{i}} = - \varepsilon_i$}, \mbox{$\delta_{\overline{\mu}} = - \delta_\mu$}.  Conveniently, when $m=2l+1$ this implies $\varepsilon_{l+1} = -\varepsilon_{l+1} = 0$.  Acting on these weights, we have the invariant bilinear form defined by:

\begin{equation*}
(\varepsilon_i, \varepsilon_j) = \delta^i_j, \quad (\delta_\mu, \delta_\nu) = -\delta^\mu_\nu, \quad (\varepsilon_i, \delta_\mu) = 0, \qquad 1 \leq i, j \leq l, \quad 1 \leq \mu, \nu \leq k.
\end{equation*}

\noindent When describing an object with unknown grading indexed by $a$ the weight will be described generically as $\varepsilon_a$.  This should not be assumed to be an even weight.

\

\noindent  The even positive roots of $osp(m|n)$ are composed entirely of the usual positive roots of $o(m)$ together with those of $sp(n)$, namely:

\begin{alignat*}{3}
&\varepsilon_i \pm \varepsilon_j, & \qquad & 1 \leq i < j \leq l, \\
&\varepsilon_i, && 1 \leq i \leq l &&\quad \text{when } m=2l+1, \\
&\delta_\mu + \delta_\nu, && 1 \leq \mu,\nu \leq k, \\
&\delta_\mu - \delta_\nu, && 1 \leq \mu < \nu \leq k. 
\end{alignat*}

\noindent  The root system also contains a set of odd positive roots, which are:

\begin{equation*}
\delta_\mu + \varepsilon_i,  \qquad  1 \leq \mu \leq k,\;1 \leq i \leq m.
\hspace{2cm}
\end{equation*}

\noindent Throughout this thesis we choose to use the following set of simple roots:

\begin{align*}
&\alpha_i = \varepsilon_i - \varepsilon_{i+1}, \hspace{11mm}  1 \leq i < l, 
  \notag \\
&\alpha_l =  
\begin{cases} 
\varepsilon_{l} + \varepsilon_{l-1},\quad & m=2l, \\
\varepsilon_l,   &m=2l+1, 
\end{cases} \notag \\
&\alpha_\mu = \delta_\mu - \delta_{\mu+1}, \hspace{9mm} 1 \leq \mu < k,\notag\\
&\alpha_s = \delta_k - \varepsilon_1.
\end{align*}

\noindent Note this choice is only valid for $m >2$. 

\

\noindent  Corresponding to these simple roots we have raising generators $e_a$, lowering generators $f_a$ and Cartan elements $h_a$ given by:

\begin{alignat*}{4}
&e_i = \sigma^i_{i+1}, &\quad& f_i = \sigma^{i+1}_i, &\quad& h_i = \sigma^i_i 
  - \sigma^{i+1}_{i+1}, &\quad&1 \leq i < l, \notag \\
&e_l = \sigma^{l-1}_{\overline{l}}, && f_l = \sigma^{\overline{l}}_{l-1}, &&
  h_l = \sigma^{l-1}_{l-1} + \sigma^l_l, &&m=2l, \notag \\
&e_l = \sigma^l_{l+1}, && f_l = \sigma^{l+1}_l, && h_l = \sigma^l_l, && m=2l+1,
  \notag \\
&e_\mu = \sigma^\mu_{\mu+1}, && f_\mu = -\sigma^{\mu+1}_\mu, && h_\mu = 
  \sigma^{\mu+1}_{\mu+1} - \sigma^\mu_{\mu}, && 1 \leq \mu < k, \notag \\
&e_s = \sigma^{\mu=k}_{i=1}, &&f_s = -\sigma^{i=1}_{\mu=k}, &&h_s = 
  -\sigma^{\mu=k}_{\mu=k} - \sigma^{i=1}_{i=1}.&&
\end{alignat*}

\noindent These automatically satisfy the defining relations of a Lie superalgebra, which are:

\begin{alignat}{2} 
&[h_a, e_b] = (\alpha_a, \alpha_b) e_b, && \notag \\
&[h_a, f_b] = - (\alpha_a, \alpha_b) f_b, &&  \notag \\
&[h_a, h_b] = 0,&& \notag \\
&[e_a, f_b] = \delta^a_b h_a, \notag \\
&[e_a, e_a] = [f_a,f_a]=0 & \quad &\text{for } (\alpha_a, \alpha_a)=0,\notag\\
&(ad\,e_b\, \circ)^{1-a_{bc}} e_c = 0 && \text{for } b \neq c, \label{Ser1} \\
&(ad\,f_b\, \circ)^{1-a_{bc}} f_c = 0 && \text{for } b \neq c, \label{Ser2}
\end{alignat}

\noindent where the $a_{bc}$ are the entries of the corresponding Cartan matrix,

\begin{equation*}
a_{bc} =
\begin{cases}
   \frac{2(\alpha_b, \alpha_c)}{(\alpha_b, \alpha_b)}, \quad & (\alpha_b, \alpha_b) \neq 0, \\
(\alpha_b, \alpha_c), &(\alpha_b, \alpha_b) = 0,
\end{cases}
\end{equation*}

\noindent and $ad$ represents the \textit{adjoint action}

\begin{equation*}
ad\, x \circ y = [x,y].
\end{equation*}

\noindent The relations \eqref{Ser1} and \eqref{Ser2} are known as the \textit{Serre relations} \cite{Serre}.  Superalgebras also have higher order defining relations, not included here, which are known as the \textit{extra Serre relations}.  They are dependent on the structure of the root system \cite{Yamane}.


\section{The $q$-Deformation: $U_q[osp(m|n)]$}

\noindent A quantum superalgebra is a more generalised version of a classical superalgebra involving a complex parameter $q$, which reduces to the classical case as $q \rightarrow 1$.  In particular, we construct $U_q[osp(m|n)]$ by \textit{$q$-deforming} the original enveloping algebra of $osp(m|n)$ so that the generators remain unchanged, but are now related by a quantised version of the defining relations.  

\

\noindent First note that in the enveloping algebra of $osp(m|n)$ the commutator is given by

\begin{equation*}
[A,B] = AB - (-1)^{[A][B]} BA.
\end{equation*}

\noindent With this operation, the defining relations for $U_q[osp(m|n)]$ are:

\begin{alignat}{2}
&[h_a, e_b] = (\alpha_a, \alpha_b) e_b, && \notag \\
&[h_a, f_b] = - (\alpha_a, \alpha_b) f_b, &&  \notag \\
&[h_a, h_b] = 0,&& \notag \\
&[e_a, f_b] = \delta^a_b \frac{(q^{h_a} - q^{-h_a})}{(q - q^{-1})},&&
   \notag \\
&[e_a, e_a] = [f_a,f_a]=0 & \quad &\text{for } (\alpha_a, \alpha_a)=0,\notag\\
&(ad\,e_b\, \circ)^{1-a_{bc}} e_c = 0 &&  \text{for }b \neq c, \label{qS1}\\
&(ad\,f_b\, \circ)^{1-a_{bc}} f_c = 0 && \text{for } b \neq c \label{qS2}.
\end{alignat}

\noindent The relations \eqref{qS1} and \eqref{qS2} are called the \textit{$q$-Serre relations}.  Again, there are also extra $q$-Serre relations which are not included here.  A complete list of them, including those for affine superalgebras, can be found in \cite{Yamane}.  Both the standard and extra $q$-Serre relations depend on the adjoint action, which is no longer simply the commutator.  To define the adjoint action for a quantum superalgebra, we first need some new operations.

\

\noindent The \textit{coproduct}, $\Delta: U_q[osp(m|n)]^{\otimes 2} \rightarrow U_q[osp(m|n)]^{\otimes 2}$, is the superalgebra homomorphism given by:

\begin{align}
&\Delta (e_a) = q^{\frac{1}{2}h_a} \otimes e_a + e_a \otimes q^{-\frac{1}{2}
  h_a}, \notag \\
&\Delta (f_a) = q^{\frac{1}{2}h_a} \otimes f_a + f_a \otimes q^{-\frac{1}{2}
  h_a},\notag \\
&\Delta (q^{\pm \frac{1}{2}h_a}) = q^{\pm \frac{1}{2}h_a} \otimes q^{\pm 
  \frac{1}{2}h_a}, \notag \\
&\Delta (ab) = \Delta(a) \Delta(b). \label{coprod}
\end{align}

\noindent  Note that in a $\mathbb{Z}_2$-graded algebra, multiplying tensor products induces a grading term, according to

\begin{equation*}
(a \otimes b) (c \otimes d) = (-1)^{[b][c]} (ac \otimes bd).
\end{equation*}

\noindent  We also require the \textit{antipode}, $S: U_q[osp(m|n)] \rightarrow U_q[osp(m|n)]$, a superalgebra anti-homomorphism defined by:

\begin{align*}
&S (e_a) = - q^{-\frac{1}{2}(\alpha_a, \alpha_a)} e_a, \notag \\
&S (f_a) = - q^{\frac{1}{2}(\alpha_a, \alpha_a)} f_a, \notag \\
&S (q^{\pm h_a}) = q^{\mp h_a}, \notag \\
&S (ab) = (-1)^{[a][b]} S(b) S(a).
\end{align*}

\noindent  It can be shown that both the coproduct and antipode are consistent with the defining relations of the superalgebra.  These mappings are necessary to define the adjoint action for a quantum superalgebra, as it can no longer be written simply in terms of the commutator.  If we adopt Sweedler's notation for the coproduct,

\begin{equation*}
\Delta(a) = \sum_{(a)} a^{(1)} \otimes a^{(2)},
\end{equation*}

\noindent the \textit{adjoint action} of $a$ on $b$ is defined to be

\begin{equation*} \label{adj}
ad\; a \circ b = \sum_{(a)} (-1)^{[b][a^{(2)}]} a^{(1)} b S(a^{(2)}).
\end{equation*}

\noindent The added $q$-factors in the defining relations ensure that working with quantum superalgebras is significantly more difficult than with their classical counterparts, even though in this case the generators and root system remain the same.  Throughout the thesis $q$ is assumed not to be a root of unity.

\

\noindent One quantity that repeatedly arises in calculations for both classical and quantum Lie superalgebras is $\rho$, the \textit{graded half-sum of positive roots}.  In the case of $U_q[osp(m|n)]$ it is given by:

\begin{equation*}
\rho = \frac{1}{2} \sum_{i=1}^l (m-2i) \varepsilon_i + \frac{1}{2} \sum_{\mu=1}^k (n-m+2-2\mu) \delta_\mu.
\end{equation*}

\noindent This satisfies the property $(\rho, \alpha) = \frac{1}{2} (\alpha, \alpha)$ for all simple roots $\alpha$.

\

\noindent As mentioned earlier, this root system and set of generators is only valid for $m > 2$.  When $m=0$, $U_q[osp(m|n)]$ is isomorphic to $U_q[sp(n)]$.  Similarly, in \cite{Zhang} it was shown that every finite dimensional representation of $U_q[osp(1|n)]$ is isomorphic to a finite dimensional representation of $U_{-q}[so(n+1)]$.  As we are only interested in finite dimensional representations, and the representation theory of these non-super quantum groups is well-understood, we need not consider the cases with $m<2$.  Thus although our root system is only valid for $m >2$, finding the Lax operator for this root system will actually complete the work for all $B$ and $D$ type quantum superalgebras.  This has, of course, already been done for the more straightforward $A$ type quantum supergroups, $U_q[gl(m|n)]$ \cite{Zhang2}.  The Lax operator has yet to be constructed for the $C$ type quantum supergroups ($U_q[osp(m|n)]$ where $m=2$), although an $R$-matrix for the vector representation is known \cite{Scheun}.

\section{$U_q[osp(m|n)]$ as a Quasi-Triangular Hopf \\
         Superalgebra}

\noindent A quantum superalgebra is actually a specific type of quasi-triangular Hopf superalgebra.  This guarantees the existence of a universal $R$-matrix, which provides a solution to the quantum Yang--Baxter equation.  Before elaborating, we need to introduce the graded twist map.

\

\noindent  The \textit{graded twist map} $T:U_q[osp(m|n)]^{\otimes 2} \rightarrow U_q[osp(m|n)]^{\otimes 2}$ is given by

\begin{equation*}
T(a \otimes b) = (-1)^{[a][b]} (b \otimes a).
\end{equation*}

\noindent For convenience $T \Delta$, the twist map applied to the coproduct, is denoted $\Delta^T$.

\
 
\noindent Then a \textit{universal $R$-matrix}, $\mathcal{R}$, is an even, non-singular element of $U_q[osp(m|n)]^{\otimes 2}$ satisfying the following properties:

\begin{align}
&\mathcal{R} \Delta (a) = \Delta^T (a)\mathcal{R}, \quad \forall a \in 
  U_q[osp(m|n)], \notag \\
&(\text{id} \otimes \Delta) \mathcal{R} = \mathcal{R}_{13} \mathcal{R}_{12},
  \notag \\
&(\Delta \otimes \text{id}) \mathcal{R} = \mathcal{R}_{13} \mathcal{R}_{23}. 
  \label{Requations}
\end{align}

\noindent Here $\mathcal{R}_{ab}$ represents a copy of $\mathcal{R}$ acting on the $a$ and $b$ components respectively of $U_1 \otimes U_2 \otimes U_3$, where each $U$ is a copy of the quantum superalgebra $U_q[osp(m|n)]$.  When $a>b$ the usual grading term from the twist map is included, so for example $\mathcal{R}_{21} =  [\mathcal{R}^T]_{12} $, where $\mathcal{R}^T = T\, \mathcal{R}$ is the \textit{opposite universal $R$-matrix}.

\

\noindent  One of the reasons $R$-matrices are so significant is that as a consequence of \eqref{Requations} they satisfy the quantum Yang--Baxter Equation, which is prominent in the study of integrable systems \cite{Baxter}:

\begin{equation*}
\mathcal{R}_{12} \mathcal{R}_{13} \mathcal{R}_{23} = \mathcal{R}_{23} 
  \mathcal{R}_{13} \mathcal{R}_{12} 
\end{equation*}

\

\noindent A superalgebra may contain many different universal $R$-matrices, but there is always a unique one belonging to $U_q[osp(m|n)]^- \otimes U_q[osp(m|n)]^+$, and its opposite $R$-matrix in $U_q[osp(m|n)]^+ \otimes U_q[osp(m|n)]^-$.  Here $U_q[osp(m|n)]^-$ is the Hopf subsuperalgebra generated by the lowering generators and Cartan elements, while $U_q[osp(m|n)]^+$ is generated by the raising generators and Cartan elements.  These particular $R$-matrices arise out of Drinfeld's double construction \cite{Drinfeld}.  In this thesis we consider the universal $R$-matrix belonging to $U_q[osp(m|n)]^- \otimes U_q[osp(m|n)]^+$.

\chapter{Construction of the Lax operator for $U_{q}[osp(m|n)]$} \label{R-mat}
	\noindent In this chapter we construct a \textit{Lax operator} for $U_q[osp(m|n)]$.  Previously this had been done only for $U_q[gl(m|n)]$ \cite{Zhang2}. Before defining a Lax operator, however, we need to introduce the vector representation.  

\

\noindent Let $\text{End} \; V$ be the set of endomorphisms of $V$, an $(m+n)$-dimensional vector space.  Then the irreducible \textit{vector representation} $\pi: U_q[osp(m|n)] \rightarrow \text{End} \; V$ is left undeformed from the classical vector representation of $osp(m|n)$, which acts on the Cartan-Weyl generators given in equation \eqref{CW} according to:

\begin{equation*}
\pi (e^a_b) = E^a_b,
\end{equation*}

\noindent where $E^a_b$ is the $(m+n) \times (m+n)$-dimensional elementary matrix with $(a,b)$ entry $1$ and zeroes elsewhere.

\

\noindent Now let $\mathcal{R}$ be a universal $R$-matrix of $U_{q}[osp(m|n)]$ and $\pi$ the vector representation.  The Lax operator associated with $\mathcal{R}$ is given by

\begin{equation*}
R = (\pi \otimes \text{id}) \mathcal{R} \in (\text{End} \; V) \otimes U_{q}[osp(m|n)].
\end{equation*}

\noindent Previously only an $R$-matrix in the vector representation, $(\pi \otimes \pi) \mathcal{R}$, has been found, with it having been calculated for both $U_q[osp(m|n)]$ and its affine extension \cite{Martins, Mehta}.  The Lax operator is significant because we can use it to calculate solutions to the quantum Yang--Baxter equation for an arbitrary finite-dimensional representation.

\

\noindent In this chapter we also sometimes make use of the bra and ket notation.  The set $\{| a \rangle \}$ is a basis for $V$ satisfying the property

\begin{equation*}
E^a_b |c \rangle = \delta^c_b |a \rangle.
\end{equation*}

\noindent The set $\{ \langle a| \}$ is the dual basis such that 

\begin{equation*}
\langle c| E^a_b = \delta^a_c \langle b| \quad \text{and} \quad \langle a|\, | b \rangle = \langle a|b \rangle = \delta^a_b.
\end{equation*}   

\section{Developing the Governing Relations} \label{Developing Relations}

\noindent As we wish to find the Lax operator belonging to $ \pi \bigl( U_q[osp(m|n)]^-\bigr)  \otimes U_q[osp(m|n)]^+$, we adopt the following ansatz for $R$:

\begin{equation*} \label{RmatrixAnsatz}
R \equiv q ^ {\underset{a}{\sum} h_{a} \otimes h^{a}} \Bigl[ I \otimes I + (q - q^{-1}) \sum_{\varepsilon_{a} < \varepsilon_{b}} (-1)^{[b]} E^a_b \otimes \hat{\sigma}_{ba} \Bigr].
\end{equation*} 

\noindent Here $\{ h_a \}$ is a basis for the Cartan subalgebra such that $h_a = h_{\varepsilon_a}$, and $\{ h^a \}$ the dual basis, so $h^a = (-1)^{[a]} h_{\varepsilon_a}$. The $\hat{\sigma}_{ba}$ are the unknown operators for which we are trying to solve.  Throughout this chapter when working in the vector representation we simply use $h_a$ rather than $\pi(h_a)$, and $e_a$ rather than $\pi(e_a)$.

\

\noindent Now $R$ must satisfy the defining relations for an $R$-matrix, which were given as equation \eqref{Requations} in the previous chapter.  In particular, we begin by considering the intertwining property for the raising generators,

\begin{equation*}
R \Delta (e_c) = \Delta^T (e_c) R.
\end{equation*}

\noindent To apply this, recall that

\vspace{-10mm}
\begin{equation*}
\Delta(e_c) = q^{\frac{1}{2} h_c} \otimes e_c + e_c \otimes q^{-\frac{1}{2} 
  h_c}.
\end{equation*}

\noindent But
\begin{alignat*}{2}
& & \qquad & [h_a, e_c] = (\alpha_c, \varepsilon_a)e_c \notag \\
&\Leftrightarrow & \qquad & e_{c} h_a = [h_a - (\alpha_c, \varepsilon_a)] e_{c}
  \notag \\
&\Leftrightarrow & \qquad & e_{c} q^{h_a} = q^{[h_a -(\alpha_c, 
  \varepsilon_a)]}e_c.
\end{alignat*}

\noindent Hence

\begin{align*}
\Delta^T(e_c) q^{\underset{a}{\sum} h_a \otimes h^a} &= (e_{c} \otimes 
   q^{\frac{1}{2} h_c} + q^{-\frac{1}{2} h_c} \otimes e_{c} ) q^{\underset{a}
   {\sum} h_a \otimes h^a} \notag \\
&= q^{\underset{a}{\sum} [h_a - (\alpha_c,\varepsilon_a)I] \otimes h^a} e_{c} 
   \otimes q^{\frac{1}{2}h_c} + q^{\underset{a}{\sum} h_a \otimes [h^a - 
   (\alpha_c, \varepsilon_a)I]} q^{-\frac{1}{2} h_c} \otimes e_{c} \notag \\
&= q^{\underset{a}{\sum} h_a \otimes h^a} \bigl[ q^{-\underset{a}{\sum} 
   (\alpha_c,\varepsilon_a) I \otimes h^a} (e_c \otimes q^{\frac{1}{2} h_c}) + 
   q^{-\underset{a}{\sum} (\alpha_c,\varepsilon_a) h_a \otimes I} 
   (q^{-\frac{1}{2} h_c} \otimes e_c) \bigr] \notag \\
&= q^{\underset{a}{\sum} h_a \otimes h^a} \bigl[ q^{-I \otimes h_c} (e_c 
   \otimes q^{\frac{1}{2} h_c}) + q^{-h_c \otimes I}(q^{-\frac{1}{2}h_c} 
   \otimes e_c) \bigr] \notag \\
&= q^{\underset{a}{\sum} h_a \otimes h^a} (e_c \otimes q^{-\frac{1}{2} h_c} + 
   q^{-\frac{3}{2} h_c} \otimes e_c).
\end{align*}

\noindent Using this, we see

\begin{align}
\Delta^T(e_c)R &= \Delta^T(e_c) q^{\underset{a}{\sum} h_{a} \otimes h^{a}} 
   \Bigl[ I \otimes I + (q - q^{-1}) \sum_{\varepsilon_{a} <\varepsilon_{b}} 
   (-1)^{[b]} E^a_b \otimes \hat{\sigma}_{ba} \Bigr] \notag \\
&= q^{\underset{a}{\sum} h_a \otimes h^a} (e_c \otimes q^{-\frac{1}{2} h_c} + 
   q^{-\frac{3}{2} h_c} \otimes e_c) \notag \\
&\hspace{4cm} \times \Bigl[ I \otimes I + (q - q^{-1}) 
   \sum_{\varepsilon_{a}< \varepsilon_{b}} (-1)^{[b]} E^a_b \otimes 
   \hat{\sigma}_{ba} \Bigr] \notag \\
&= q^{\underset{a}{\sum} h_a \otimes h^a} \biggl\{ e_c \otimes q^{-\frac{1}{2} 
   h_c} + q^{-\frac{3}{2} h_c} \otimes e_c \notag \\
&\quad + (q - q^{-1}) \sum_{\varepsilon_{a} < \varepsilon_{b}} 
   (-1)^{[b]} \Bigl[ e_c E^a_b \otimes q^{-\frac{1}{2} h_c} \hat{\sigma}_{ba} 
   \notag \\
& \hspace{50mm} + (-1)^{([a] + [b])[c]} q^{-\frac{3}{2} (\alpha_c, 
   \varepsilon_a)} E^a_b \otimes e_c \hat{\sigma}_{ba}  \Bigr] \biggr\}.
   \label{delR} 
\end{align}

\noindent Also,

\begin{align}
R \Delta (e_c) = q^{\underset{a}{\sum} h_a \otimes h^a} \biggl\{ q^{\frac{1}{2}
   h_c} \otimes e_c + e_c \otimes q^{-\frac{1}{2} h_c} \quad& \notag \\
+ (q - q^{-1}) \sum_{\varepsilon_{a} < \varepsilon_{b}} (-1)^{[b]} &\Bigl[ 
   q^{\frac{1}{2} (\alpha_c, \varepsilon_b)} E^a_b \otimes \hat{\sigma}_{ba} 
   e_c \notag \\
&\quad+ (-1)^{([a]+[b])[c]}  E^a_b e_c \otimes \hat{\sigma}_{ba} 
   q^{-\frac{1}{2} h_c} \Bigr] \biggr\} \label{Rdel}.
\end{align}


\noindent  Hence to apply the intertwining property we simply equate (\ref{delR}) and (\ref{Rdel}).   First note that $R$ is weightless, so $\hat{\sigma}_{ba}$ has weight $\varepsilon_b - \varepsilon_{a}$, and thus

\begin{equation*}
q^{-\frac{1}{2} h_c} \hat{\sigma}_{ba} = q^{-\frac{1}{2} (\alpha_c, 
  \varepsilon_b - \varepsilon_a)} \hat{\sigma}_{ba} q^{-\frac{1}{2} h_c}.
\end{equation*}

\noindent Then, equating those terms with zero weight in the first element of the tensor product, we obtain

\begin{multline}
(q^{\frac{1}{2} h_c} - q^{-\frac{3}{2} h_c}) \otimes e_c \\
  =  (q-q^{-1}) \sum_{\varepsilon_{b} - \varepsilon_{a} = \alpha_c} (-1)^{[b]}
    \bigl( q^{-\frac{1}{2} (\alpha_c, \alpha_c)} e_c E^a_b - (-1)^{[c]} 
    E^a_b e_c \bigr) \otimes \hat{\sigma}_{ba} q^{-\frac{1}{2} h_c}.\label{eq1}
\end{multline}

\noindent Comparing the remaining terms, we also find

\begin{multline} \label{eq2}
\underset{\varepsilon_{b} - \varepsilon_{a} \neq \alpha_c} {\sum_{\varepsilon_
  {a} < \varepsilon_{b}}} (-1)^{[b]} \bigl( q^{-\frac{1}{2}(\alpha_c, 
  \varepsilon_{b} - \varepsilon_{a})} e_c E^a_b - (-1)^{([a] + [b])([c])} E^a_b
  e_c \bigr) \otimes  \hat{\sigma}_{ba} q^{-\frac{1}{2} h_c} \\
= \sum_{\varepsilon_{a} < \varepsilon_{b}} (-1)^{[b]} E^a_b \otimes 
   \bigl( q^{\frac{1}{2} (\alpha_c, \varepsilon_b)} \hat{\sigma}_{ba} e_c - 
   (-1)^{([a]+[b])[c]} q^{-\frac{3}{2} (\alpha_c, \varepsilon_a)} e_c 
   \hat{\sigma}_{ba} \bigr).
\end{multline}



\

\noindent  From the first of these equations we can deduce certain fundamental values of $\hat{\sigma}_{ba}$; from the second, relations involving all the $\hat{\sigma}_{ba}$.  Before doing so, however, it is convenient to define a new set, $\overline{\Phi}^+$.

\begin{definition}
The extended system of positive roots, $\overline{\Phi}^+$, is defined by 

\begin{equation*}
\overline{\Phi}^+ \equiv \{\varepsilon_b - \varepsilon_a| \varepsilon_b > \varepsilon_a\} = \Phi^+ \cup \{2 \varepsilon_i | 1 \leq i \leq l\}
\end{equation*}

\noindent where $\Phi^+$ is the usual system of positive roots.
\end{definition}

\

\noindent  Now consider equation (\ref{eq2}).  In the case when $\varepsilon_b - \varepsilon_a + \alpha_c \notin \overline{\Phi}^+$, by collecting the terms of weight $\varepsilon_b - \varepsilon_a + \alpha_c$ in the second half of the tensor product we find:

\begin{equation} \label{*}
q^{\frac{1}{2}(\alpha_c, \varepsilon_{b})} \hat{\sigma}_{ba} e_c - 
    (-1)^{([a]+[b])[c]} q^{-\frac{3}{2} (\alpha_c, \varepsilon_a)} e_c 
    \hat{\sigma}_{ba} = 0.
\end{equation}

\noindent  Similarly, when $\varepsilon_b > \varepsilon_a$ and $\varepsilon_b - \varepsilon_a + \alpha_c = \varepsilon_{b'} - \varepsilon_{a'} \in \overline{\Phi}^+$ we find:

\begin{multline*}
\underset{\varepsilon_{b} - \varepsilon_{a} + \alpha_c = \varepsilon_{b'} - 
  \varepsilon_{a'}} {\sum_{\varepsilon_{a'} < \varepsilon_{b'}}} (-1)^{[b']} 
  \bigl(q^{-\frac{1}{2}(\alpha_c, \varepsilon_{b'}-\varepsilon_{a'})} e_c 
  E^{a'}_{b'} - (-1)^{([a']+[b'])[c]} E^{a'}_{b'} e_c \bigr) \otimes 
  \hat{\sigma}_{b'\!a'} q^{-\frac{1}{2} h_c} \\
= (-1)^{[b]} E^a_b \otimes \bigl( q^{\frac{1}{2} (\alpha_c, \varepsilon_b)} 
  \hat{\sigma}_{ba} e_c - (-1)^{([a]+[b])[c])} q^{-\frac{3}{2}
  (\alpha_c, \varepsilon_a)} e_c \hat{\sigma}_{ba} \bigr).
\end{multline*}

\noindent However $e_c E^{a'}_{b'}$ and $E^a_{b}$ are linearly independent unless $b=b'$, as are $E^{a'}_{b'} e_c$ and $E^{a}_{b}$ for $a \neq a'$, and thus this equation reduces to

\begin{align*}
\underset{\varepsilon_{a'} = \varepsilon_{a} - \alpha_c}{\sum_{\varepsilon_{a'}
  < \varepsilon_{b}}} (-1)^{[b]} q^{-\frac{1}{2} (\alpha_c, \varepsilon_b - 
  \varepsilon_{a'})} e_c E^{a'}_b \otimes &\hat{\sigma}_{ba'} q^{-\frac{1}{2} 
  h_c} \notag \\
-  \underset{\varepsilon_{b'} = \varepsilon_b + \alpha_c} {\sum_{\varepsilon_
  {b'} > \varepsilon_a}} (-1)^{[b'] +([a]+[b'])[c]} &E^a_{b'}e_c \otimes \hat
  {\sigma}_{b'a} q^{-\frac{1}{2} h_c} \notag \\
=  (-1)^{[b]} E^a_b \otimes &\bigl( q^{\frac{1}{2} (\alpha_c,\varepsilon_b)}
   \hat{\sigma}_{ba} e_c - (-1)^{([a]+[b])[c]} q^{-\frac{3}{2} 
   (\alpha_c,\varepsilon_a)} e_c \hat{\sigma}_{ba} \bigr).
\end{align*}

\noindent This can also be written as

\begin{align*}
q^{-\frac{1}{2} (\alpha_c, \varepsilon_b - \varepsilon_{a} + \alpha_c)} 
   e_c &E^{a'}_{b} \otimes \hat{\sigma}_{ba'} q^{-\frac{1}{2} h_c} 
   \Big\vert_{\varepsilon_{a'} = \varepsilon_a - \alpha_c} 
-  (-1)^{([a]+[b])[c]} E^{a}_{b'}e_c \otimes \hat{\sigma}_{b'a} q^{-\frac{1}{2}
   h_c} \Big\vert_{\varepsilon_{b'} = \varepsilon_b + \alpha_c} \notag \\
=  &\; E^{a}_{b} \otimes \bigl( q^{\frac{1}{2} (\alpha_c,\varepsilon_b)} 
   \hat{\sigma}_{ba} e_c - (-1)^{([a]+[b])[c]} q^{-\frac{3}{2} 
   (\alpha_c,\varepsilon_a)} e_c \hat{\sigma}_{ba} \bigr), \qquad 
   \varepsilon_b > \varepsilon_a.
\end{align*}

\noindent This equation then implies

\begin{align*}
q^{-\frac{1}{2} (\alpha_c, \varepsilon_b - \varepsilon_{a} + \alpha_c)} 
   &\langle a|e_c|a' \rangle \hat{\sigma}_{ba'} q^{-\frac{1}{2} h_c} 
-  (-1)^{([a]+[b])[c]} \langle b'|e_c|b \rangle \hat{\sigma}_{b'\!a} 
   q^{-\frac{1}{2} h_c} \notag \\
&=  q^{\frac{1}{2} (\alpha_c,\varepsilon_b)} \hat{\sigma}_{ba} e_c - (-1)^
  {([a]+[b])[c]} q^{-\frac{3}{2} (\alpha_c,\varepsilon_a)}e_c 
  \hat{\sigma}_{ba}, \qquad \varepsilon_b > \varepsilon_a.
\end{align*}

\noindent A more useful form of these relations is:

\begin{align}
q^{-\frac{1}{2} (\alpha_c,  \alpha_c - \varepsilon_{a})} \langle a|&e_c|a' 
  \rangle \hat{\sigma}_{ba'}-(-1)^{([a]+[b])[c]} q^{\frac{1}{2} (\alpha_c,
  \varepsilon_b)} \langle b'|e_c|b \rangle \hat{\sigma}_{b'a} \notag \\
&= q^{(\alpha_c,\varepsilon_b)}\hat{\sigma}_{ba} e_c q^{\frac{1}{2}
   h_c} - (-1)^ {([a]+[b])[c]} q^{-\frac{3}{2} (\alpha_c,\varepsilon_a) + 
   \frac{1}{2} (\alpha_c,\varepsilon_b)}  e_c \hat{\sigma}_{ba} q^{\frac{1}{2} 
   h_c} \notag \\
&= q^{(\alpha_c,\varepsilon_b)} \hat{\sigma}_{ba}e_c q^{\frac{1}{2} h_c} - 
   (-1)^ {([a]+[b])[c]} q^{-(\alpha_c,\varepsilon_a)} e_c q^{\frac{1}{2} h_c} 
   \hat{\sigma}_{ba} \label{**}
\end{align}

\noindent for $\varepsilon_b > \varepsilon_a$.  All the necessary information is contained within these relations and equation \eqref{eq1}.  To construct the Lax operator $R = (\pi \otimes 1) \mathcal{R}$ first we use equation \eqref{eq1} to find the solutions for $\hat{\sigma}_{ba}$ associated with the simple roots $\alpha_c$.  Then we apply the recursion relations arising from \eqref{**} to find the remaining values of $\hat{\sigma}_{ba}$.


\section{Fundamental Values} \label{iv}
    \noindent In this section we solve equation \eqref{eq1}, rewritten below, to find the fundamental values of $\hat{\sigma}_{ba}$, namely those for which $\varepsilon_b - \varepsilon_a$ is a simple root.

\begin{multline*} 
(q^{\frac{1}{2} h_c} - q^{-\frac{3}{2} h_c}) \otimes e_c \\
=  (q-q^{-1}) \sum_{\varepsilon_{b} - \varepsilon_{a} = \alpha_c} (-1)^{[b]} 
   \bigl( q^{-\frac{1}{2} (\alpha_c, \alpha_c)} e_c E^a_b - (-1)^{[c]} E^a_b 
   e_c \bigr) \otimes \hat{\sigma}_{ba} q^{-\frac{1}{2} h_c}. \tag{\ref{eq1}}
\end{multline*}

\

\noindent  To solve this we must consider the various simple roots individually.

\

\noindent \underline{Solution for $\alpha_i = \varepsilon_i - \varepsilon_{i+1},\; 1 \leq i < l$}

\noindent In the vector representation $e_i = E^i_{i+1} - E^{\overline{i+1}}_{\overline{i}}$ and $h_i = E^i_i - E^{i+1}_{i+1} + E^{\overline{i+1}}_{\overline{i+1}} - E^{\overline{i}}_{\overline{i}}$.

\noindent Hence the left-hand side of (\ref{eq1}) becomes:

\begin{align*}
LHS &= (q^{\frac{1}{2} h_i} - q^{-\frac{3}{2} h_i}) \otimes e_i \notag \\
    &= \bigl\{ (q^{\frac{1}{2}} - q^{-\frac{3}{2}})(E^i_i + E^{\overline{i+1}}
    _{\overline{i+1}}) + (q^{-\frac{1}{2}} - q^{\frac{3}{2}})(E^{i+1}_{i+1} 
    + E^{\overline{i}}_{\overline{i}}) \bigr\} \otimes e_i \notag \\
    &= (q-q^{-1}) \bigl\{ q^{-\frac{1}{2}} (E^i_i + E^{\overline{i+1}}_
    {\overline{i+1}}) - q^{\frac{1}{2}} (E^{i+1}_{i+1} + E^{\overline{i}}
    _{\overline{i}}) \bigr\} \otimes e_i,
\end{align*}

\noindent whereas the right-hand side is:

\begin{align*}
RHS &= (q-q^{-1}) \sum_{\varepsilon_{b} - \varepsilon_{a} = \alpha_i} \bigl( 
   q^{-1} e_i E^a_b - E^a_b e_i \bigr) \otimes \hat{\sigma}_{ba} 
   q^{-\frac{1}{2}  h_i} \notag \\
&= (q-q^{-1}) \bigl\{ (q^{-1} E^i_i - E^{i+1}_{i+1}) \otimes \hat{\sigma}_{i\,
   i+1} q^{-\frac{1}{2} h_i} - (q^{-1} E^{\overline{i+1}}_{\overline{i+1}} - 
   E^{\overline{i}}_{\overline{i}}) \otimes \hat{\sigma}_{\overline{i+1}\, 
   \overline{i}} q^{-\frac{1}{2} h_i} \bigr\}.
\end{align*}

\noindent Equating these gives

\begin{equation*}
\boxed{\hat{\sigma}_{i\, i+1} = - \hat{\sigma}_{\overline{i+1}\, \overline{i}} 
 = q^{\frac{1}{2}} e_i q^{\frac{1}{2} h_i}, \quad 1 \leq i < l.}
\end{equation*}


\noindent \underline{Solution for $\alpha_l = \varepsilon_{l-1} + \varepsilon_{l},\; m=2l$}

\noindent Here $e_l = E_{\overline{l}}^{l-1} - E_{\overline{l-1}}^{l}$ and $h_l = E^l_l - E^{\overline{l}}_{\overline{l}} + E^{l-1}_{l-1} - E^{\overline{l-1}}_{\overline{l-1}}$.  Substituting these into equation \eqref{eq1} gives:

\begin{align*}
LHS &= \bigl\{ (q^{\frac{1}{2}} - q^{-\frac{3}{2}}) (E^l_l + E^{l-1}_{l-1}) + 
   (q^{-\frac{1}{2}} - q^{\frac{3}{2}}) (E^{\overline{l}}_{\overline{l}} +  
   E^{\overline{l-1}}_{\overline{l-1}}) \bigr\} \otimes e_l \notag \\
&= (q - q^{-1}) \bigl\{ q^{-\frac{1}{2}} (E^l_l + E^{l-1}_{l-1}) - 
   q^{\frac{1}{2}} (E^{\overline{l}}_{\overline{l}} +  E^{\overline{l-1}}_
   {\overline{l-1}}) \bigr\} \otimes e_l 
\end{align*}

\noindent and 

\begin{equation*}
RHS = (q - q^{-1})  \bigl\{ (q^{-1} E^{l-1}_{l-1} - E^{\overline{l}}_
  {\overline{l}}) \otimes \hat{\sigma}_{l-1 \, \overline{l}} q^{-\frac{1}{2} 
  h_l} - (q^{-1} E^l_l - E^{\overline{l-1}}_{\overline{l-1}}) \otimes 
  \hat{\sigma}_{l \, \overline{l-1}}q^{-\frac{1}{2} h_l} \bigr\} .
\end{equation*}

\noindent Thus

\begin{equation*}
\boxed{\hat{\sigma}_{l-1\, \overline{l}} = - \hat{\sigma}_{l\, \overline{l-1}} 
 = q^{\frac{1}{2}} e_l q^{\frac{1}{2} h_l}, \quad m=2l.}
\end{equation*}


\noindent \underline{Solution for $\alpha_l = \varepsilon_l,\; m=2l+1$}

\noindent  In this case $e_l = E^l_{l+1} - E^{l+1}_{\overline{l}}$ whereas 
$h_l = E^l_l - E^{\overline{l}}_{\overline{l}}$, so we obtain:

\begin{align*}
LHS &= \bigl\{ (q^{\frac{1}{2}} - q^{-\frac{3}{2}}) E^l_l + (q^{-\frac{1}{2}} -
  q^{\frac{3}{2}}) E^{\overline{l}}_{\overline{l}} \bigr\} \otimes e_l \notag\\
&=  (q-q^{-1}) ( q^{-\frac{1}{2}} E^l_l - q^{\frac{1}{2}} 
    E^{\overline{l}}_{\overline{l}} ) \otimes e_l, \\
\notag \\
RHS &= (q-q^{-1}) \bigl\{ (q^{-\frac{1}{2}} E^l_l - E^{l+1}_{l+1}) \otimes \hat{\sigma}_{l\, l+1} q^{-\frac{1}{2} h_l} - (q^{-\frac{1}{2}} E^{l+1}_{l+1} - E^{\overline{l}}_{\overline{l}})  \otimes \hat{\sigma}_{l+1 \, \overline{l}} q^{-\frac{1}{2} h_l} \bigr\}.
\end{align*}

\noindent Together these imply

\begin{equation*}
\boxed{\hat{\sigma}_{l\, l+1} = - q^{-\frac{1}{2}} \hat{\sigma}_{l+1\, 
 \overline{l}} = e_l q^{\frac{1}{2} h_l}, \quad m=2l+1.}
\end{equation*}


\noindent \underline{Solution for $\alpha_{\mu} = \delta_{\mu} - \delta_{\mu +1},\; 1 \leq \mu < k$}

\noindent Here $e_\mu = E^\mu_{\mu + 1} + E^{\overline{\mu + 1}}_{\overline{\mu}}$ and $h_\mu = E^{\mu + 1}_{\mu + 1} + E^{\overline{\mu}}_{\overline{\mu}} - E^\mu_\mu - E^{\overline{\mu + 1}}_{\overline{\mu +1}}$, giving:

\begin{align*}
LHS &= \bigl\{ (q^{\frac{1}{2}} - q^{-\frac{3}{2}}) (E^{\mu + 1}_{\mu + 1} + 
   E^{\overline{\mu}}_{\overline{\mu}}) + (q^{-\frac{1}{2}} - q^{\frac{3}{2}}) 
   (E^\mu_\mu + E^{\overline{\mu + 1}}_{\overline{\mu +1}}) \bigr\} \otimes 
   e_\mu \notag \\
&= (q - q^{-1}) \bigl\{ q^{-\frac{1}{2}} (E^{\mu + 1}_{\mu + 1} + 
   E^{\overline{\mu}}_{\overline{\mu}}) - q^{\frac{1}{2}} (E^\mu_\mu + 
   E^{\overline{\mu + 1}}_{\overline{\mu +1}}) \bigr\} \otimes e_\mu,  \\
\notag \\
RHS &= -(q - q^{-1}) \bigl\{ (q E^\mu_\mu - E^{\mu +1}_{\mu +1}) \otimes \hat{\sigma}_{\mu\, \mu+1} q^{-\frac{1}{2} h_\mu} + (q E^{\overline{\mu+1}}_{\overline{\mu+1}} - E^{\overline{\mu}}_{\overline{\mu}}) \otimes \hat{\sigma}_{\overline{\mu+1}\, \overline{\mu}} q^{-\frac{1}{2} h_\mu} \bigr\},
\end{align*}

\noindent  and hence

\begin{equation*}
\boxed{\hat{\sigma}_{\mu \, \mu +1} = \hat{\sigma}_{\overline{\mu+1} \, 
  \overline{\mu}} = q^{-\frac{1}{2}} e_\mu q^{\frac{1}{2} h_\mu}, \quad 1 
  \leq \mu <k.}
\end{equation*}


\noindent \underline{Solution for $\alpha_s = \delta_k - \varepsilon_1$}

\noindent In this case $e_s = E^{\mu = k}_{i=1} + (-1)^k E^{\overline{i}= \overline{1}}_{\overline{\mu}= \overline{k}}$ and $h_s = E^{\overline{i}= \overline{1}}_{\overline{i}= \overline{1}} - E^{i=1}_{i=1} + E^{\overline{\mu} = \overline{k}}_{\overline{\mu} = \overline{k}} - E^{\mu =k}_{\mu =k}$.  Substituting these into equation \eqref{eq1} produces:

\begin{align*}
LHS &= \bigl\{ (q^{\frac{1}{2}} - q^{-\frac{3}{2}}) (E^{\overline{i}=\overline
  {1}}_{\overline{i}= \overline{1}} + E^{\overline{\mu} = \overline{k}}
  _{\overline{\mu} = \overline{k}}) + (q^{-\frac{1}{2}} - q^{\frac{3}{2}}) 
  (E^{i=1}_{i=1} +  E^{\mu =k}_{\mu =k})\bigr\} \otimes e_s \notag \\
&=  (q- q^{-1}) \bigl\{ q^{-\frac{1}{2}} (E^{\overline{i}= \overline{1}}_
   {\overline{i}= \overline{1}} + E^{\overline{\mu} = \overline{k}}_{\overline
   {\mu} = \overline{k}}) - q^{\frac{1}{2}} (E^{i=1}_{i=1} +  E^{\mu =k}_{\mu =
   k}) \bigr\} \otimes e_s, \\
\notag \\
RHS &= (q- q^{-1}) \bigl\{ -(E^{\mu =k}_{\mu =k} + E^{i=1}_{i=1}) \otimes 
  \hat{\sigma}_{\mu =k \, i=1} q^{-\frac{1}{2} h_s} \notag \\
& \hspace{5cm} + (-1)^k (E^{\overline{i}= \overline{1}}_{\overline{i}= 
  \overline{1}} + E^{\overline{\mu} = \overline{k}}_{\overline{\mu} = \overline
  {k}}) \otimes \hat{\sigma}_{\overline{i}= \overline{1} \, \overline{\mu} = 
  \overline{k}} q^{-\frac{1}{2} h_s} \bigr\},
\end{align*}

\noindent and thus

\begin{equation*}
\boxed{\hat{\sigma}_{\mu=k \, i=1} = (-1)^k q\, \hat{\sigma}_{\overline{i} = 
  \overline{1} \, \overline{\mu} = \overline{k}} = q^{\frac{1}{2}} e_s 
  q^{\frac{1}{2} h_s}.}
\end{equation*}

\noindent These values for $\hat{\sigma}_{ba}$ form the basis for finding $R$, as from these all the others can be explicitly determined in any given representation.

\section{Constructing the Non-Simple Values} 
  \label{alphai}
  
  \noindent Now we develop the recurrence relations required to calculate the remaining values of $\hat{\sigma}_{ba}$.  Recall that for $\varepsilon_b > \varepsilon_a$,

\begin{multline} \label{**2}
q^{-\frac{1}{2} (\alpha_c,  \alpha_c - \varepsilon_{a})} \langle a|e_c|a' 
  \rangle \hat{\sigma}_{ba'}-(-1)^{([a]+[b])[c]} q^{\frac{1}{2} (\alpha_c, 
   \varepsilon_b)} \langle b'|e_c|b \rangle \hat{\sigma}_{b'a} \\
= q^{(\alpha_c,\varepsilon_b)} \hat{\sigma}_{ba} e_c q^{\frac{1}{2} h_c} - 
   (-1)^ {([a]+[b])[c]} q^{-(\alpha_c,\varepsilon_a)} e_c q^{\frac{1}{2} h_c} 
   \hat{\sigma}_{ba}.
\end{multline}

\noindent  To extract the recurrence relations to be applied to the fundamental values of $\hat{\sigma}_{ba}$, we must again consider the simple roots individually.  We begin with the case $\alpha_i = \varepsilon_i - \varepsilon_{i+1}$, so \hbox{$e_i = \sigma^i_{i+1} \equiv E^i_{i+1} - E^{\overline{i+1}}_{\overline{i}}$}.  Now

\begin{equation*}
\langle a|e_i = \delta_{ai} \langle i+1| - \delta_{a\, \overline{i+1}} 
  \langle \overline{i}|,\qquad
e_i |b \rangle = \delta_{b\, i+1} |i \rangle - \delta_{b \overline{i}} 
  |\overline{i+1} \rangle.
\end{equation*}

\noindent We then apply that to equation \eqref{**2} to obtain:

\begin{align*}
q^{-\frac{1}{2}(\alpha_i, \alpha_i)} \bigl\{ \delta_{a i} 
   q^{\frac{1}{2}(\alpha_i, \varepsilon_i)} \hat{\sigma}_{b\, i+1}- \delta_{a\,
   \overline{i+1}} &q^{-\frac{1}{2}(\alpha_i, \varepsilon_{i+1})} \hat{\sigma}_
   {b\, \overline{i}} \bigr\} \notag \\ 
- \bigl\{ \delta_{b\, i+1} q^{\frac{1}{2} 
   (\alpha_i, \varepsilon_{i+1})}&\hat{\sigma}_{ia} - \delta_{b\, \overline{i}}
   q^{-\frac{1}{2} (\alpha_i, \varepsilon_i)} \hat{\sigma}_{\overline{i+1}\,a} 
   \bigr\} \notag \\
= &q^{(\alpha_i,\varepsilon_b)} \hat{\sigma}_{ba}e_i q^{\frac{1}{2}
   h_i} - q^{-(\alpha_i,\varepsilon_a) } 
   e_i q^{\frac{1}{2} h_i} \hat{\sigma}_{ba}, \quad \varepsilon_b >
   \varepsilon_a.
\end{align*}

\noindent This simplifies to

\begin{multline*}
q^{-\frac{1}{2}} \bigl\{ \delta_{ai} \hat{\sigma}_{b\, i+1}- \delta_{a\, 
  \overline{i+1}} \hat{\sigma}_{b\, \overline{i}} - \delta_{b\, i+1} 
  \hat{\sigma}_{ia} + \delta_{b\,\overline{i}} \hat{\sigma}_{\overline{i+1}\,a}
  \bigr\} \\
= q^{(\alpha_i,\varepsilon_b)} \hat{\sigma}_{ba}e_i q^{\frac{1}{2} h_i} - 
  q^{-(\alpha_i,\varepsilon_a)} e_i q^{\frac{1}{2} h_i} \hat{\sigma}_{ba}, 
  \quad \varepsilon_b > \varepsilon_a,
\end{multline*}

\noindent which, recalling that $\hat{\sigma}_{i\, i+1} = - \hat{\sigma}_{\overline{i+1}\, \overline{i}} = q^{\frac{1}{2}} e_i q^{\frac{1}{2} h_i},\; 1 \leq i < l$, reduces to

\begin{align*}
\delta_{ai} \hat{\sigma}_{b\, i+1}- \delta_{a\, 
  \overline{i+1}} \hat{\sigma}_{b\, \overline{i}} - \delta_{b\, i+1} 
  \hat{\sigma}_{ia} + \delta_{b\,\overline{i}} \hat{\sigma}_{\overline{i+1}\,a}
&=q^{(\alpha_i,\varepsilon_b)} \hat{\sigma}_{ba} \hat{\sigma}_{i\, i+1} - 
  q^{-(\alpha_i,\varepsilon_a)} \hat{\sigma}_{i\, i+1} \hat{\sigma}_{ba} \notag
  \\
&=q^{-(\alpha_i,\varepsilon_a)} \hat{\sigma}_{\overline{i+1}\, \overline{i}}
  \hat{\sigma}_{ba} - q^{(\alpha_i,\varepsilon_b)} \hat{\sigma}_{ba} 
  \hat{\sigma}_{\overline{i+1}\, \overline{i}}.
\end{align*}

\noindent From this we can deduce the following relations for $1 \leq i <l$:

\begin{alignat}{2}
&\hat{\sigma}_{b\, i+1} = \hat{\sigma}_{b\,i} \hat{\sigma}_{i\, i+1} - q^{-1} 
   \hat{\sigma}_{i\, i+1} \hat{\sigma}_{b\, i}, & \qquad &\varepsilon_b > 
   \varepsilon_i, \notag \\
&\hat{\sigma}_{\overline{i+1}\, a} = \hat{\sigma}_{\overline{i+1}\, 
    \overline{i}} \hat{\sigma}_{\overline{i}\, a} - q^{-1} \hat{\sigma}_
    {\overline{i}\, a} \hat{\sigma}_{\overline{i+1}\, \overline{i}}, & &
    \varepsilon_a < -\varepsilon_i, \notag \\
&\hat{\sigma}_{b\, \overline{i}} = q^{(\alpha_i, \varepsilon_b)} \hat{\sigma}_
   {b\, \overline{i+1}} \hat{\sigma}_{\overline{i+1}\, \overline{i}} - q^{-1} 
   \hat{\sigma}_{\overline{i+1}\,\overline{i}} \hat{\sigma}_{b\,\overline{i+1}}
   , & & \varepsilon_b > - \varepsilon_{i+1},\, b \neq i+1, \notag \\
&\hat{\sigma}_{i\,a} = q^{-(\alpha_i, \varepsilon_a)} \hat{\sigma}_{i\, i+1} 
   \hat{\sigma}_{i+1\,a} - q^{-1} \hat{\sigma}_{i+1\,a} \hat{\sigma}_{i\,i+1}, 
   & & \varepsilon_a < \varepsilon_{i+1}, \; a \neq \overline{i+1},\notag \\
&\hat{\sigma}_{i\, \overline{i+1}} + \hat{\sigma}_{i+1\,\overline{i}} = q^{-1} 
   \bigl[ \hat{\sigma}_{i\,i+1}, \hat{\sigma}_{i+1\,\overline{i+1}} \bigr],&&
   \label{icom} \\
&q^{(\alpha_i, \varepsilon_b)} \hat{\sigma}_{ba} \hat{\sigma}_{i\,i+1} - 
   q^{-(\alpha_i, \varepsilon_a)} \hat{\sigma}_{i\, i+1} \hat{\sigma}_{ba} = 0,
   && \varepsilon_b > \varepsilon_a,\; a \neq i, \overline{i+1}\text{ and }b 
   \neq i+1, \overline{i}. \notag
\end{alignat}

\noindent We then follow the same procedure to find the relations associated with the other simple roots.  A detailed derivation of these relations is included in Appendix A, with a complete list of the relations derived in this manner given in Tables \ref{list}, \ref{even} and \ref{odd} on pages \pageref{list} and \pageref{odd}.

\

\noindent  Although the list of relations is very long, they can be summarised in a compact form.  There are two different types of relations; recursive and $q$-commutative.  The latter can be condensed into:

\vspace{-2mm}
\begin{equation} \label{commutationrelations}
\boxed{q^{(\alpha_c, \varepsilon_b)} \hat{\sigma}_{ba} e_c q^{\frac{1}{2} h_c} 
      - (-1)^{([a]+[b])[c]} q^{-(\alpha_c, \varepsilon_a)} e_c 
      q^{\frac{1}{2}h_c} \hat{\sigma}_{ba} = 0, \quad \varepsilon_b > 
      \varepsilon_a}
\end{equation}
      
\noindent where neither $\varepsilon_a - \alpha_c$ nor $\varepsilon_b + \alpha_c$ equals any $\varepsilon_x$.  Note this is almost the same as equation \eqref{*}, with slightly softer restrictions on $a$ and $b$ the only difference.  

\

\noindent  The recursion relations can, in the case $m=2l+1$, be summarised as:

\vspace{-2mm}
\begin{equation} \label{nice}
\boxed{\hat{\sigma}_{ba} = q^{-(\varepsilon_b,\varepsilon_a)} \hat{\sigma}_{bc}
  \hat{\sigma}_{ca} - q^{-(\varepsilon_c, \varepsilon_c)} (-1)^{([b]+[c])
  ([a]+[c])} \hat{\sigma}_{ca} \hat{\sigma}_{bc}, \quad \varepsilon_b
  > \varepsilon_c > \varepsilon_a}
\end{equation}

\noindent where $ c \neq \overline{b} \text{ or } \overline{a}$.  

\

\noindent In the case $m=2l$, these include all the information except that contained in equation \eqref{icom} for the case $i = l-1$.  Hence when $m=2l$ one extra relation is required, namely:

\vspace{-2mm}
\begin{equation} \label{l-1lbar}
\boxed{\hat{\sigma}_{l-1\, \overline{l}} + \hat{\sigma}_{l\, \overline{l-1}} =
  q^{-1} [ \hat{\sigma}_{l-1\, l}, \hat{\sigma}_{l \overline{l}} ],\quad m=2l.}
\end{equation}

\

\noindent  It is not difficult to see that all these recursion relations can be obtained from those listed on pages \pageref{list} and \pageref{odd}.  To show the reverse is tedious, but straightforward.  The only relations on pages \pageref{list} and \pageref{odd} which are not clearly of this form are those involving commutators.  As an example of how these can be obtained from equation \eqref{nice}, consider equation \eqref{icom}, which arose from considering $\alpha_i = \varepsilon_i - \varepsilon_{i+1}$:

\begin{equation*}
\hat{\sigma}_{i\, \overline{i+1}} + \hat{\sigma}_{i+1\,\overline{i}} = q^{-1} 
   \bigl[ \hat{\sigma}_{i\,i+1}, \hat{\sigma}_{i+1\,\overline{i+1}} \bigr].
\end{equation*}

\noindent Using equation \eqref{nice}, we can say that:

\begin{alignat*}{2}
&\hat{\sigma}_{i\, \overline{i+1}} = \hat{\sigma}_{i\, \overline{i+2}} 
  \hat{\sigma}_{\overline{i+2} \, \overline{i+1}} - q^{-1} \hat{\sigma}_
  {\overline{i+2}\, \overline{i+1}} \hat{\sigma}_{i\, \overline{i+2}},& \quad
  &i < l-1, \\
&\hat{\sigma}_{i+1\, \overline{i}} = \hat{\sigma}_{i+1\, \overline{i+2}}
  \hat{\sigma}_{\overline{i+2} \, \overline{i}} - q^{-1}\hat{\sigma}_{\overline
  {i+2} \, \overline{i}} \hat{\sigma}_{i+1\, \overline{i+2}}, && i < l-1,\\
&\hat{\sigma}_{\overline{i+2}\, \overline{i}} = \hat{\sigma}_{\overline{i+2}\,
  \overline{i+1}} \hat{\sigma}_{\overline{i+1}\, \overline{i}} - q^{-1}
  \hat{\sigma}_{\overline{i+1}\, \overline{i}} \hat{\sigma}_{\overline{i+2}\,
  \overline{i+1}},&& i< l-1, \\
&\hat{\sigma}_{i\, \overline{i+2}} = \hat{\sigma}_{i\, i+1} \hat{\sigma}_{i+1\,
  \overline{i+2}} - q^{-1} \hat{\sigma}_{i+1\, \overline{i+2}} \hat{\sigma}_
  {i\, i+1}, && i< l-1, \\
&\hat{\sigma}_{i+1\, \overline{i+1}} = q \hat{\sigma}_{i+1\, \overline{i+2}} 
  \hat{\sigma}_{\overline{i+2}\, \overline{i+1}} - q^{-1} \hat{\sigma}_
  {\overline{i+2}\, \overline{i+1}} \hat{\sigma}_{i+1\, \overline{i+2}},
  && i< l-1 .
\end{alignat*}

\noindent Combining these, we find that for $i<l-1$

\begin{align*}
\hat{\sigma}_{i\, \overline{i+1}} + \hat{\sigma}_{i+1\, \overline{i}} &=
\hat{\sigma}_{i\, \overline{i+2}} \hat{\sigma}_{\overline{i+2}\,\overline{i+1}}
+\hat{\sigma}_{i+1\,\overline{i+2}} \hat{\sigma}_{\overline{i+2}\,\overline{i}}
- q^{-1} \hat{\sigma}_{\overline{i+2}\, \overline{i+1}} \hat{\sigma}_{i\, 
  \overline{i+2}} - q^{-1} \hat{\sigma}_{\overline{i+2} \, \overline{i}} 
  \hat{\sigma}_{i+1\, \overline{i+2}}  \\
&=\hat{\sigma}_{i\, \overline{i+2}} \hat{\sigma}_{\overline{i+2}\,
  \overline{i+1}} + \hat{\sigma}_{i+1\,\overline{i+2}} \bigl( \hat{\sigma}_
  {\overline{i+2}\,\overline{i+1}} \hat{\sigma}_{\overline{i+1}\, \overline{i}}
  - q^{-1} \hat{\sigma}_{\overline{i+1}\, \overline{i}} \hat{\sigma}_
  {\overline{i+2}\, \overline{i+1}} \bigr) \\
& \qquad - q^{-1} \hat{\sigma}_{\overline{i+2}\, \overline{i+1}} \bigl( \hat{
  \sigma}_{i\, i+1} \hat{\sigma}_{i+1\,\overline{i+2}} - q^{-1} \hat{\sigma}_
  {i+1\,\overline{i+2}} \hat{\sigma}_{i\, i+1} \bigr) - q^{-1} \hat{\sigma}_
  {\overline{i+2}\, \overline{i}} \hat{\sigma}_{i+1\, \overline{i+2}}  \\
&=\bigl( \hat{\sigma}_{i\, \overline{i+2}} - q^{-1} \hat{\sigma}_{i+1\,
  \overline{i+2}} \hat{\sigma}_{\overline{i+1}\, \overline{i}} \bigr) 
  \hat{\sigma}_{\overline{i+2}\,\overline{i+1}} - q^{-1} \bigl(
  \hat{\sigma}_{\overline{i+2}\, \overline{i}} + \hat{\sigma}_{\overline{i+2}\,
  \overline{i+1}} \hat{\sigma}_{i\, i+1} \bigr) \hat{\sigma}_{i+1\, 
  \overline{i+2}}  \\
& \qquad + \hat{\sigma}_{i+1\,\overline{i+2}} \hat{\sigma}_{\overline{i+2}\,
  \overline{i+1}} \hat{\sigma}_{\overline{i+1}\, \overline{i}} + q^{-2}
  \hat{\sigma}_{\overline{i+2}\, \overline{i+1}} \hat{\sigma}_{i+1\,
  \overline{i+2}} \hat{\sigma}_{i\, i+1} \\
&=\bigl( \hat{\sigma}_{i\, \overline{i+2}} + q^{-1} \hat{\sigma}_{i+1\,
  \overline{i+2}} \hat{\sigma}_{i\, i+1} \bigr) \hat{\sigma}_{\overline{i+2}\,
  \overline{i+1}} - q^{-1} \bigl( \hat{\sigma}_{\overline{i+2}\, \overline{i}} 
  - \hat{\sigma}_{\overline{i+2}\, \overline{i+1}} \hat{\sigma}_{\overline{i+1}
  \, \overline{i}} \bigr) \hat{\sigma}_{i+1\, \overline{i+2}} \\
& \qquad - \bigl( \hat{\sigma}_{i+1\,\overline{i+2}} \hat{\sigma}_
  {\overline{i+2}\, \overline{i+1}} - q^{-2} \hat{\sigma}_{\overline{i+2}\, 
  \overline{i+1}} \hat{\sigma}_{i+1\, \overline{i+2}} \bigr) \hat{\sigma}_
  {i\,i+1} \\
&= \hat{\sigma}_{i\, i+1} \hat{\sigma}_{i+1\, \overline{i+2}} \hat{\sigma}_
  {\overline{i+2}\,\overline{i+1}} + q^{-2} \hat{\sigma}_{\overline{i+1}\, 
  \overline{i}} \hat{\sigma}_{\overline{i+2}\,\overline{i+1}} \hat{\sigma}_
  {i+1\, \overline{i+2}} - q^{-1} \hat{\sigma}_{i+1\, \overline{i+1}} \hat
  {\sigma}_{i\, i+1} \\
&= \hat{\sigma}_{i\, i+1} \bigl( \hat{\sigma}_{i+1\, \overline{i+2}} \hat
  {\sigma}_{\overline{i+2}\,\overline{i+1}} - q^{-2} \hat{\sigma}_{\overline
  {i+2}\,\overline{i+1}} \hat{\sigma}_{i+1\, \overline{i+2}} \bigr) - q^{-1} 
  \hat{\sigma}_{i+1\, \overline{i+1}} \hat{\sigma}_{i\, i+1} \\
&= q^{-1} \hat{\sigma}_{i\, i+1} \hat{\sigma}_{i+1\, \overline{i+1}} - q^{-1} 
  \hat{\sigma}_{i+1\, \overline{i+1}} \hat{\sigma}_{i\, i+1} \\
&= q^{-1} [\hat{\sigma}_{i\, i+1}, \hat{\sigma}_{i+1\, \overline{i+1}}] 
\end{align*}

\noindent as required.  Note that this working also holds true for $i=l-1$ in the case $m=2l+1$, as then the conditions on equation \eqref{nice} are all met.  In the case $m=2l$, however, there is no way to meet the conditions to find an expansion of $\hat{\sigma}_{l-1\, \overline{l}}$ or $\hat{\sigma}_{l\, \overline{l-1}}$, which is why the extra relation  \eqref{l-1lbar} must be included.  The technique used above can be applied almost identically to find the commutation relations that arose from the roots $\alpha_\mu$ and $\alpha_s$.  Hence equations \eqref{commutationrelations}, \eqref{nice} and \eqref{l-1lbar} are equivalent to the complete set of $q$-commutation and recursion relations derived in Appendix A.

\

\noindent In the case $m=2l$, we can find an alternative extra relation to equation \eqref{l-1lbar}.  For such $m$, consider the relations involving $\hat{\sigma}_{l\, \overline{l}}$.

\

\noindent Firstly, we have

\begin{alignat*}{2}
&q^{(\alpha_i, \varepsilon_b)} \hat{\sigma}_{ba} \hat{\sigma}_{i\,i+1} - 
   q^{-(\alpha_i, \varepsilon_a)} \hat{\sigma}_{i\, i+1} \hat{\sigma}_{ba} = 0,
  \quad && \varepsilon_a < \varepsilon_b,\; a \neq i, \overline{i+1}
  \text{ and }b \neq \overline{i}, i+1 \notag \\
\Rightarrow  \quad &[\hat{\sigma}_{i\, i+1}, \hat{\sigma}_{l\, \overline{l}}] 
  = 0, && i<l-1
\end{alignat*}

\noindent and

\begin{align*}
&q^{(\alpha_l, \varepsilon_b)} \hat{\sigma}_{ba} \hat{\sigma}_{l-1 \,
  \overline{l}} - q^{-(\alpha_l, \varepsilon_a)} \hat{\sigma}_{l-1 \,
  \overline{l}} \hat{\sigma}_{ba} = 0, \quad \varepsilon_a < \varepsilon_b,\; 
  a \neq l-1,l  \text{ and }b \neq \overline{l-1},\overline{l}  \notag \\
\Rightarrow  \quad &[\hat{\sigma}_{l-1\, \overline{l}}, 
  \hat{\sigma}_{l\, \overline{l}}] = 0.
\end{align*}

\noindent Moreover, we know

\begin{alignat*}{2}
 & [\hat{\sigma}_{l-1\, l}, \hat{\sigma}_{l\, \overline{l}}] =
  q (\hat{\sigma}_{l-1\, \overline{l}} + \hat{\sigma}_{l\, \overline{l-1}})&&
  \notag \\
\Leftrightarrow \quad & [\hat{\sigma}_{l-1\, l}, \hat{\sigma}_{l\, 
  \overline{l}}] = 0
\end{alignat*}

\noindent from the simple generators in Section \ref{iv}. Similarly, $\hat{\sigma}_{l\,\overline{l}}$ can be shown to commute with the remaining simple generators.  Together these imply $\hat{\sigma}_{l\,\overline{l}}$ is an invariant of the system.  It cannot, therefore, have weight $2 \varepsilon_l$, as it would if it were non-zero, so

\begin{equation*} 
\boxed{\hat{\sigma}_{l\, \overline{l}} = 0, \qquad m=2l.}
\end{equation*}

\noindent This equation is a convenient alternative to \eqref{l-1lbar} in the unified form of the relations.

\

\noindent Hence we have found the following result:

\

\begin{lemma}

\noindent There is a unique matrix in $(\text{End }V) \otimes U_q[osp(m|n)]^+$  of the form

\begin{equation*}
R = q^{\underset{a}{\sum} h_{a} \otimes h^{a}} \Bigl[I \otimes I + (q - q^{-1}) \sum_{\varepsilon_{a} < \varepsilon_{b}} (-1)^{[b]} E^a_b \otimes \hat{\sigma}_{ba} \Bigr],
\end{equation*} 

\noindent satisfying \mbox{$R \Delta (e_c) = \Delta^T (e_c) R$}. The fundamental values of $\hat{\sigma}_{ba}$ for that matrix are given by:

\begin{alignat}{2}
&\hat{\sigma}_{i\, i+1} = - \hat{\sigma}_{\overline{i+1}\, \overline{i}} = 
  q^{\frac{1}{2}} e_i q^{\frac{1}{2} h_i},&&1 \leq i < l, \notag \\
&\hat{\sigma}_{l-1\, \overline{l}} = - \hat{\sigma}_{l\, \overline{l-1}} = 
  q^{\frac{1}{2}} e_l q^{\frac{1}{2} h_l}, && m=2l, \notag \\
&\hat{\sigma}_{l\, \overline{l}} = 0,  && m=2l, \notag \\
&\hat{\sigma}_{l\, l+1} = - q^{-\frac{1}{2}} \hat{\sigma}_{l+1\, \overline{l}} 
  = e_l q^{\frac{1}{2} h_l}, && m=2l+1, \notag \\
&\hat{\sigma}_{\mu \, \mu +1} = \hat{\sigma}_{\overline{\mu+1}\,\overline{\mu}}
  = q^{-\frac{1}{2}} e_\mu q^{\frac{1}{2} h_\mu}, &\qquad&1 \leq \mu <k, \notag
  \\
&\hat{\sigma}_{\mu=k \, i=1} = (-1)^k q\, \hat{\sigma}_{i = \overline{1} \, 
  \overline{\mu} = \overline{k}} = q^{\frac{1}{2}} e_s q^{\frac{1}{2} h_s};
\label{ivalues}
\end{alignat}

\newpage

\noindent and the remaining values can be calculated using

\

\noindent (i) the $q$-commutation relations

\begin{equation} \label{qcom}
q^{(\alpha_c, \varepsilon_b)} \hat{\sigma}_{ba} e_c q^{\frac{1}{2} h_c} - 
  (-1)^{([a]+[b])[c]} q^{-(\alpha_c, \varepsilon_a)} e_c q^{\frac{1}{2}h_c} 
  \hat{\sigma}_{ba} = 0, \quad \varepsilon_b > \varepsilon_a,
\end{equation}
      
\noindent where neither $\varepsilon_a - \alpha_c$ nor $\varepsilon_b + \alpha_c$ equals any $\varepsilon_x$; and

\

\noindent (ii) the induction relations

\begin{equation} \label{indrel}
\hat{\sigma}_{ba} = q^{-(\varepsilon_b,\varepsilon_a)} \hat{\sigma}_{bc}
  \hat{\sigma}_{ca} - q^{-(\varepsilon_c, \varepsilon_c)} (-1)^{([b]+[c])
  ([a]+[c])} \hat{\sigma}_{ca} \hat{\sigma}_{bc}, \quad \varepsilon_b
  > \varepsilon_c > \varepsilon_a,
\end{equation}

\noindent where $ c \neq \overline{b} \text{ or } \overline{a}$.  
\end{lemma}

\

\

\noindent  This matrix can also be written in a slightly different form.  As we are working in the $(\pi \otimes \text{id})$ representation, we have

\begin{align*}
q^{\underset{a}{\sum} h_a \otimes h^a} &= (\sum_a E^a_a \otimes I) 
   q^{\underset{b}{\sum} h_b \otimes h^b} \notag \\
&= (\sum_a E^a_a \otimes I) q^{\underset{b}{\sum}(\varepsilon_a, \varepsilon_b)
   I \otimes h^b} \notag \\
&= \sum_a E^a_a \otimes q^{h_{\varepsilon_a}}.
\end{align*}

\noindent  Hence an alternative way of expressing $R$ is

\begin{equation*}
R = \sum_a E^a_a \otimes q^{h_{\varepsilon_a}} + (q-q^{-1}) \sum_{\varepsilon_a
  < \varepsilon_b} (-1)^{[b]} E^a_b \otimes q^{h_{\varepsilon_a}} \hat{\sigma}_{ba},
\end{equation*}

\noindent with the $\hat{\sigma}_{ba}$ as given before.

\chapter{A Closer Look at the Lax operator} \label{close}
\label{closerlook}
        
\noindent We have found a set of fundamental values and relations which uniquely define the unknowns $\hat{\sigma}_{ba}$.  Theoretically the resultant matrix $R$ must be a Lax operator, as we know there is one of the given form.  It seems advisable, however, to check this by verifying that $R$ satisfies the remaining $R$-matrix properties.  These are 

\begin{equation} \label{1tensordel}
(\text{id} \otimes \Delta) R = R_{13} R_{12}
\end{equation}

\noindent and the intertwining property for the remaining generators,

\begin{equation*}
R \Delta(a) = \Delta^T(a) R, \qquad \forall a \in U_q[osp(m|n)].
\end{equation*}

\noindent  In this chapter we confirm that $R$ satisfies both these properties.  We also calculate the opposite Lax operator $R^T$, and briefly examine whether the defining relations for the $\hat{\sigma}_{ba}$ incorporate the $q$-Serre relations for $U_q[osp(m|n)]$.

\section{Calculating the Coproduct}
    \noindent We begin by considering the first of these defining properties, equation \eqref{1tensordel}.  In order to evaluate $(\text{id} \otimes \Delta) R$, however, we need to know $\Delta (\hat{\sigma}_{ba})$.  Now at the end of the previous chapter we showed that

\begin{equation*}
R= \sum_a E^a_a \otimes q^{h_{\varepsilon_a}} \; + \; (q-q^{-1}) \sum_{\varepsilon_{b} > \varepsilon_{a}} (-1)^{[b]} E^a_b \otimes q^{h_{\varepsilon_a}} \hat{\sigma}_{ba}.
\end{equation*}

\noindent Using this form for $R$, we find

\begin{align*}
R_{13} R_{12} &= \Bigl( \sum_a E^a_a \otimes I \otimes q^{h_{\varepsilon_a}} \;
    + \; (q-q^{-1}) \sum_{\varepsilon_{b} > \varepsilon_{a}} (-1)^{[b]} E^a_b 
    \otimes I \otimes q^{h_{\varepsilon_a}} \hat{\sigma}_{ba} \Bigr) \notag \\
& \hspace{12mm} \Bigl( \sum_c E^c_c \otimes q^{h_{\varepsilon_c}} \otimes I \; 
    + \; (q-q^{-1}) \sum_{\varepsilon_{d} > \varepsilon_{c}} (-1)^{[d]} E^c_d 
    \otimes q^{h_{\varepsilon_c}} \hat{\sigma}_{dc} \otimes I \Bigr) \notag \\
&= \sum_a \sum_c E^a_a E^c_c \otimes q^{h_{\varepsilon_c}} \otimes 
    q^{h_{\varepsilon_a}} \notag \\ 
& \qquad + (q-q^{-1}) \sum_a \sum_{\varepsilon_{d} >
    \varepsilon_{c}} (-1)^{[d]} E^a_a E^c_d \otimes q^{h_{\varepsilon_c}} 
    \hat{\sigma}_{dc} \otimes q^{h_{\varepsilon_a}} \notag \\
& \qquad + (q-q^{-1}) \sum_c \sum_{\varepsilon_{b} > \varepsilon_{a}} 
    (-1)^{[b]} E^a_b E^c_c \otimes q^{h_{\varepsilon_c}} \otimes 
    q^{h_{\varepsilon_a}} \hat{\sigma}_{ba}\notag \\
& \qquad + (q-q^{-1})^2 \sum_{\varepsilon_{b} > \varepsilon_{a}} 
    \sum_{\varepsilon_{d} > \varepsilon_{c}} (-1)^{[b]+[d]} E^a_b E^c_d \otimes
    q^{h_{\varepsilon_c}} \hat{\sigma}_{dc} \otimes q^{h_{\varepsilon_a}} 
    \hat{\sigma}_{ba} \notag \\
&= \sum_a E^a_a \otimes q^{h_{\varepsilon_a}} \otimes q^{h_{\varepsilon_a}} \;
    + \; (q-q^{-1}) \sum_{\varepsilon_{d} >\varepsilon_{c}} (-1)^{[d]} E^c_d 
    \otimes q^{h_{\varepsilon_c}} \hat{\sigma}_{dc} \otimes
    q^{h_{\varepsilon_c}} \notag \\
& \qquad + (q-q^{-1}) \sum_{\varepsilon_{b} > \varepsilon_{a}} (-1)^{[b]} E^a_b
    \otimes q^{h_{\varepsilon_b}} \otimes q^{h_{\varepsilon_a}} \hat{\sigma}
    _{ba} \notag \\
& \qquad + (q-q^{-1})^2 \sum_{\varepsilon_{d} > \varepsilon_{c} >\varepsilon_a}
    (-1)^{[c]+[d]} E^a_d \otimes q^{h_{\varepsilon_c}} \hat{\sigma}_{dc} 
    \otimes q^{h_{\varepsilon_a}} \hat{\sigma}_{ca} \notag \\
&= \sum_a E^a_a \otimes q^{h_{\varepsilon_a}} \otimes q^{h_{\varepsilon_a}} 
   \notag \\
& \qquad + (q-q^{-1}) \sum_{\varepsilon_{b} > \varepsilon_{a}} (-1)^{[b]} E^a_b
    \otimes \bigl( q^{h_{\varepsilon_a}} \hat{\sigma}_{ba} \otimes 
    q^{h_{\varepsilon_a}} + q^{h_{\varepsilon_b}} \otimes q^{h_{\varepsilon_a}}
    \hat{\sigma}_{ba} \bigr) \notag \\ 
& \qquad + (q-q^{-1})^2 \sum_{\varepsilon_{b} > \varepsilon_{c} >\varepsilon_a}
    (-1)^{[b]+[c]} E^a_b \otimes q^{h_{\varepsilon_c}} \hat{\sigma}_{bc} 
    \otimes q^{h_{\varepsilon_a}} \hat{\sigma}_{ca}. \label{R13R12}
\end{align*}

\noindent Also, the coproduct properties \eqref{coprod} imply

\begin{equation*} 
(\text{id} \otimes \Delta) R = \sum_a E^a_a \otimes q^{h_{\varepsilon_a}} \otimes q^{h_{\varepsilon_a}} + (q-q^{-1}) \sum_{\varepsilon_{b} > \varepsilon_{a}} (-1)^{[b]} E^a_b \otimes (q^{h_{\varepsilon_a}} \otimes q^{h_{\varepsilon_a}}) \Delta(\hat{\sigma}_{ba}).
\end{equation*}

\noindent  Hence $R$ will satisfy equation \eqref{1tensordel} if and only if $\Delta(\hat{\sigma}_{ba})$ is given by:

\begin{equation*}
\Delta(\hat{\sigma}_{ba}) = \hat{\sigma}_{ba} \otimes I + q^{h_{\varepsilon_b} 
  - h_{\varepsilon_a}} \otimes \hat{\sigma}_{ba} + (q-q^{-1}) 
  \sum_{\varepsilon_{b} > \varepsilon_{c} >\varepsilon_a} (-1)^{[c]} 
  q^{h_{\varepsilon_c} - h_{\varepsilon_a}} \hat{\sigma}_{bc} \otimes 
  \hat{\sigma}_{ca}.
\end{equation*}

\noindent  Now we use the fundamental values of $\hat{\sigma}_{ba}$ \eqref{ivalues} and the inductive relations \eqref{indrel} to calculate $\Delta(\hat{\sigma}_{ba})$, and show that it is indeed of this form.  First set

\begin{equation*} 
h_{ba} \equiv h_{\varepsilon_b}- h_{\varepsilon_a},
\end{equation*}

\noindent so we need to show 

\begin{equation*} 
\Delta(\hat{\sigma}_{ba}) = \hat{\sigma}_{ba} \otimes I + q^{h_{ba}} \otimes 
  \hat{\sigma}_{ba} + (q-q^{-1}) \sum_{\varepsilon_{b} > \varepsilon_{c} >
  \varepsilon_a}(-1)^{[c]} q^{h_{ca}}\hat{\sigma}_{bc} \otimes\hat{\sigma}_{ca}.
\end{equation*}

\

\noindent Consider the non-zero fundamental values of $\hat{\sigma}_{ba}$, given in equation \eqref{ivalues}.  For these values $\alpha_b = \varepsilon_b - \varepsilon_a$ is a simple root.  Note that in each case $\hat{\sigma}_{ba} = A e_b q^{\frac{1}{2} h_{ba}}$ or $A e_{\overline{a}} q^{\frac{1}{2} h_{ba}}$ for some constant $A$.  Then

\begin{alignat*}{2}
\Delta(\hat{\sigma}_{ba}) &= A \Delta(e_c) \Delta(q^{\frac{1}{2} h_{ba}}) 
  , &&c = b \text{ or } \overline{a} \notag \\
&= A (q^{\frac{1}{2}h_{ba}} \otimes e_c + e_c \otimes q^{-\frac{1}{2} h_{ba}}) 
  (q^{\frac{1}{2} h_{ba}} \otimes q^{\frac{1}{2} h_{ba}})&\qquad& \notag \\
&= q^{h_{ba}} \otimes A e_c q^{\frac{1}{2} h_b} + A e_c q^{\frac{1}{2} h_b} 
  \otimes I&& \notag \\
&= q^{h_{ba}} \otimes \hat{\sigma}_{ba} + \hat{\sigma}_{ba} \otimes I.&&
\end{alignat*}

\noindent In the case of a simple root there is usually no $c$ satisfying $\varepsilon_b > \varepsilon_c > \varepsilon_a$, so this is the expected result.  The only exceptions to that generalisation are $\hat{\sigma}_{l\, \overline{l-1}}$ and $\hat{\sigma}_{l-1\, \overline{l}}$ where $m=2l$.  In both those cases, however, the sum in our expression for $\Delta (\hat{\sigma}_{ba})$ still disappears as it becomes a single term containing $\hat{\sigma}_{l\, \overline{l}}$, which we know equals 0.  Hence our formula for the expected coproduct is correct for the non-zero fundamental values of $\sigma_{ba}$.

\

\noindent Also, when $m=2l$

\begin{equation*}
\hat{\sigma}_{l\overline{l}} \otimes I + q^{h_{l\overline{l}}} \otimes 
  \hat{\sigma}_{l\overline{l}} + (q-q^{-1}) \sum_{\varepsilon_l > 
  \varepsilon_c > -\varepsilon_l} (-1)^{[c]} q^{h_{c\overline{l}}} \hat{\sigma}
  _{lc} \otimes \hat{\sigma}_{c\overline{l}} = 0+0+0 = \Delta(\hat{\sigma}_
  {l \overline{l}} =0)
\end{equation*}

\noindent as required.  Thus we have verified the formula for the coproduct for all the fundamental values of $\hat{\sigma}_{ba}$ given in equation \eqref{ivalues}.

\

\noindent To find the coproduct for the remaining values of $\hat{\sigma}_{ba}$ we use the inductive relations \eqref{indrel}:

\begin{equation*}
\hat{\sigma}_{ba} = q^{-(\varepsilon_b,\varepsilon_a)} \hat{\sigma}_{bc} 
  \hat{\sigma}_{ca}-q^{-(\varepsilon_c,\varepsilon_c)}(-1)^{([b]+[c])([a]+[c])}
  \hat{\sigma}_{ca} \hat{\sigma}_{bc}, \quad \varepsilon_b > \varepsilon_c > 
  \varepsilon_a,
\end{equation*}

\noindent  where $c \neq \overline{b}$ or $\overline{a}$.  We assume our formula for the coproduct holds for $\hat{\sigma}_{bc}$ and $\hat{\sigma}_{ca}$, where $\varepsilon_b > \varepsilon_c > \varepsilon_a$, and then show it is also true for $\hat{\sigma}_{ba}$.

\

\noindent We can always choose $c$ satisfying the conditions such that either $\varepsilon_b - \varepsilon_c$ or $\varepsilon_c - \varepsilon_a$ is a simple root.  First consider $\varepsilon_b - \varepsilon_c$ is a simple root, denoted by either $\alpha_b$ or $\alpha_{\overline{c}}$ depending on circumstance, so $\hat{\sigma}_{bc} = A e_b q^{\frac{1}{2} h_b}$ or $A e_{\overline{c}} q^{\frac{1}{2} h_{\overline{c}}}$ for some constant $A$.

\

\noindent The coproduct is an algebra homomorphism, so for $\varepsilon_b > \varepsilon_c > \varepsilon_a$, $c \neq \overline{a}$ or $ \overline{b}$, we have

\begin{equation*}
\Delta(\hat{\sigma}_{ba}) = q^{-(\varepsilon_b,\varepsilon_a)} \Delta(\hat
  {\sigma}_{bc}) \Delta(\hat{\sigma}_{ca}) -q^{-(\varepsilon_c, \varepsilon_c)}
  (-1)^{([b]+[c])([a]+[c])} \Delta(\hat{\sigma}_{ca})\Delta(\hat{\sigma}_{bc}).
\end{equation*}

\noindent Substituting in our expression for the coproduct gives:

\begin{align*}
\Delta( \hat{\sigma}_{ba}) = &q^{-(\varepsilon_b,\varepsilon_a)} (\hat{\sigma}
   _{bc} \otimes I + q^{h_{bc}} \otimes \hat{\sigma}_{bc}) \notag \\
&\qquad \bigl( \hat{\sigma}_{ca} \otimes I + q^{h_{ca}} \otimes 
   \hat{\sigma}_{ca} + (q-q^{-1}) \sum_{\varepsilon_{c} > \varepsilon_{d} >
   \varepsilon_a} (-1)^{[d]} q^{h_{da}} \hat{\sigma}_{cd} \otimes 
   \hat{\sigma}_{da} \bigr) \notag \\
&- q^{-(\varepsilon_c, \varepsilon_c)} (-1)^{([b]+[c])([a]+[c])}\notag\\
&\qquad \bigl(\hat{\sigma}_{ca} \otimes I + q^{h_{ca}} \otimes 
   \hat{\sigma}_{ca} + (q-q^{-1}) \sum_{\varepsilon_{c} > \varepsilon_{d} >
   \varepsilon_a} (-1)^{[d]} q^{h_{da}} \hat{\sigma}_{cd} \otimes 
   \hat{\sigma}_{da} \bigr) \notag \\
&\qquad (\hat{\sigma}_{bc} \otimes I + q^{h_{bc}} \otimes \hat{\sigma}_{bc}). 
\end{align*}

\noindent Expanding, we obtain

\begin{align*}
\Delta( \hat{\sigma}_{ba})&=(q^{-(\varepsilon_b,\varepsilon_a)} \hat{\sigma}_{bc} \hat{\sigma}_{ca} - 
   q^{-(\varepsilon_c, \varepsilon_c)} (-1)^{([b]+[c])([a]+[c])} \hat{\sigma}_
   {ca} \hat{\sigma}_{bc}) \otimes I \\
& \quad+ q^{h_{ba}} \otimes (q^{-(\varepsilon_b, \varepsilon_a)} \hat{\sigma}_
   {bc} \hat{\sigma}_{ca} - q^{-(\varepsilon_c, \varepsilon_c)} (-1)^{([b]+[c])
   ([a]+[c])} \hat{\sigma}_{ca} \hat{\sigma}_{bc}) \\
& \quad+ (q^{-(\varepsilon_b,\varepsilon_a)} q^{-(\varepsilon_c- \varepsilon_a,
   \varepsilon_b - \varepsilon_c)} - q^{-(\varepsilon_c, \varepsilon_c)}
   [(-1)^{([b]+[c])([a]+[c])}]^2) q^{h_{ca}} \hat{\sigma}_{bc} \otimes 
   \hat{\sigma}_{ca} \\
& \quad+ \bigl(q^{-(\varepsilon_b,\varepsilon_a)} (-1)^{([b]+[c])([a]+[c])} \\
& \hspace{3cm} -q^{-(\varepsilon_c, \varepsilon_c)}(-1)^{([b]+[c])([a]+[c])} 
   q^{-(\varepsilon_b- \varepsilon_c, \varepsilon_c - \varepsilon_a)}\bigr)
   q^{h_{bc}} \hat{\sigma}_{ca} \otimes \hat{\sigma}_{bc} \\
& \quad+ (q-q^{-1}) \sum_{\varepsilon_{c} > \varepsilon_{d} >\varepsilon_a} 
   (-1)^{[d]} q^{h_{da}} q^{-(\varepsilon_b,\varepsilon_a)} q^{-(\varepsilon_d
   - \varepsilon_a, \varepsilon_b - \varepsilon_c)} \hat{\sigma}_{bc} 
   \hat{\sigma}_{cd} \otimes \hat{\sigma}_{da} \\
& \quad- (q-q^{-1}) \sum_{\varepsilon_{c} > \varepsilon_{d} >\varepsilon_a} 
   (-1)^{[d]} q^{h_{da}} q^{-(\varepsilon_c, \varepsilon_c)} (-1)^{([b]+[c])
   ([c]+[d])} \hat{\sigma}_{cd} \hat{\sigma}_{bc} \otimes \hat{\sigma}_{da} \\
& \quad+ (q-q^{-1}) \sum_{\varepsilon_{c} > \varepsilon_{d} >\varepsilon_a} 
   (-1)^{[d]} q^{h_{bc}+h_{da}} q^{-(\varepsilon_b, \varepsilon_a)} 
   (-1)^{([b]+[c])([c]+[d])} \hat{\sigma}_{cd} \otimes \hat{\sigma}_{bc} 
   \hat{\sigma}_{da} \\
& \quad- (q-q^{-1}) \sum_{\varepsilon_{c} > \varepsilon_{d} >\varepsilon_a} 
   (-1)^{[d]} q^{h_{bc}+h_{da}} q^{-(\varepsilon_c, \varepsilon_c)}
   (-1)^{([b]+[c])([a]+[c])}\\
&\hspace{8cm} \times q^{-(\varepsilon_b-\varepsilon_c, 
   \varepsilon_c - \varepsilon_d)}\hat{\sigma}_{cd} \otimes \hat{\sigma}_{da} 
   \hat{\sigma}_{bc}.
\end{align*}

\noindent Since $\varepsilon_b > \varepsilon_c > \varepsilon_a$ and $c \neq \overline{a}$ or $\overline{b}$, we know that $(\varepsilon_a, \varepsilon_c) = (\varepsilon_b, \varepsilon_c) = 0$.  Using this, we simplify the above expression to:

\begin{align}
\Delta( \hat{\sigma}_{ba}) = &\hat{\sigma}_{ba} \otimes I + q^{h_{ba}} \otimes 
   \hat{\sigma}_{ba} + (q^{(\varepsilon_c,\varepsilon_c)} - q^{-(\varepsilon_c,
   \varepsilon_c)}) q^{h_{ca}} \hat{\sigma}_{bc} \otimes \hat{\sigma}_{ca} 
   \notag \\
&+ (q^{-(\varepsilon_b, \varepsilon_a)} - q^{(\varepsilon_b, 
   \varepsilon_a)}) (-1)^{([b]+[c])([a]+[c])} q^{h_{bc}} \hat{\sigma}_{ca} 
   \otimes \hat{\sigma}_{bc} \notag \\
&+ (q-q^{-1}) \sum_{\varepsilon_{c} > \varepsilon_{d} > \varepsilon_a} 
   (-1)^{[d]} q^{h_{da}} \bigl(q^{-(\varepsilon_d, \alpha_b)} \hat{\sigma}_{bc}
   \hat{\sigma}_{cd}\notag  \\
& \hspace{58mm}- q^{-(\varepsilon_c, \varepsilon_c)} (-1)^{([b]+[c])
   ([c]+[d])} \hat{\sigma}_{cd} \hat{\sigma}_{bc} \bigr) \otimes 
   \hat{\sigma}_{da} \notag \\
&+ (q-q^{-1}) \sum_{\varepsilon_{c} > \varepsilon_{d} >\varepsilon_a} 
    (-1)^{[d]} q^{h_{da} + h_{bc}} (-1)^{([b]+[c])([a]+[c])} \hat{\sigma}_{cd} 
   \notag \\
& \hspace{35mm}\otimes\bigl( q^{-(\varepsilon_b-\varepsilon_c, \varepsilon_a)} 
    (-1)^{([b]+[c])([a]+[d])} \hat{\sigma}_{bc} \hat{\sigma}_{da} - 
    q^{(\varepsilon_d, \alpha_b)} \hat{\sigma}_{da} \hat{\sigma}_{bc}\bigr).
   \label{rhs2}
\end{align}

\noindent But when $d \neq \overline{b}$ the $q$-commutation relations \eqref{qcom} can be used to show

\begin{align*}
q^{-(\alpha_b, \varepsilon_a)} (-1)^{([b]+[c])([a]+[d])} 
&  \hat{\sigma}_{bc} \hat{\sigma}_{da} - q^{(\varepsilon_d, \alpha_b)} 
   \hat{\sigma}_{da} \hat{\sigma}_{bc} \notag \\
= &-A (q^{(\varepsilon_d, \alpha_b)} \hat{\sigma}_{da} e_b q^{\frac{1}{2} h_b}
   - (-1)^{([a]+[d])[\alpha_b]} q^{-(\alpha_b, \varepsilon_a)} e_b 
   q^{\frac{1}{2} h_b} \hat{\sigma}_{da}) \notag \\
\bigl(\text{or }&-A(q^{(\varepsilon_d,\alpha_{\overline{c}})} \hat{\sigma}_{da}
  e_{\overline{c}} q^{\frac{1}{2} h_{\overline{c}}}
  - (-1)^{([a]+[d])[\alpha_{\overline{c}}]} q^{-(\alpha_{\overline{c}}, 
  \varepsilon_a)} e_b q^{\frac{1}{2} h_{\overline{c}}} \hat{\sigma}_{da})\bigr)
   \notag \\
= &\;0.
\end{align*}

\noindent  And in the case $c \neq \overline{d}$, we have 

\begin{equation*}
q^{-(\varepsilon_d, \alpha_b)} \hat{\sigma}_{bc} \hat{\sigma}_{cd} - 
q^{-(\varepsilon_c, \varepsilon_c)} (-1)^{([b]+[c])([c]+[d])} \hat{\sigma}_{cd}
\hat{\sigma}_{bc} = \hat{\sigma}_{bd} 
\end{equation*}

\noindent from the inductive relations \eqref{indrel}.  Moreover, when $c = \overline{d}$ (so $\varepsilon_c > 0$) we note from the relations in Table \ref{list} that

\begin{align*}
q^{-(\varepsilon_d, \alpha_b)} \hat{\sigma}_{bc} \hat{\sigma}_{cd} &- 
  q^{-(\varepsilon_c, \varepsilon_c)} (-1)^{([b]+[c])([c]+[d])} 
  \hat{\sigma}_{cd} \hat{\sigma}_{bc} \notag \\
&= q^{-(\varepsilon_c, \varepsilon_c)} [\hat{\sigma}_{bc}, \hat{\sigma}_{c\, 
  \overline{c}}] \notag \\
&= \hat{\sigma}_{b \overline{c}} + (-1)^{[b]} [(-1)^k q]^{\delta^c_{i=1}} 
  \hat{\sigma}_{c \overline{b}} \notag \\
&= \hat{\sigma}_{bd} + (-1)^{[b]} [(-1)^k q]^{\delta^c_{i=1}} \hat{\sigma}_{c
  \overline{b}},
\end{align*}

\noindent so for all $d$ satisfying $\varepsilon_c > \varepsilon_d > \varepsilon_a$ we have

\begin{equation*}
q^{-(\varepsilon_d, \alpha_b)} \hat{\sigma}_{bc} \hat{\sigma}_{cd} - 
  q^{-(\varepsilon_c, \varepsilon_c)} (-1)^{([b]+[c])([c]+[d])} 
  \hat{\sigma}_{cd} \hat{\sigma}_{bc} = \hat{\sigma}_{bd} + \delta^
  {\overline{c}}_d (-1)^{[b]} [(-1)^k q]^{\delta^c_{i=1}} \hat{\sigma}_{c
  \overline{b}}.
\end{equation*}

\noindent We also introduce a new function $\theta_{xy}$, defined by

\begin{equation*}
\theta_{xy} =
  \begin{cases}
  1, \quad & \varepsilon_x > \varepsilon_y, \\
  0, & \varepsilon_x \leq \varepsilon_y.
  \end{cases}
\end{equation*}

\noindent Combining all this information, we simplify equation \eqref{rhs2} to:

\begin{align}
\Delta( \hat{\sigma}_{ba})= &\hat{\sigma}_{ba} \otimes I + q^{h_{ba}} \otimes 
   \hat{\sigma}_{ba} + (1 - \delta^c_{i= l+1}) (q-q^{-1}) (-1)^{[c]} q^{h_{ca}}
   \hat{\sigma}_{bc} \otimes \hat{\sigma}_{ca} \notag  \\
&+ \delta^a_{\overline{b}} (q-q^{-1}) (-1)^{[c]} q^{h_{bc}} 
   \hat{\sigma}_{ca} \otimes \hat{\sigma}_{bc} \notag \\
&+ (q-q^{-1}) \sum_{\varepsilon_{c} > \varepsilon_{d} >\varepsilon_a} 
   (-1)^{[d]} q^{h_{da}} \hat{\sigma}_{bd} \otimes \hat{\sigma}_{da} \notag \\
&+ (q-q^{-1}) (-1)^{[c]+[b]} \theta_{c\overline{c}} \theta_{\overline{c}a} 
   q^{h_{\overline{c}a}}[(-1)^kq]^{\delta^c_{i=1}} \hat{\sigma}_{c\overline{b}}
   \otimes \hat{\sigma}_{\overline{c} a} \notag \\
&+ (q-q^{-1}) (-1)^{[c]} \theta_{c\overline{b}} \theta_{\overline{b}a} 
   q^{h_{\overline{c}a}} \hat{\sigma}_{c\overline{b}} \otimes(\hat{\sigma}_{bc}
   \hat{\sigma}_{\overline{b}a} -(-1)^{([b]+[c])([a]+[b])} q^{-(\varepsilon_b,
   \varepsilon_b)} \hat{\sigma}_{\overline{b}a} \hat{\sigma}_{bc}).\label{rhs3}
\end{align}

\noindent  While this currently does not look much like the expected formula for $\Delta(\hat{\sigma}_{ba})$, it can be further simplified.  First note that since $\varepsilon_b - \varepsilon_c = \varepsilon_{\overline{c}} - \varepsilon_{\overline{b}}$ is a simple root, 

\begin{align*}
\theta_{c\overline{c}} \theta_{\overline{c}a} &= \theta_{c\overline{c}} \theta_
  {\overline{c}\overline{b}} \theta_{\overline{b}a} + \delta^a_{\overline{b}} 
  \theta_{c\overline{c}} \theta_{\overline{c}a} + \theta_{c\overline{c}} \theta
  _{\overline{c}a} \theta_{a\overline{b}} \notag \\
& = \theta_{c\overline{c}} \theta_{\overline{b}a} + \delta^a_{\overline{b}} 
  \theta_{c\overline{c}} 
\end{align*}

\noindent and 

\begin{align*}
\theta_{c\overline{b}} \theta_{\overline{b}a} &= \theta_{c\overline{c}} 
  \theta_{\overline{c} \overline{b}} \theta_{\overline{b}a} + \delta^c_{
  \overline{c}} \theta_{c\overline{b}} \theta_{\overline{b}a} + \theta_
  {\overline{c}c} \theta_{c\overline{b}} \theta_{\overline{b}a} \notag \\
&= \theta_{c\overline{c}} \theta_{\overline{b}a} + \delta^c_{i=l+1} 
  \theta_{\overline{l}a} + \delta^b_{i=l-1} \delta^c_{j=\overline{l}} 
  \theta_{\overline{l-1}\,a}.
\end{align*}

\noindent Also, looking back at the formulae for $\hat{\sigma}_{bc}$ associated with the simple roots, we see that if $\delta^c_{i=l+1}=1$ we have $\sigma_{bc} = - q^{-\frac{1}{2}} \sigma_{\overline{c}\overline{b}}$; if \hbox{$\delta^b_{i=l-1} \delta^c_{j=\overline{l}}=1$} then \hbox{$\sigma_{bc} = - \sigma_{\overline{c}\overline{b}}$}; and if \hbox{$\theta_{c\overline{c}} =1$} then \hbox{$\sigma_{bc} = -(-1)^{[b]} [(-1)^kq]^{\delta^c_{i=1}} \sigma_{\overline{c}\overline{b}}$}.  Moreover, if $\theta_{c\overline{b}} \theta_{\overline{b}a} =1$ then $b \neq \overline{c},a$ so we can simplify the final term in \eqref{rhs3} using

\begin{equation*}
\hat{\sigma}_{\overline{c} a} = \hat{\sigma}_{\overline{c}\, \overline{b}} \hat{\sigma}_{\overline{b}a} - (-1)^{([b]+[c])([a]+[b])} q^{-(\varepsilon_b,\varepsilon_b)} \hat{\sigma}_{\overline{b}a} \hat{\sigma}_{\overline{c}\, \overline{b}}.
\end{equation*}

\noindent  Applying all this gives:

\begin{align*}
\Delta( \hat{\sigma}_{ba}) = & \hat{\sigma}_{ba} \otimes I + q^{h_{ba}} \otimes
  \hat{\sigma}_{ba} + (q-q^{-1}) \sum_{\varepsilon_{c} \geq \varepsilon_{d} >
  \varepsilon_a} (-1)^{[d]} q^{h_{da}} \hat{\sigma}_{bd} \otimes 
  \hat{\sigma}_{da} \notag \\
&-\delta^c_{i= l+1} (q-q^{-1}) q^{h_{ca}} \hat{\sigma}_{bc} \otimes 
  \hat{\sigma}_{ca} + \delta^a_{\overline{b}} (q-q^{-1}) (-1)^{[c]} 
  q^{h_{bc}} \hat{\sigma}_{ca} \otimes \hat{\sigma}_{bc} \notag \\
&+(q-q^{-1}) (-1)^{[b]+[c]} (\theta_{c\overline{c}} \theta_{\overline{b}a} + 
   \delta^a_{\overline{b}} \theta_{c\overline{c}}) q^{h_{\overline{c}a}}
   [(-1)^kq]^{\delta^c_{i=1}} \hat{\sigma}_{c\overline{b}} \otimes 
   \hat{\sigma}_{\overline{c} a} \notag \\
&-(q-q^{-1}) (-1)^{[b]+[c]} [(-1)^kq]^{\delta^c_{i=1}} \theta_{c\overline{c}} 
   \theta_{\overline{b}a} q^{h_{\overline{c}a}} \hat{\sigma}_{c\overline{b}}
   \otimes \hat{\sigma}_{\overline{c}a} \notag \\
&-\delta^c_{i=l+1} (q-q^{-1}) q^{-\frac{1}{2}} \theta_{\overline{l}a} q^{h_
   {ca}} \hat{\sigma}_{\overline{c}\,\overline{b}} \otimes \hat{\sigma}_{ca} 
   \notag \\
&-\delta^b_{i=l-1} \delta^c_{j=\overline{l}} (q-q^{-1}) 
   \theta_{\overline{l-1}a} q^{h_{la}} \hat{\sigma}_{\overline{l}\, 
   \overline{l-1}} \otimes \hat{\sigma}_{la} \notag \\
= & \hat{\sigma}_{ba} \otimes I + q^{h_{ba}} \otimes \hat{\sigma}_{ba} + 
   (q-q^{-1}) \sum_{\varepsilon_{c} \geq \varepsilon_{d} >\varepsilon_a} 
   (-1)^{[d]} q^{h_{da}} \hat{\sigma}_{bd} \otimes \hat{\sigma}_{da} \notag \\
&-\delta^c_{i= l+1} \delta^a_{\overline{b}} (q-q^{-1}) (-1)^{[c]} q^{h_{bc}} 
   \hat{\sigma}_{ca} \otimes \hat{\sigma}_{bc} + \delta^a
   _{\overline{b}} (q-q^{-1}) (-1)^{[c]} q^{h_{bc}} \hat{\sigma}_{ca} \otimes 
   \hat{\sigma}_{bc} \notag \\
&-\delta^a_{\overline{b}} (q-q^{-1}) (-1)^{[c]} \theta_{c\overline{c}}
  q^{h_{bc}} \hat{\sigma}_{ca} \otimes \hat{\sigma}_{bc} \notag \\
&+\delta^b_{i=l-1} \delta^c_{j=\overline{l}} (q-q^{-1}) \theta_{\overline
  {l-1}a} q^{h_{la}} \hat{\sigma}_{l-1\,l} \otimes \hat{\sigma}_{la} .
\end{align*}

\noindent  Now note that 

\begin{equation*}
\delta^a_{\overline{b}} (1 - \delta^c_{i=l+1} - \theta_{c\overline{c}}) =
\delta^a_{\overline{b}} \delta^b_{i=l-1} \delta^c_{j=\overline{l}}
\end{equation*}

\noindent and that $\hat{\sigma}_{l\overline{l}} = 0$.  Then our formula for $\Delta(\hat{\sigma}_{ba})$ becomes:

\begin{align*}
\Delta( \hat{\sigma}_{ba}) = & \hat{\sigma}_{ba} \otimes I + q^{h_{ba}} \otimes
  \hat{\sigma}_{ba} + (q-q^{-1}) \sum_{\varepsilon_{c} \geq \varepsilon_{d} >
  \varepsilon_a} (-1)^{[d]} q^{h_{da}} \hat{\sigma}_{bd} \otimes 
  \hat{\sigma}_{da} \notag \\
&+\delta^a_{\overline{b}} \delta^b_{i=l-1} \delta^c_{j=\overline{l}} (q-q^{-1})
  (-1)^{[c]} q^{h_{bc}} \hat{\sigma}_{ca} \otimes \hat{\sigma}_{bc}  \notag \\
&+\delta^b_{i=l-1} \delta^c_{j=\overline{l}} (q-q^{-1}) \theta_{\overline{l-1}
  a} q^{h_{la}} \hat{\sigma}_{l-1\,l} \otimes \hat{\sigma}_{la} \notag \\
&+ \delta^b_{i=l} \delta^c_{j=\overline{l-1}} (q-q^{-1}) \theta_{\overline{l}a}
  q^{h_{\overline{l}a}} \hat{\sigma}_{l\overline{l}} \otimes 
  \hat{\sigma}_{\overline{l}a} \notag \\
= &\hat{\sigma}_{ba} \otimes I + q^{h_{ba}} \otimes \hat{\sigma}_{ba} + 
  (q-q^{-1}) \sum_{\varepsilon_c\geq \varepsilon_d > \varepsilon_a} (-1)^{[d]} 
  q^{h_{da}} \hat{\sigma}_{bd} \otimes \hat{\sigma}_{da} \notag \\
&+ \delta^b_{i=l-1} \delta^c_{j=\overline{l}} (q-q^{-1}) \theta_{\overline{l}a}
  q^{h_{la}} \hat{\sigma}_{l-1\,l} \otimes \hat{\sigma}_{la} \notag \\
&+ \delta^b_{i=l} \delta^c_{j=\overline{l-1}} (q-q^{-1}) \theta_{\overline{l-1}
  a} q^{h_{\overline{l}a}} \hat{\sigma}_{l\overline{l}} \otimes \hat{\sigma}_
  {\overline{l}a} \notag \\
=&\hat{\sigma}_{ba} \otimes I + q^{h_{ba}} \otimes 
  \hat{\sigma}_{ba} + (q-q^{-1}) \sum_{\varepsilon_b > \varepsilon_d > 
  \varepsilon_a} (-1)^{[d]} q^{h_{da}} \hat{\sigma}_{bd} \otimes \hat{\sigma}
  _{da}
\end{align*}

\noindent as required.

\

\noindent  To verify our formula for the coproduct it is also necessary to consider the case when $\hat{\sigma}_{ca}$ is a fundamental value.  The calculations, however, are extremely similar to those where $\hat{\sigma}_{bc}$ is a fundamental value, so they are not included.  Suffice it to say that they give the expected result.  Moreover, as a check, it has also been shown directly that the coproduct is consistent with the commutation relations \eqref{qcom}, although again the calculations are rather tedious and have been omitted.

\

\noindent Thus we have shown the coproduct of the operators $\hat{\sigma}_{ba}$ is given by

\begin{equation*}
\Delta(\hat{\sigma}_{ba}) = \hat{\sigma}_{ba} \otimes I + q^{h_{ba}} \otimes 
  \hat{\sigma}_{ba} + (q-q^{-1}) \sum_{\varepsilon_{b} > \varepsilon_{c} >
  \varepsilon_a}(-1)^{[c]} q^{h_{ca}}\hat{\sigma}_{bc} \otimes\hat{\sigma}_{ca},
\end{equation*}

\noindent and consequently that the matrix $R$ found in the previous chapter satisfies the property

\begin{equation*}
(\text{id} \otimes \Delta) R = R_{13} R_{12}.
\end{equation*}

\section{The Intertwining Property}
    
\noindent  To confirm that we have a Lax operator we need to check one last relation, namely the intertwining property for the other generators. 

\begin{equation} \label{intertwine}
R \Delta(a) = \Delta^T (a) R, \qquad \forall a \in U_q[osp(m|n)].
\end{equation}

\noindent Now $R$ is weightless, so it commutes with all the Cartan elements.  Moreover, $\Delta(q^{h_a}) = \Delta^T (q^{h_a}), \forall h_a \in H$, so the Cartan elements will automatically satisfy equation \eqref{intertwine}.  Thus it remains only to verify the intertwining property for the lowering generators, $f_a$.  Unfortunately, knowing the raising generators satisfy the intertwining property does not appear helpful.  Instead, we start by assuming the form of the Lax operator and that it satisfies the intertwining property for the lowering generators, and then proceed as in the previous chapter. Provided the relations and fundamental values obtained are consistent with those already developed, we will have confirmed that the matrix $R$ constructed in the previous chapter is a Lax operator.  Initially the process mirrors that in Section \ref{Developing Relations}, so some of the detail is omitted.  

\

\noindent Now we know

\begin{equation*}
R \equiv q ^ {\underset{a}{\sum} h_{a} \otimes h^{a}} \Bigl[ I \otimes I + (q - q^{-1}) \sum_{\varepsilon_{a} < \varepsilon_{b}} (-1)^{[b]} E^a_b \otimes \hat{\sigma}_{ba} \Bigr]
\end{equation*} 

\noindent and 

\begin{equation*}
\Delta(f_c) = q^{\frac{1}{2}h_c} \otimes f_c + f_c \otimes q^{-\frac{1}{2}h_c}.
\end{equation*}

\noindent Moreover,

\begin{align*}
\Delta^T(f_c) q^{\underset{a}{\sum} h_a \otimes h^a} 
&= (q^{-\frac{1}{2}h_c} \otimes f_c + f_c \otimes q^{\frac{1}{2}h_c}) 
   q^{\underset{a}{\sum} h_a \otimes h^a} \notag \\
&= q^{\underset{a}{\sum} h_a \otimes h^a} (q^{\frac{1}{2}h_c} \otimes f_c + f_c
   \otimes q^{\frac{3}{2}h_c}).
\end{align*}

\noindent Therefore

\begin{align}
\Delta^T&(f_c) R \notag \\
&= q^{\underset{a}{\sum} h_a \otimes h^a} (q^{\frac{1}{2}h_c} \otimes f_c + f_c
   \otimes q^{\frac{3}{2}h_c}) \Bigl[ I \otimes I + (q - q^{-1}) \sum_
   {\varepsilon_{a} < \varepsilon_{b}} (-1)^{[b]} E^a_b \otimes 
   \hat{\sigma}_{ba} \Bigr]\notag \\
&= q^{\underset{a}{\sum} h_a \otimes h^a} \Bigl\{ q^{\frac{1}{2}h_c}\otimes f_c
   + f_c \otimes q^{\frac{3}{2} h_c} \notag \\
&\quad+ (q - q^{-1}) \sum_{\varepsilon_{a} < \varepsilon_{b}} (-1)^{[b]}
   \bigl[ (-1)^{([a] + [b])[c]} q^{\frac{1}{2} (\alpha_c, \varepsilon_a)} E^a_b
   \otimes f_c \hat{\sigma}_{ba} + f_c E^a_b \otimes q^{\frac{3}{2} h_c}
   \hat{\sigma}_{ba} \bigr] \Bigr\},  \label{DelRf}
\end{align}

\noindent while

\begin{align}
R \Delta (f_c) = q^{\underset{a}{\sum} h_a \otimes h^a} \Bigl\{q^{\frac{1}{2} 
   h_c} \otimes f_c + f_c \otimes q^{-\frac{1}{2} h_c} \qquad &\notag \\
+ (q - q^{-1}) \sum_{\varepsilon_{a} < \varepsilon_{b}} (-1)^{[b]} \bigl[ 
   &q^{\frac{1}{2} (\alpha_c, \varepsilon_b)} E^a_b \otimes \hat{\sigma}_{ba} 
   f_c \notag \\ 
&+ (-1)^{([a]+[b])[c]}  E^a_b f_c \otimes \hat{\sigma}_{ba} 
   q^{-\frac{1}{2} h_c} \bigr] \Bigr\}. \label{Rdelf} 
\end{align}

\noindent Equating (\ref{DelRf}) and (\ref{Rdelf}), we find

\begin{multline}
f_c \otimes (q^{\frac{3}{2}h_c} - q^{-\frac{1}{2} h_c}) \\
= (q - q^{-1}) 
   \sum_{\varepsilon_{a} < \varepsilon_{b}} (-1)^{[b]} \Bigl[ q^{\frac{1}{2} 
   (\alpha_c, \varepsilon_b)} E^a_b \otimes \hat{\sigma}_{ba} f_c - (-1)^{
   ([a] + [b])[c]} q^{\frac{1}{2} (\alpha_c, \varepsilon_a)} E^a_b \otimes f_c 
   \hat{\sigma}_{ba}\Bigr] \\
+ (q - q^{-1}) \sum_{\varepsilon_{a} < \varepsilon_{b}} (-1)^{[b]} \Bigl[
   (-1)^{([a]+[b])[c]}  E^a_b f_c \otimes \hat{\sigma}_{ba} q^{-\frac{1}{2} 
   h_c} -  f_c E^a_b \otimes q^{\frac{3}{2} h_c} \hat{\sigma}_{ba} \Bigr]. 
 \label{f}
\end{multline}

\noindent Taking the terms with zero weight on the right-hand side of the tensor product gives

\begin{multline}
f_c \otimes (q^{\frac{3}{2}h_c} - q^{-\frac{1}{2} h_c})  \\
= (q - q^{-1}) \sum_{\varepsilon_{b} - \varepsilon_{a} = \alpha_c} (-1)^{[b]}
    E^a_b \otimes \bigl( q^{\frac{1}{2} (\alpha_c, \varepsilon_b)} 
    \hat{\sigma}_{ba} f_c - (-1)^{[c]} q^{\frac{1}{2} (\alpha_c, 
    \varepsilon_a)} f_c \hat{\sigma}_{ba} \bigr). \label{basef}
\end{multline}

\noindent This can be used to find the fundamental values of $\hat{\sigma}_{ba}$, and check that they agree with those in Section \ref{iv}.

\

\noindent Similarly, taking the terms of Equation (\ref{f}) with non-zero weight on the right-hand side of the tensor product, we find

\begin{multline*}
\underset{\varepsilon_{b} - \varepsilon_{a} \neq \alpha_c} {\sum_{\varepsilon_
  {b} > \varepsilon_{a}}} (-1)^{[b]} E^a_b \otimes \bigl( q^{\frac{1}{2} 
  (\alpha_c, \varepsilon_b)} \hat{\sigma}_{ba} f_c - (-1)^{([a]+[b])[c]} 
  q^{\frac{1}{2} (\alpha_c, \varepsilon_a)} f_c \hat{\sigma}_{ba} \bigr) \\
= \sum_{\varepsilon_{b} > \varepsilon_{a}} (-1)^{[b]} \bigl( f_c E^a_b \otimes 
   q^{\frac{3}{2} h_c} \hat{\sigma}_{ba} - (-1)^{([a]+[b])[c]}  E^a_b f_c 
   \otimes \hat{\sigma}_{ba} q^{-\frac{1}{2} h_c} \bigr).
\end{multline*}

\noindent When $\varepsilon_b - \varepsilon_a - \alpha_c \notin \overline{\Phi}^+$ (recalling that $\overline{\Phi}^+ = \{ \varepsilon_b - \varepsilon_a | \varepsilon_b > \varepsilon_a \}$), this gives

\begin{equation*}
q^{\frac{1}{2} (\alpha_c, \varepsilon_b)} \hat{\sigma}_{ba} f_c - (-1)^{([a]+ 
  [b])[c]} q^{\frac{1}{2} (\alpha_c, \varepsilon_a)} f_c \hat{\sigma}_{ba} = 0.
\end{equation*}

\noindent Conversely, when $\varepsilon_b - \varepsilon_a - \alpha_c = \varepsilon_{b'} - \varepsilon_{a'}$ we obtain

\begin{multline*} 
\underset{\varepsilon_b-\varepsilon_a -\alpha_c = \varepsilon_{b'}-\varepsilon_
  {a'}}{\sum_{\varepsilon_{b'} > \varepsilon_{a'}}} (-1)^{[b']} \bigl( f_c 
  E^{a'}_{b'} \otimes q^{\frac{3}{2} h_c} \hat{\sigma}_{b'\!a'} - 
  (-1)^{([a']+[b'])[c]}  E^{a'}_{b'} f_c \otimes \hat{\sigma}_{b'\!a'} 
  q^{-\frac{1}{2} h_c} \bigr) \\
= (-1)^{[b]} E^a_b \otimes \bigl( q^{\frac{1}{2} (\alpha_c, \varepsilon_b)} 
  \hat{\sigma}_{ba} f_c - (-1)^{([a] + [b])[c]} q^{\frac{1}{2} (\alpha_c, 
  \varepsilon_a)} f_c \hat{\sigma}_{ba} \bigr), \quad \varepsilon_b > 
  \varepsilon_a. 
\end{multline*}

\noindent However $E^a_b$ and $f_c E^{a'}_{b'}$ are linearly independent unless $b=b'$, as are $E^a_b$ and $E^{a'}_{b'}f_c$ when $a \neq a'$. Hence we can simplify this equation to

\begin{multline*}
\underset{\varepsilon_{a'} = \varepsilon_a + \alpha_c}{\sum_{\varepsilon_{b} > 
  \varepsilon_{a'}}} (-1)^{[b]} f_c E^{a'}_{b} \otimes q^{\frac{3}{2} h_c} 
  \hat{\sigma}_{ba'} 
- \underset{\varepsilon_{b'} = \varepsilon_b-\alpha_c}{\sum_{\varepsilon_{b'} 
  > \varepsilon_{a}}} (-1)^{[b']}(-1)^{([a]+[b'])[c]}  E^{a}_{b'} f_c \otimes 
  \hat{\sigma}_{b'\!a} q^{-\frac{1}{2} h_c} \\
= (-1)^{[b]} E^a_b \otimes \bigl( q^{\frac{1}{2} (\alpha_c, \varepsilon_b)} 
  \hat{\sigma}_{ba} f_c - (-1)^{([a]+[b])[c]} q^{\frac{1}{2} (\alpha_c, 
  \varepsilon_a)} f_c \hat{\sigma}_{ba} \bigr), \quad \varepsilon_b > 
  \varepsilon_a.
\end{multline*}

\noindent This then reduces to

\begin{multline*}
(-1)^{[b]} q^{\frac{3}{2} (\varepsilon_b -\varepsilon_a -\alpha_c, \alpha_c)}
  f_c E^{a'}_{b} \otimes \hat{\sigma}_{ba'} q^{\frac{3}{2} h_c} \Big\vert_
  {\varepsilon_{a'} = \varepsilon_a + \alpha_c} \\
- (-1)^{[b]+[c]} (-1)^{([a]+[b]+[c])[c]}  E^{a}_{b'} f_c \otimes \hat{\sigma}_
  {b'a} q^{-\frac{1}{2} h_c} \Big\vert_{\varepsilon_{b'} = \varepsilon_b - 
  \alpha_c} \\
= (-1)^{[b]} E^a_b \otimes \bigl( q^{\frac{1}{2} (\alpha_c, \varepsilon_b)} 
  \hat{\sigma}_{ba} f_c - (-1)^{([a]+[b])[c]} q^{\frac{1}{2} (\alpha_c, 
  \varepsilon_a)} f_c \hat{\sigma}_{ba} \bigr) 
\end{multline*}

\noindent for $\varepsilon_b > \varepsilon_a$, which can, in turn, be simplified to

\begin{multline}
q^{\frac{3}{2} (\varepsilon_b -\varepsilon_a -\alpha_c, \alpha_c)} \langle 
  a|f_c|a' \rangle \hat{\sigma}_{ba'} q^{\frac{3}{2} h_c} - (-1)^{([a]+[b])[c]}
  \langle b'|f_c|b \rangle \hat{\sigma}_{b'a} q^{-\frac{1}{2} h_c} \\
= q^{\frac{1}{2} (\alpha_c, \varepsilon_b)} \hat{\sigma}_{ba} f_c - 
  (-1)^{([a]+[b])[c]} q^{\frac{1}{2} (\alpha_c, \varepsilon_a)} f_c 
  \hat{\sigma}_{ba}, \quad \varepsilon_b > \varepsilon_a. \label{eqf}
\end{multline}

\noindent We now test whether equations (\ref{basef}) and (\ref{eqf}) are consistent with the $\hat{\sigma}_{ba}$ found in the preceding chapter.  Firstly, we use the former to check the fundamental values of $\hat{\sigma}_{ba}$.

\begin{multline*}
f_c \otimes (q^{\frac{3}{2}h_c} - q^{-\frac{1}{2} h_c})\\ 
= (q - q^{-1}) \sum_{\varepsilon_{b} - \varepsilon_{a} = \alpha_c} (-1)^{[b]}  
  E^a_b \otimes \bigl( q^{\frac{1}{2} (\alpha_c, \varepsilon_b)} \hat{\sigma}_
  {ba} f_c - (-1)^{[c]} q^{\frac{1}{2} (\alpha_c, \varepsilon_a)} f_c 
  \hat{\sigma}_{ba} \bigr).
\end{multline*}

\noindent  Consider the case of the root $\alpha_i = \varepsilon_i - \varepsilon_{i+1}$, $1 \leq i < l$, so $f_i \equiv E^{i+1}_{i} - E^{\overline{i}}_{\overline{i+1}}$.  Then the equation becomes:

\begin{align*}
(E^{i+1}_{i} - E_{\overline{i+1}}^{\overline{i}}) \otimes (q^{\frac{3}{2}h_i} 
  - q^{-\frac{1}{2} h_i}) = \quad 
&(q - q^{-1}) E^{i+1}_{i} \otimes (
  q^{\frac{1}{2}} \hat{\sigma}_{i\, i+1} f_i - q^{-\frac{1}{2}} f_i 
  \hat{\sigma}_{i\, i+1}) \notag \\ 
+ \;&(q - q^{-1}) E^{\overline{i}}_{\overline{i+1}} \otimes (q^{\frac{1}{2}}
  \hat{\sigma}_{\overline{i+1}\, \overline{i}} f_i - q^{-\frac{1}{2}} f_i 
  \hat{\sigma}_{\overline{i+1}\, \overline{i}}).
\end{align*}

\noindent Hence we can see immediately that $\hat{\sigma}_{i\, i+1} = - \hat{\sigma}_{\overline{i+1}\, \overline{i}}$, and that

\begin{alignat*}{2}
&& q^{h_i} - q^{- h_i} &= (q - q^{-1}) (q^{\frac{1}{2}} \hat{\sigma}_{i\, i+1} 
  f_i q^{-\frac{1}{2}h_i} - q^{-\frac{1}{2}} f_i \hat{\sigma}_{i\, i+1} 
  q^{-\frac{1}{2} h_i}) \notag \\
&\therefore \qquad& \frac{q^{h_i} - q^{- h_i}}{q - q^{-1}} &= q^{-\frac{1}{2}} 
  \hat{\sigma}_{i\, i+1} q^{-\frac{1}{2}h_i} f_i - q^{-\frac{1}{2}} f_i 
  \hat{\sigma}_{i\, i+1} q^{-\frac{1}{2} h_i} \notag \\
&\therefore& [e_i, f_i] &= [q^{-\frac{1}{2}} \hat{\sigma}_{i\, i+1} 
  q^{-\frac{1}{2}h_i}, f_i].
\end{alignat*}

\noindent  This is certainly consistent with

\begin{equation*}
\hat{\sigma}_{i\, i+1} = - \hat{\sigma}_{\overline{i+1}\,\overline{i}} =
  q^{\frac{1}{2}} e_i q^{\frac{1}{2}h_i},
\end{equation*}

\noindent  the formula obtained in Section \ref{iv}.  Similarly, we can check all the other fundamental values using the same method, and in each case they are consistent with those previously obtained.  Thus it only remains to check that the relations arising out of the equation 

\begin{multline*}
q^{\frac{3}{2} (\varepsilon_b -\varepsilon_a -\alpha_c, \alpha_c)} \langle 
  a|f_c|a' \rangle \hat{\sigma}_{ba'} q^{\frac{3}{2} h_c} - (-1)^{[c]([a]+[b])}
  \langle b'|f_c|b \rangle \hat{\sigma}_{b'a} q^{-\frac{1}{2} h_c} \\
= q^{\frac{1}{2} (\alpha_c, \varepsilon_b)} \hat{\sigma}_{ba} f_c - 
  (-1)^{[c]([a] + [b])} q^{\frac{1}{2} (\alpha_c, \varepsilon_a)} f_c 
  \hat{\sigma}_{ba}, \quad \varepsilon_b > \varepsilon_a
\end{multline*}

\noindent  are consistent with relations \eqref{qcom} and \eqref{indrel} from the previous chapter.

\

\noindent  Again, consider the root $\alpha_i = \varepsilon_i - \varepsilon_{i+1}$, $1 \leq i < l$. Here $f_i = \sigma^{i+1}_{i} \equiv E^{i+1}_{i} - E^{\overline{i}}_{\overline{i+1}}$. Then

\begin{equation*}
\langle a|f_i = \delta^a_{i+1} \langle i| - \delta^a_{\overline{i}} \langle 
  \overline{i+1}|,\qquad f_i |b \rangle = \delta^b_i |i+1 \rangle - 
 \delta^b_{\overline{i+1}} |\overline{i} \rangle.
\end{equation*}

\noindent Hence our equation becomes:

\begin{multline*}
q^{\frac{3}{2} [(\varepsilon_b,\alpha_i)-1]} (\delta^a_{i+1} \hat{\sigma}_{bi}
  - \delta^a_{\overline{i}} \hat{\sigma}_{b\,\overline{i+1}}) q^{\frac{3}{2} 
  h_i} - (\delta^b_i \hat{\sigma}_{i+1\,a} - \delta^b_{\overline{i+1}} 
  \hat{\sigma}_{\overline{i}\,a}) q^{-\frac{1}{2} h_i} \\
= q^{\frac{1}{2} (\alpha_i, \varepsilon_b)} \hat{\sigma}_{ba} f_i -  
  q^{\frac{1}{2} (\alpha_i, \varepsilon_a)} f_i \hat{\sigma}_{ba}, \quad 
  \varepsilon_b > \varepsilon_a.
\end{multline*}

\noindent From this we can deduce the following relations:

\begin{alignat}{2}
&q^{-\frac{3}{2}} \hat{\sigma}_{bi} q^{\frac{3}{2} h_i} = \hat{\sigma}_{b\,i+1}
  f_i - q^{-\frac{1}{2}} f_i \hat{\sigma}_{b\, i+1}, && \varepsilon_b > 
  \varepsilon_i, \label{bi} \\
&q^{-\frac{3}{2}} q^{\frac{3}{2} (\varepsilon_b,\alpha_i)} \hat{\sigma}_{b\,
  \overline{i+1}} q^{\frac{3}{2} h_i} = q^{-\frac{1}{2}} f_i \hat{\sigma}_{b\,
  \overline{i}} - q^{\frac{1}{2} (\alpha_i, \varepsilon_b)} \hat{\sigma}_{b\,
  \overline{i}} f_i, &\qquad& b \neq i, \varepsilon_b> - \varepsilon_{i+1}, 
  \notag \\
&\hat{\sigma}_{i+1\, a} q^{-\frac{1}{2} h_i} = q^{\frac{1}{2} (\alpha_i, 
  \varepsilon_a)} f_i \hat{\sigma}_{ia} - q^{\frac{1}{2}} \hat{\sigma}_{ia} 
  f_i, && a \neq \overline{i}, \varepsilon_a < \varepsilon_{i+1}, \notag \\
&\hat{\sigma}_{\overline{i}\, a} q^{-\frac{1}{2} h_i} = q^{\frac{1}{2}} 
  \hat{\sigma}_{\overline{i+1}\, a} f_i - f_i \hat{\sigma}_{\overline{i+1}\, a}
  , && \varepsilon_a < -\varepsilon_i, \notag \\
&\hat{\sigma}_{i\,\overline{i+1}} q^{\frac{3}{2} h_i} + \hat{\sigma}_{i+1\, 
  \overline{i}} q^{-\frac{1}{2} h_i} = q^{-\frac{1}{2}} f_i \hat{\sigma}_{i\,
  \overline{i}} -q^{\frac{1}{2}} \hat{\sigma}_{i\,\overline{i}} f_i,&&\notag\\
&q^{\frac{1}{2} (\alpha_i, \varepsilon_b)} \hat{\sigma}_{ba} f_i -  
  q^{\frac{1}{2} (\alpha_i, \varepsilon_a)} f_i \hat{\sigma}_{ba} =0, &&
  \varepsilon_b > \varepsilon_a; a \neq i+1, \overline i;\, b \neq i, 
  \overline{i+1}. \label{baf} 
\end{alignat}

\noindent  Unlike the relations obtained in the previous chapter, these cannot be used to inductively construct the $\hat{\sigma}_{ba}$.  Also, there is no simple general form.  Neither of these is a problem, however, since we only need to confirm that these relations are consistent with those in Chapter \ref{R-mat}.  

\

\noindent For instance, consider relation (\ref{bi}).  Previously we found 

\begin{equation*}
\hat{\sigma}_{b\, i+1} = \hat{\sigma}_{bi} \hat{\sigma}_{i\, i+1} - q^{-1} 
\hat{\sigma}_{i\, i+1} \hat{\sigma}_{bi}.
\end{equation*}

\noindent  Using this, we find that

\begin{align*}
RHS &= \hat{\sigma}_{b\,i+1} f_i - q^{-\frac{1}{2}} f_i \hat{\sigma}_{b\,i+1}\\
&= (\hat{\sigma}_{bi} \hat{\sigma}_{i\, i+1} - q^{-1} \hat{\sigma}_{i\, i+1} 
   \hat{\sigma}_{bi}) f_i - q^{-\frac{1}{2}} f_i (\hat{\sigma}_{bi} 
   \hat{\sigma}_{i\, i+1} - q^{-1} \hat{\sigma}_{i\, i+1} \hat{\sigma}_{bi})
   \\
&= q^{\frac{1}{2}} \hat{\sigma}_{bi} e_i q^{\frac{1}{2} h_i} f_i - 
   q^{-\frac{1}{2}} e_i q^{\frac{1}{2} h_i} \hat{\sigma}_{bi} f_i - f_i 
   \hat{\sigma}_{bi} e_i q^{\frac{1}{2} h_i} + q^{-1} f_i e_i q^{\frac{1}{2} 
   h_i} \hat{\sigma}_{bi}.
\end{align*}

\noindent Note from equation (\ref{baf}) that whenever $\varepsilon_b > \varepsilon_i,$

\begin{equation*}
\hat{\sigma}_{bi} f_i = q^{\frac{1}{2}} f_i \hat{\sigma}_{bi}.
\end{equation*}

\noindent Applying this together with the usual commutation relations, we see

\begin{align*}
RHS &= q^{-\frac{1}{2}} \hat{\sigma}_{bi} e_i f_i q^{\frac{1}{2} h_i} - e_i 
  q^{\frac{1}{2} h_i} f_i \hat{\sigma}_{bi} -q^{-\frac{1}{2}} \hat{\sigma}_{bi}
  f_i e_i q^{\frac{1}{2} h_i} \\
& \qquad + q^{-1} \Bigl(e_i f_i - \frac{q^{h_i} -q^{-h_i}}
  {q-q^{-1}} \Bigr) q^{\frac{1}{2} h_i} \hat{\sigma}_{bi} \\
&= q^{-\frac{1}{2}} \hat{\sigma}_{bi} e_i f_i q^{\frac{1}{2} h_i} - q^{-1} e_i
  f_i q^{\frac{1}{2} h_i} \hat{\sigma}_{bi} - q^{-\frac{1}{2}} \hat{\sigma}_
  {bi} \Bigl( e_i f_i - \frac{q^{h_i} - q^{-h_i}}{q-q^{-1}} \Bigr) 
  q^{\frac{1}{2} h_i} \\
& \qquad  + q^{-1} \Bigl( e_i f_i - \frac{q^{h_i} - q^{-h_i}}{q-q^{-1}} \Bigr)
  q^{\frac{1}{2} h_i} \hat{\sigma}_{bi} \\
&= \frac{1}{q-q^{-1}} \Bigl[ q^{-\frac{1}{2}} \hat{\sigma}_{bi} (q^{\frac{3}{2}
  h_i} - q^{-\frac{1}{2}h_i})  - q^{-1} (q^{\frac{3}{2}h_i} - q^{-\frac{1}{2}
  h_i}) \hat{\sigma}_{bi} \Bigr] \\
&= \frac{1}{q-q^{-1}} \Bigl[ q^{-\frac{1}{2}} \hat{\sigma}_{bi} (q^{\frac{3}{2}
  h_i} - q^{-\frac{1}{2}h_i})  - q^{-1} \hat{\sigma}_{bi} (q^{\frac{3}{2} 
  (\alpha_i, \varepsilon_b - \varepsilon_i)} q^{\frac{3}{2}h_i} - 
  q^{-\frac{1}{2} (\alpha_i, \varepsilon_b - \varepsilon_i)} q^{-\frac{1}{2}
  h_i}) \Bigr] \\
&= \frac{q^{-\frac{1}{2}} - q^{-\frac{5}{2}}}{q-q^{-1}} \hat{\sigma}_{bi}
  q^{\frac{3}{2}h_i} \\
&= q^{-\frac{3}{2}} \hat{\sigma}_{bi} q^{\frac{3}{2}h_i}
\end{align*}

\noindent as expected.  Hence equation (\ref{bi}) is consistent with the defining relations for $\hat{\sigma}_{ba}$ found in the previous chapter. 

\
 
\noindent Although time-consuming, it can be confirmed that all the other relations generated by equation (\ref{eqf}) are similarly consistent, regardless of which root is chosen.  Thus we have verified that the matrix $R$ constructed in the previous chapter satisfies the intertwining property 

\begin{equation*}
R \Delta(a) = \Delta^T(a) R
\end{equation*}

\noindent for all elements $a \in U_q[osp(m|n)]$.

\section{The Lax Operator}

\noindent We have now proven, as expected, that the matrix $R$ found in the previous chapter satisfies both the intertwining property and $(\text{id}\, \otimes\, \Delta)R = R_{13} R_{12}.$  The other $R$-matrix property, containing $(\Delta\, \otimes\, \text{id})R$, is clearly not applicable here.  It is not necessary, however, as we know there is a Lax operator belonging to $\pi \bigl( U_q[osp(m|n)]^-\bigr) \otimes U_q[osp(m|n)]^+$, and we have shown there is only one such possibility.  Thus the work in this chapter confirms the following theorem:

\begin{theorem}

The Lax operator, $R = (\pi \otimes \text{id}) \mathcal{R}$ for the quantum superalgebra $U_q[osp(m|n)]$, where $\mathcal{R} \in U_q[osp(m|n)]^- \otimes U_q[osp(m|n)]^+$ and $m >2$, is given by

\begin{align*}
R &= q ^ {h_{x} \otimes h^{x}} \Bigl[ I \otimes I + (q - q^{-1}) \sum_{\varepsilon_{a} < 
  \varepsilon_{b}} (-1)^{[b]} E^a_b \otimes \hat{\sigma}_{ba} \Bigr] \\
&= \sum_a E^a_a \otimes q^{h_{\varepsilon_a}} + (q-q^{-1}) \sum_{\varepsilon_a
  < \varepsilon_b} (-1)^{[b]} E^a_b \otimes q^{h_{\varepsilon_a}} 
  \hat{\sigma}_{ba},
\end{align*}

\noindent where the operators $\hat{\sigma}_{ba}$ satisfy:

\

\noindent (i) the $q$-commutation relations

\begin{equation*}
q^{(\alpha_c, \varepsilon_b)} \hat{\sigma}_{ba} e_c q^{\frac{1}{2} h_c} 
      - (-1)^{([a]+[b])[c]} q^{-(\alpha_c, \varepsilon_a)} e_c 
      q^{\frac{1}{2}h_c} \hat{\sigma}_{ba} = 0, \quad \varepsilon_b > 
      \varepsilon_a
\end{equation*}

\noindent when neither $\varepsilon_a - \alpha_c$ nor $\varepsilon_b + \alpha_c$ equals any $\varepsilon_x$; and 

\

\noindent (ii) the recursion relations

\begin{equation*}
\hat{\sigma}_{ba} = q^{-(\varepsilon_b,\varepsilon_a)} \hat{\sigma}_{bc}
  \hat{\sigma}_{ca} - q^{-(\varepsilon_c, \varepsilon_c)} (-1)^{([b]+[c])
  ([a]+[c])} \hat{\sigma}_{ca} \hat{\sigma}_{bc}, \quad \varepsilon_b
  > \varepsilon_c > \varepsilon_a
\end{equation*}

\noindent when $c \neq \overline{b}$ or $\overline{a}$; and with initial values given by:

\begin{alignat*}{4}
&\hat{\sigma}_{i\, i+1}&&= - \hat{\sigma}_{\overline{i+1}\, \overline{i}}&&= 
  q^{\frac{1}{2}} e_i q^{\frac{1}{2} h_i},& \quad& 1 \leq i < l, \\
&\hat{\sigma}_{l-1\, \overline{l}}&& = - \hat{\sigma}_{l\, \overline{l-1}}&&= 
  q^{\frac{1}{2}} e_l q^{\frac{1}{2} h_l}, && m=2l, \\
&\hat{\sigma}_{l\, l+1}&& = - q^{-\frac{1}{2}} \hat{\sigma}_{l+1\,\overline{l}}
  && = e_l q^{\frac{1}{2} h_l}, &&m=2l+1, \\
&\hat{\sigma}_{\mu \,\mu+1} &&= \hat{\sigma}_{\overline{\mu+1}\,\overline{\mu}}
  &&= q^{-\frac{1}{2}} e_\mu q^{\frac{1}{2} h_\mu}, && 1 \leq \mu < k, \\
&\hat{\sigma}_{\mu=k\,i=1} &&= (-1)^k q\, \hat{\sigma}_{\overline{i} = 
  \overline{1} \, \overline{\mu} = \overline{k}} &&= q^{\frac{1}{2}} e_s 
  q^{\frac{1}{2} h_s},\\
&\hat{\sigma}_{l\, \overline{l}} &&= 0, &&&& m=2l.
\end{alignat*}

\end{theorem}

\

\noindent As an aside, the two properties verified directly are sufficient to prove $R$ satisfies the Yang-Baxter equation.  For using only those, we see

\begin{align*}
R_{23} R_{13} R_{12} &= R_{23} (\text{id} \otimes \Delta)R \\
&= [(\text{id} \otimes \Delta ^T )R] R_{23} \\
&= [(\text{id} \otimes T) ((\text{id} \otimes \Delta)R] R_{23} \\
&= [(\text{id} \otimes T) R_{13} R_{12}] R_{23} \\
&= R_{12} R_{13} R_{23}
\end{align*}

\noindent as required. 

\

\noindent  It is very surprising that there is a unique solution to 

\begin{equation*}
R \Delta (e_c) = \Delta^T (e_c) R,
\end{equation*}

\noindent even given we restricted ourselves to matrices in $\pi \bigl( U_q[osp(m|n)]^-\bigr) \otimes U_q[osp(m|n)]^+$.  While it is reassuring that the solution is a Lax operator, it means the remaining $R$-matrix relations were redundant, which raises the question of why.  It suggests there may be some underlying symmetries in the system; some way in which the other $R$-matrix properties can be derived from the one used.  If so, however, they are not obvious.

\section{The Opposite Lax Operator} \label{opposite}
    
\noindent Having found the Lax operator $R = (\pi \otimes \text{id}) \mathcal{R}$, we wish to use that result to find its opposite $R^T = (\pi \otimes \text{id}) \mathcal{R}^T$, where $\mathcal{R}^T$ is the opposite universal $R$-matrix of $U_q[osp(m|n)]$.  We begin by showing that $\mathcal{R}^T$ is in fact equal to $\mathcal{R}^\dagger$, where $^\dagger$ represents graded conjugation, defined below.

\

\noindent A graded conjugation on $U_q[osp(m|n)]$ is defined on the simple generators by:

\begin{equation*}
e_a^\dagger = f_a, \qquad
f_a^\dagger = (-1)^{[a]} e_a, \qquad
h_a^\dagger = h_a.
\end{equation*}

\noindent It is consistent with the coproduct and extends naturally to all remaining elements of $U_q[osp(m|n)]$, satisfying the properties:

\begin{align*}
&(\sigma^a_b)^\dagger = (-1)^{[a]([a]+[b])} \sigma^b_a,\\
&(ab)^\dagger = (-1)^{[a][b]} b^\dagger a^\dagger, \\
&(a \otimes b)^\dagger = a^\dagger \otimes b^\dagger, \\
&\Delta (a)^\dagger = \Delta (a^\dagger).
\end{align*}

\noindent  Returning to the universal $R$-matrix $\mathcal{R}$, we know

\begin{alignat*}{3}
&&&\mathcal{R} \Delta(a) = \Delta^T(a) \mathcal{R}, &\qquad &\forall a \in U_q
[osp(m|n)], \notag \\
&\Rightarrow& \quad & \Delta(a)^\dagger \mathcal{R}^\dagger = \mathcal{R}^
  \dagger\Delta^T(a)^\dagger \notag \\
&\Rightarrow && \Delta(a^\dagger) \mathcal{R}^\dagger = \mathcal{R}^\dagger
  \Delta^T(a^\dagger) \notag \\
&\Rightarrow && \Delta(a) \mathcal{R}^\dagger = \mathcal{R}^\dagger
  \Delta^T(a),&& \forall a \in U_q[osp(m|n)].
\end{alignat*}

\noindent Similarly, $\mathcal{R}^\dagger$ satisfies the other $R$-matrix properties \eqref{Requations}.  As there is a unique universal $R$-matrix belonging to $U_q[osp(m|n)]^+ \otimes U_q[osp(m|n)]^-$, the only possibility is $\mathcal{R}^T = \mathcal{R}^\dagger$.

\

\noindent  Now it is known that the vector representation is superunitary.  A discussion of superunitary representations is given in \cite{LG}, where they are called grade star representations, but for this thesis we need only note this implies

\begin{equation*}
\pi (a^\dagger) = \pi (a)^\dagger, \qquad \forall a \in U_q[osp(m|n)].
\end{equation*}

\noindent Hence

\begin{align*}
R^T &= (\pi \otimes \text{id}) \mathcal{R}^\dagger \notag \\
&= [(\pi \otimes \text{id}) \mathcal{R}]^\dagger \notag \\
&= R^\dagger.
\end{align*}

\noindent Thus we can find the opposite Lax operator $R^T$ simply by using the usual rules for graded conjugation.  As $R$ is given by

\begin{equation*}
R = \sum_a E^a_a \otimes q^{h_{\varepsilon_a}} + (q-q^{-1}) \sum_{\varepsilon_b > \varepsilon_a} (-1)^{[b]} E^a_b \otimes q^{h_{\varepsilon_a}} \hat{\sigma}_{ba},
\end{equation*}

\noindent we obtain

\begin{equation*}
R^T = \sum_a E^a_a \otimes q^{h_{\varepsilon_a}} + (q-q^{-1}) \sum_{\varepsilon_b > \varepsilon_a} (-1)^{[b]} (E^a_b)^\dagger \otimes (\hat{\sigma}_{ba})^\dagger q^{h_{\varepsilon_a}}.
\end{equation*}

\noindent  As $(E^a_b)^\dagger = (-1)^{[a]([a]+[b])} E^b_a$, set

\begin{equation*}
\hat{\sigma}_{ab} = (-1)^{[b]([a]+[b])} \hat{\sigma}_{ba}^\dagger, \quad \varepsilon_b > \varepsilon_a.
\end{equation*}

\noindent Then the opposite Lax operator $R^T$ can be written as

\begin{equation} \label{RT}
R^T = \sum_a E^a_a \otimes q^{h_{\varepsilon_a}} + (q-q^{-1}) \sum_{\varepsilon_b > \varepsilon_a} (-1)^{[a]} E^b_a \otimes \hat{\sigma}_{ab} q^{h_{\varepsilon_a}},
\end{equation}

\noindent where the operators $\hat{\sigma}_{ab}$ can be calculated from $\hat{\sigma}_{ba}$ using the usual graded conjugation rules.

\section{q-Serre Relations}
    \noindent Having shown that the relations found in Chapter \ref{R-mat} define a Lax operator, we also wish to see if they incorporate the $q$-Serre relations.  It is too time-consuming to verify all of these, so we will merely provide a couple of examples, including the extra $q$-Serre relations.

\

\noindent First recall that if $\varepsilon_b - \varepsilon_a$ is a simple root, then $\hat{\sigma}_{ba} \propto e_c q^{\frac{1}{2}h_c}$ for either $c = b$ or $c=\overline{a}$.  Then setting $E_a = e_a q^{\frac{1}{2} h_a}$, we see from the definitions on page \pageref{coprod} that:

\begin{alignat}{2}
&&\Delta (E_a) &= q^{h_a} \otimes E_a + E_a \otimes 1 \notag \\
&&S(E_a) &= -q^{-\frac{1}{2}(\alpha_a, \alpha_a)} q^{-\frac{1}{2}h_a}e_a \notag
  \\
&&&= -  q^{-h_a} E_a \notag \\
&\therefore \qquad & ad\, E_a \circ b &= - (-1)^{[a][b]} q^{h_a}b 
  q^{-h_a}E_a + E_a b \notag \\
&&&= E_a b - (-1)^{[a][b]} q^{(\alpha_a,\varepsilon_b)} bE_a. \label{adjoint}
\end{alignat}

\noindent Now consider the simple generators $\hat{\sigma}_{i\,i+1}$ and $\hat{\sigma}_{i+1\, i+2}$.

\begin{alignat*}{2}
(ad\: \hat{\sigma}_{i\,i+1}\: \circ)^2 \hat{\sigma}_{i+1\, i+2} &= ad\: \hat
  {\sigma}_{i\,i+1} \circ (\hat{\sigma}_{i\,i+1} \hat{\sigma}_{i+1\, i+2} - 
  q^{-1} \hat{\sigma}_{i+1\, i+2} \hat{\sigma}_{i\,i+1}) &&\notag \\
&= ad\: \hat{\sigma}_{i\,i+1} \circ \hat{\sigma}_{i\,i+2} &&\notag \\
&= \hat{\sigma}_{i\,i+1} \hat{\sigma}_{i\,i+2} - q  \hat{\sigma}_{i\,i+2}
   \hat{\sigma}_{i\,i+1} &&\notag \\
&= 0&& \text{from } \eqref{qcom}.
\end{alignat*}

\noindent This is equivalent to the $q$-Serre relation $(ad\: e_b\: \circ)^{1-a_{bc}} e_c = 0$ for this pair of simple operators.  In a similar way, we can verify this relation for any $b \neq c$.  The defining relations for the $\hat{\sigma}_{ba}$, therefore, incorporates all the standard $q$-Serre relations for raising generators.

\

\noindent This still leaves the extra $q$-Serre relations, which involve the odd root.  There are only two of these for our choice of simple roots \cite{Yamane}.  Explicitly, taking into account the different conventions, the relevant extra $q$-Serre relations for $U_q[osp(m|n)]$ can be written as  

\begin{align} \label{q1}
&\bigl[ \hat{\sigma}_{\mu=k\,i=1}, \bigl[ \hat{\sigma}_{\nu=k-1\,\mu=k}, [\hat{\sigma}_{\mu=k\,i=1},\hat{\sigma}_{i=1\,j=2} ]_q \, ]_q \,] = 0 \\
&\bigl[ \hat{\sigma}_{\mu=k\,i=1}, \bigl[ \hat{\sigma}_{i=1\,j=2}, [\hat{\sigma}_{\mu=k\,i=1}, \hat{\sigma}_{\nu=k-1\,\mu=k}]_q \, ]_q \,] = 0, \label{q2}
 \end{align}

\noindent where $[x,y]_q$ represents the adjoint action $ad\, x \circ y$.  

\

\noindent Consider equation \eqref{q1}.  Using the defining relations \eqref{qcom} and \eqref{indrel} for the $\hat{\sigma}_{ba}$ together with the adjoint action as given in equation \eqref{adjoint}, we find:

\begin{align*}
&[\hat{\sigma}_{\mu=k\,i=1}, [ \hat{\sigma}_{\nu=k-1\,\mu=k}, 
  [\hat{\sigma}_{\mu=k\,i=1},\hat{\sigma}_{i=1\,j=2} ]_q \,  ]_q \,]\notag\\
=\:&[ \hat{\sigma}_{\mu=k\,i=1}, [ \hat{\sigma}_{\nu=k-1\,\mu=k},
  (\hat{\sigma}_{\mu=k\,i=1} \hat{\sigma}_{i=1\,j=2} - q^{-1} 
  \hat{\sigma}_{i=1\,j=2} \hat{\sigma}_{\mu=k\,i=1}) ]_q \,]\notag\\
=\:&[ \hat{\sigma}_{\mu=k\,i=1}, [ \hat{\sigma}_{\nu=k-1\,\mu=k},
  \hat{\sigma}_{\mu=k\, j=2} ]_q \,] \notag \\
=\:&[ \hat{\sigma}_{\mu=k\,i=1}, (\hat{\sigma}_{\nu=k-1\,\mu=k} 
  \hat{\sigma}_{\mu=k\, j=2} - q \hat{\sigma}_{\mu=k\, j=2} 
  \hat{\sigma}_{\nu=k-1\,\mu=k})] \notag \\
=\:&[ \hat{\sigma}_{\mu=k\,i=1}, \hat{\sigma}_{\nu=k-1\,j=2}] \notag \\
=\:& 0
\end{align*}

\noindent as required.  It is equally straightforward to show that equation \eqref{q2} arises from the defining relations of the $\hat{\sigma}_{ba}$.  Hence these compact defining relations for the $\hat{\sigma}_{ba}$ incorporate not only the standard $q$-Serre relations for the raising generators, but also the extra ones.  This is quite interesting, as the equivalent $q$-Serre relations were not used in the derivation.

\chapter{The $R$-matrix for the Vector Representation} \label{vector}
        \noindent The Lax operator can be used to explicitly calculate an $R$-matrix for any representation stemming from the $\pi \otimes \text{id}$ representation.  In particular, it provides a more straightforward method of calculating $R$ for the tensor product of the vector representation, $\pi \otimes \pi$, than previously found \cite{Mehta}.

\

\noindent By specifically constructing the $R$-matrix for the vector representation, we also illustrate concretely the way the recursion relations can be applied to find the $R$-matrix for an arbitrary representation.  Although the values for $\hat{\sigma}_{ba}$ obtained will change for each representation, they can always be constructed by applying the same equations in the same order.  We could choose to use only the relations listed in the tables in the appendix, but using the general form of the inductive relations shortens and simplifies the process.


\section{Fundamental values of $\hat{\sigma}_{ba}$} \label{base}

\noindent  The first step is to calculate the values of $\hat{\sigma}_{ba}$ where $\varepsilon_b - \varepsilon_a$ is a simple root, using the formulae derived in Section \ref{iv}.  As before, we use $e_c$ and $h_c$ to denote the image of the raising generators and Cartan elements in the vector representation, with the $\pi$ being implicit.

\

\noindent Now recall that in the vector representation $e_i = E^i_{i+1} - E^{\overline{i+1}}_{\overline{i}}$ and $h_i = E^i_i - E^{i+1}_{i+1} + E^{\overline{i+1}}_{\overline{i+1}}- E^{\overline{i}}_{\overline{i}}$.  Then for $i<l$

\vspace{-2mm}

\begin{align*}
\hat{\sigma}_{i\, i+1} = -\hat{\sigma}_{\overline{i+1}\, \overline{i}} 
&= q^{\frac{1}{2}} e_i q^{\frac{1}{2} h_i} \notag \\
&= q^{\frac{1}{2}}  (E^i_{i+1} - E^{\overline{i+1}}_{\overline{i}}) 
   [I + (q^{\frac{1}{2}}-1) (E^i_i + E^{\overline{i+1}}_{\overline{i+1}}) + 
   (q^{-\frac{1}{2}}-1) (E^{i+1}_{i+1} + E^{\overline{i}}_{\overline{i}})] 
   \notag\\
&= E^i_{i+1} - E^{\overline{i+1}}_{\overline{i}}.
\end{align*}

\

\noindent In the case $m=2l$, we have $e_l = E^{l-1}_{\overline{l}} - E^l_{\overline{l-1}}$ and $h_l = E^l_l - E^{\overline{l}}_{\overline{l}} + E^{l-1}_{l-1} - E^{\overline{l-1}}_{\overline{l-1}}$.  So

\begin{align*}
\hat{\sigma}_{l-1\, \overline{l}} = -\hat{\sigma}_{l\, \overline{l-1}} 
&= q^{\frac{1}{2}} e_l q^{\frac{1}{2} h_l} \notag \\
&= q^{\frac{1}{2}} (E^{l-1}_{\overline{l}} - E^l_{\overline{l-1}})  
   q^{\frac{1}{2} h_l} \notag \\
&=  E^{l-1}_{\overline{l}} - E^l_{\overline{l-1}}.
\end{align*}

\

\noindent When $m=2l+1$, $e_l = E^l_{l+1} - E^{l+1}_{\overline{l}}$, while $h_l = E^l_l - E^{\overline{l}}_{\overline{l}}$.  Thus

\begin{align*}
\hat{\sigma}_{l\, l+1} = - q^{-\frac{1}{2}} \hat{\sigma}_{l+1\, \overline{l}}
&= e_l q^{\frac{1}{2} h_l} \notag \\
&= (E^l_{l+1} - E^{l+1}_{\overline{l}}) q^{\frac{1}{2} h_l} \notag \\
&= E^l_{l+1} - q^{-\frac{1}{2}} E^{l+1}_{\overline{l}}.
\end{align*}

\noindent  Similarly, $e_\mu = E^\mu_{\mu+1} + E^{\overline{\mu+1}}_{\overline{\mu}}$ and $h_\mu = E^{\mu+1}_{\mu+1} + E^{\overline{\mu}}_{\overline{\mu}} - E^\mu_\mu - E^{\overline{\mu+1}}_{\overline{\mu+1}}$. These give:

\begin{align*}
\hat{\sigma}_{\mu\, \mu+1} = \hat{\sigma}_{\overline{\mu+1}\, \overline{\mu}}
&= q^{-\frac{1}{2}} e_\mu q^{\frac{1}{2} h_\mu} \notag \\
&= q^{-\frac{1}{2}} (E^\mu_{\mu+1} + E^{\overline{\mu+1}}_{\overline{\mu}}) 
   q^{\frac{1}{2} h_\mu} \notag \\
&= E^\mu_{\mu+1} + E^{\overline{\mu+1}}_{\overline{\mu}}.
\end{align*}

\

\noindent  Lastly, in the case of the odd root remember that $e_s = E^{\mu=k}_{i=1} + (-1)^k E^{\overline{i} = \overline{1}}_{\overline{\mu} = \overline{k}}$, whereas \\
$h_s = E^{\overline{i}=\overline{1}}_{\overline{i}=\overline{1}} - E^{i=1}_{i=1} + E^{\overline{\mu} = \overline{k}}_{\overline{\mu} = \overline{k}} - E^{\mu=k}_{\mu=k}$.  Applying these, we find

\begin{align*}
\hat{\sigma}_{\mu=k\,i=1} = (-1)^kq\, \hat{\sigma}_{\overline{i}=\overline{1}\,
   \overline{\mu} = \overline{k}} &= q^{\frac{1}{2}} e_s q^{\frac{1}{2} h_s} 
   \notag \\
&= q^{\frac{1}{2}} (E^{\mu=k}_{i=1} + (-1)^k E^{\overline{i} = \overline{1}}_
   {\overline{\mu} = \overline{k}}) q^{\frac{1}{2} h_s}  \notag \\
&=E^{\mu=k}_{i=1} + (-1)^k q E^{\overline{i} = \overline{1}}_{\overline{\mu} = 
   \overline{k}}.
\end{align*}

\

\noindent This completes the calculation of the fundamental values of $\hat{\sigma}_{ba}$ in the vector representation. They are summarised in Table \ref{basevector}.

\

\begin{table}[ht] \label{basevector}
\caption{The fundamental values for $\hat{\sigma}_{ba}$ in the vector representation.}
\begin{tabular}{|l|lll|} \hline
\multicolumn{1}{|c|}{Simple Root}& \multicolumn{3}{c|}{Corresponding 
  $\hat{\sigma}_{ba}$} \\ \hline
$\alpha_i = \varepsilon_i - \varepsilon_{i+1},\,i < l$ & $\hat{\sigma}_{i\,
  i+1}$&$ = -\hat{\sigma}_{\overline{i+1}\,\overline{i}}$&$ = E^i_{i+1} - E^
  {\overline{i+1}}_{\overline{i}}$ \\ 
$\alpha_l = \varepsilon_{l-1} + \varepsilon_l,\,m = 2l$ & $\hat{\sigma}_
  {l-1\,\overline{l}}$&$ = -\hat{\sigma}_{l\, \overline{l-1}}$&$ = E^{l-1}_
  {\overline{l}} - E^l_{\overline{l-1}}$ \\
$\alpha_l = \varepsilon_l,\,m=2l+1$ & $\hat{\sigma}_{l\, l+1}$&$ = -q^{-\frac
  {1}{2}} \hat{\sigma}_{l+1\, \overline{l}}$&$ = E^l_{l+1} - q^{-\frac{1}{2}} 
  E^{l+1}_{\overline{l}}$ \\
$\alpha_\mu = \delta_\mu - \delta_{\mu+1},\,\mu < k$ & $\hat{\sigma}_
  {\mu\, \mu+1} $&$ = \hat{\sigma}_{\overline{\mu+1}\, \overline{\mu}} $&$= 
  E^\mu_{\mu+1} + E^{\overline{\mu+1}}_{\overline{\mu}}$ \\
$\alpha_s = \delta_k - \varepsilon_1,$ & $\hat{\sigma}_{\mu=k\, i=1} $&$ = 
  (-1)^k q \hat{\sigma}_{\overline{i}=\overline{1}\, \overline{\mu} = 
  \overline{k}} $&$=  E^{\mu=k}_{i=1} + (-1)^k q E^{\overline{i} = 
  \overline{1}}_{\overline{\mu} = \overline{k}}$ \\ \hline
\end{tabular}
\end{table}


\section{Calculating $\hat{\sigma}_{ji},\hat{\sigma}_{\overline{i}\, \overline{j}}$} \label{ij}

\noindent Now that the fundamental values of $\hat{\sigma}_{ba}$ for the vector representation have been explicitly calculated, the remaining values can be found by applying the various inductive relations.  There is no one correct way of doing this, with several equivalent methods giving the same result.  We choose to begin by finding the remaining operators of the form $\hat{\sigma}_{ji}$ and $\hat{\sigma}_{\overline{i}\, \overline{j}}$.  As mentioned earlier, the same process can be applied to any representation.

\

\noindent In the previous section we found that $\hat{\sigma}_{i\, i+1} =  E^i_{i+1} - E^{\overline{i+1}}_{\overline{i}}$ for $i < l$.  We also know from Chapter \ref{R-mat} that 

\begin{equation*}
\hat{\sigma}_{b\,i+1} = \hat{\sigma}_{bi} \hat{\sigma}_{i\, i+1} - q^{-1} 
  \hat{\sigma}_{i\, i+1} \hat{\sigma}_{bi}, \qquad  i<l, \; \varepsilon_b > 
  \varepsilon_i.
\end{equation*}

\noindent Combining these, we find

\begin{alignat*}{2}
\hat{\sigma}_{i-1\, i+1} &= \hat{\sigma}_{i-1\, i} \hat{\sigma}_{i\, i+1} - 
  q^{-1} \hat{\sigma}_{i\, i+1} \hat{\sigma}_{i-1\, i} && \\
&= (E^{i-1}_{i} - E^{\overline{i}}_{\overline{i-1}}) (E^i_{i+1} - 
  E^{\overline{i+1}}_{\overline{i}}) - q^{-1} (E^i_{i+1} - E^{\overline{i+1}}_
  {\overline{i}}) (E^{i-1}_{i} - E^{\overline{i}}_{\overline{i-1}})&& \\
&= E^{i-1}_{i+1} - q^{-1} E^{\overline{i+1}}_{\overline{i-1}},&&1<i<l
\end{alignat*}

\noindent and

\begin{alignat*}{2}
\hat{\sigma}_{i-2\, i+1} &= \hat{\sigma}_{i-2\, i} \hat{\sigma}_{i\, i+1} -
  q^{-1} \hat{\sigma}_{i\, i+1} \hat{\sigma}_{i-2\, i}&& \\
&= (E^{i-2}_{i} - q^{-1} E^{\overline{i}}_{\overline{i-2}}) (E^i_{i+1} - 
  E^{\overline{i+1}}_{\overline{i}}) - q^{-1} (E^i_{i+1} - E^{\overline{i+1}}_
  {\overline{i}}) (E^{i-2}_{i} - q^{-1}&& E^{\overline{i}}_{\overline{i-2}})\\
&= E^{i-2}_{i+1} - q^{-2} E^{\overline{i+1}}_{\overline{i-2}},&&2<i<l.
\end{alignat*}

\noindent We postulate that $\hat{\sigma}_{ji} = E^j_i - q^{j-i+1} E^{\overline{i}}_{\overline{j}}$, where $1 \leq j < i \leq l$.  Clearly this is true when $i-j \leq 3$.  So we assume it is true for $i-j = x$ for some $x \geq 1$, and try to show it holds for $i-j = x+1$:

\begin{alignat*}{2}
\hat{\sigma}_{ji} &= \hat{\sigma}_{j\, i-1} \hat{\sigma}_{i-1\, i} -   q^{-1} 
  \hat{\sigma}_{i-1\, i} \hat{\sigma}_{ji}, && i-j =x+1 \notag \\
&= (E^j_{i-1} - q^{j-i+2} E^{\overline{i-1}}_{\overline{j}}) (E^{i-1}_i - 
  E^{\overline{i}}_{\overline{i-1}}) \notag \\
& \qquad - q^{-1} (E^{i-1}_i - E^{\overline{i}}_{\overline{i-1}}) (E^j_{i-1} - 
  q^{j-i+2} E^{\overline{i-1}}_{\overline{j}}) & \quad &\text{by inductive 
  hypothesis} \notag \\
&= E^j_i - q^{j-i+1} E^{\overline{i}}_{\overline{j}}, &&
\end{alignat*}

\noindent as required. 

\

\noindent Recalling that $\rho = \frac{1}{2} \sum_{i=1}^l (m-2i) \varepsilon_i + \frac{1}{2} \sum_{\mu=1}^k (n-m+2-2\mu) \delta_\mu$, we write this as

\begin{equation*}
\hat{\sigma}_{ji} = E^j_i - q^{(\rho, \varepsilon_i - \varepsilon_j)+1} 
E^{\overline{i}}_{\overline{j}}, \quad 1 \leq j < i \leq l.
\end{equation*}

\

\noindent This includes all the values for $\hat{\sigma}_{ji}$ except for $i=l+1$ in the case $m=2l+1$.  To find these remaining values, recall that $\hat{\sigma}_{l\, l+1} = E^l_{l+1} - q^{-\frac{1}{2}} E^{l+1}_{\overline{l}}$, and that 

\begin{equation*}
\hat{\sigma}_{i\, l+1} = \hat{\sigma}_{il} \hat{\sigma}_{l\, l+1} - q^{-1} \hat{\sigma}_{l\, l+1} \hat{\sigma}_{il}, \qquad i<l.
\end{equation*}

\noindent Therefore

\begin{align*}
\hat{\sigma}_{j\, l+1} &= \hat{\sigma}_{jl} \hat{\sigma}_{l\, l+1} - q^{-1} 
  \hat{\sigma}_{l\, l+1} \hat{\sigma}_{jl} \notag \\
&= (E^j_l - q^{j-l+1} E^{\overline{l}}_{\overline{j}}) (E^l_{l+1} - 
  q^{-\frac{1}{2}} E^{l+1}_{\overline{l}}) - q^{-1} (E^l_{l+1} - 
  q^{-\frac{1}{2}} E^{l+1}_{\overline{l}}) (E^j_l - q^{j-l+1} 
  E^{\overline{l}}_{\overline{j}}) \notag \\
&= E^j_{l+1} - q^{j-l-\frac{1}{2}}  E^{l+1}_{\overline{j}} \\
&= E^j_{l+1} - q^{(\rho, \varepsilon_{l+1} - \varepsilon_j)} E^{l+1}_
   {\overline{j}} .
\end{align*}

\noindent Unifying this with the previous result, we have shown

\begin{equation} \label{sigji}
\boxed{\hat{\sigma}_{ji} = E^j_i - q^{(\rho, \varepsilon_i - \varepsilon_j)
  +(\varepsilon_i, \varepsilon_i)}  E^{i}_{\overline{j}}, \quad 1 \leq j <i
  \leq \lceil \tfrac{m}{2} \rceil.}
\end{equation}

\

\noindent Similarly, we know that $\hat{\sigma}_{\overline{i+1}\, \overline{i}} =  E^{\overline{i+1}}_{\overline{i}}- E^i_{i+1}$ for $i<l$, and that

\begin{equation*} 
\hat{\sigma}_{\overline{i+1}\, a} = \hat{\sigma}_{\overline{i+1}\,\overline{i}}
  \hat{\sigma}_{\overline{i}\,a} - q^{-1} \hat{\sigma}_{\overline{i}\,a}
  \hat{\sigma}_{\overline{i+1}\,\overline{i}}, \qquad i<l, \,\varepsilon_a < - 
  \varepsilon_i.
\end{equation*}

\noindent  Thus

\begin{alignat*}{2}
\hat{\sigma}_{\overline{i+1}\, \overline{i-1}}  &= \hat{\sigma}_{\overline{i+1}
  \,\overline{i}} \hat{\sigma}_{\overline{i}\,\overline{i-1}} - q^{-1} 
  \hat{\sigma}_{\overline{i}\, \overline{i-1}} \hat{\sigma}_{\overline{i+1}\,
  \overline{i}} && \notag \\
&= (E^{\overline{i+1}}_{\overline{i}}- E^i_{i+1}) (E^{\overline{i}}_
  {\overline{i-1}}- E^{i-1}_i) - q^{-1} (E^{\overline{i}}_{\overline{i-1}} - 
  E^{i-1}_i) (E^{\overline{i+1}}_{\overline{i}}- E^i_{i+1})&& \notag \\ 
&= E^{\overline{i+1}}_{\overline{i-1}} - q^{-1} E^{i-1}_{i+1},&&1 < i < l. 
\end{alignat*}

\noindent This time we try $\hat{\sigma}_{\overline{i}\, \overline{j}} = E^{\overline{i}}_{\overline{j}} - q^{j-i+1} E^j_i$.  Clearly the initial case of $i=j+1$ is satisfied, so assume that it holds for some $\sigma_{\overline{i}\, \overline{j}},\, j<i<l$, and we will show that hence it is true for $\sigma_{\overline{i+1}\, \overline{j}}$:

\begin{align*}
\sigma_{\overline{i+1}\, \overline{j}} &= \hat{\sigma}_{\overline{i+1}\, 
  \overline{i}} \hat{\sigma}_{\overline{i} \overline{j}} - q^{-1} 
  \hat{\sigma}_{\overline{i} \overline{j}} \hat{\sigma}_{\overline{i+1}\,
  \overline{i}} \\
&= (E^{\overline{i+1}}_{\overline{i}} - E^i_{i+1}) (E^{\overline{i}}_
  {\overline{j}}- q^{j-i+1} E^j_i) - q^{-1} (E^{\overline{i}}_{\overline{j}}- 
  q^{j-i+1} E^j_i) (E^{\overline{i+1}}_{\overline{i}} - E^i_{i+1}) \\
&= E^{\overline{i+1}}_{\overline{j}} - q^{j-(i+1)+1} E^{j}_{i+1}.
\end{align*}

\noindent Thus we have inductively found $\hat{\sigma}_{\overline{i}\, \overline{j}}$ for all $1 \leq j < i \leq l$ in the vector representation.  Writing it in terms of the half-sum of positive roots, we have

\begin{equation*}
\hat{\sigma}_{\overline{i}\, \overline{j}} = E^{\overline{i}}
  _{\overline{j}} - q^{(\rho,\varepsilon_{\overline{j}} - 
  \varepsilon_{\overline{i}}) +1}  E^j_i, \quad 1 \leq j < i \leq l.
\end{equation*}

\

\noindent When $m=2l+1$, we must also find $\hat{\sigma}_{\overline{l+1}\, \overline{j}} \;(= \hat{\sigma}_{l+1\, \overline{j}})$ for $1 \leq j \leq l$. From Table \ref{basevector} we have $\hat{\sigma}_{l+1\, \overline{l}} = E^{l+1}_{\overline{l}} - q^{\frac{1}{2}} E^l_{l+1}$.  Using this, we see that for all $j<l$

\begin{align*}
\hat{\sigma}_{l+1\, \overline{j}} &= \hat{\sigma}_{l+1\, \overline{l}} 
  \hat{\sigma}_{\overline{l}\,\overline{j}} - q^{-1} \hat{\sigma}_{\overline{l}
  \, \overline{j}} \hat{\sigma}_{l+1\, \overline{l}} \notag \\
&= (E^{l+1}_{\overline{l}} - q^{\frac{1}{2}} E^l_{l+1}) (E^{\overline{l}}_
  {\overline{j}} - q^{j-l+1} E^j_l) - q^{-1} (E^{\overline{l}}_{\overline{j}} -
  q^{j-l+1} E^j_l) (E^{l+1}_{\overline{l}} - q^{\frac{1}{2}} E^l_{l+1})\notag\\
&= E^{l+1}_{\overline{j}} - q^{j-l+\frac{1}{2}} E^j_{l+1}.
\end{align*}

\noindent Hence

\begin{equation}
\boxed{\hat{\sigma}_{\overline{i}\, \overline{j}} = E^{\overline{i}}_
  {\overline{j}} - q^{(\rho, \varepsilon_{\overline{j}}- \varepsilon_
  {\overline{i}}) +1} E^j_{i}, \quad 1 \leq j <i \leq \lceil \tfrac{m}{2} 
  \rceil.}
\end{equation}


\section{Calculating $\hat{\sigma}_{\nu\mu},\hat{\sigma}_{\overline{\mu}\, \overline{\nu}}$} \label{munu}

\noindent The next step is to construct the operators of the form $\hat{\sigma}_{\nu\mu}$ and $\hat{\sigma}_{\overline{\mu}\, \overline{\nu}}$, an equally straightforward process.

\

\noindent We have already found that $\hat{\sigma}_{\mu\, \mu+1} = E^\mu_{\mu+1} + E^{\overline{\mu+1}}_{\overline{\mu}}$, for $\mu < k$.  Also, we know that

\begin{equation*}
\hat{\sigma}_{\nu\, \mu+1} = \hat{\sigma}_{\nu \mu} \hat{\sigma}_{\mu\, \mu+1}
  - q \hat{\sigma}_{\mu\, \mu+1}  \hat{\sigma}_{\nu \mu}, \quad \nu < \mu <k.
\end{equation*}

\noindent Combining these, we find that for $1<\mu < k:$

\begin{align*}
\hat{\sigma}_{\mu-1 \, \mu+1} 
&= \hat{\sigma}_{\mu-1\, \mu} \hat{\sigma}_{\mu\,
  \mu+1} - q \hat{\sigma}_{\mu\, \mu+1} \hat{\sigma}_{\mu-1\, \mu} \\
&= (E^{\mu-1}_{\mu} + E^{\overline{\mu}}_{\overline{\mu -1}}) (E^\mu_{\mu+1} 
  + E^{\overline{\mu+1}}_{\overline{\mu}}) - q (E^\mu_{\mu+1} + 
  E^{\overline{\mu+1}}_{\overline{\mu}}) (E^{\mu-1}_{\mu} + 
  E^{\overline{\mu}}_{\overline{\mu -1}}) \\
&= E^{\mu-1}_{\mu+1} - q E^{\overline{\mu+1}}_{\overline{\mu-1}}.
\end{align*}

\noindent  We hypothesise that $\hat{\sigma}_{\nu \mu} = E^{\nu}_{\mu} - (-1)^{\mu + \nu} q^{\mu - \nu - 1} E^{\overline{\mu}}_{\overline{\nu}}, \; \nu < \mu \leq k$.  Clearly this holds when $\mu = \nu +1$.  Hence we suppose it is true for $\hat{\sigma}_{\nu \mu}$ where $1 \leq \nu < \mu < k$ and show it holds for $\hat{\sigma}_{\nu\,  \mu+1}$:

\begin{align*}
 \hat{\sigma}_{\nu\, \mu+1} &= \hat{\sigma}_{\nu \,\mu} \hat{\sigma}_{\mu\,
  \mu+1} - q \hat{\sigma}_{\mu\, \mu+1} \hat{\sigma}_{\nu\, \mu} \notag \\
&= (E^{\nu}_{\mu} - (-1)^{\mu + \nu} q^{\mu-\nu-1} E^{\overline{\mu}}
  _{\overline{\nu}}) (E^{\mu}_{\mu+1} + E^{\overline{\mu+1}}_{\overline{\mu}})
  \\ 
& \qquad - q  (E^{\mu}_{\mu+1} + E^{\overline{\mu+1}}_{\overline{\mu}})
  (E^{\nu}_{\mu} - (-1)^{\mu + \nu} q^{\mu-\nu-1} E^{\overline{\mu}}
  _{\overline{\nu}}) \\
&= E^{\nu}_{\mu+1} - (-1)^{(\mu+1)+ \nu} q^{(\mu+1)-\nu-1} E^{\overline{\mu+1}}
  _{\overline{\nu}}.
\end{align*}

\noindent Thus we have proven, by induction, that

\begin{equation}
\boxed{\hat{\sigma}_{\nu \mu} = E^{\nu}_{\mu} - (-1)^{\mu + \nu} q^{(\rho, \delta_\mu - \delta_\nu) -1} E^{\overline{\mu}}_{\overline{\nu}}, \quad 1 \leq \nu < \mu \leq k.}
\end{equation}

\

\noindent Similarly, we have $\hat{\sigma}_{\overline{\mu+1}\, \overline{\mu}} = E^{\overline{\mu+1}}_{\overline{\mu}} + E^\mu_{\mu+1}$, and

\begin{equation*}
\hat{\sigma}_{\overline{\mu+1}\, \overline{\nu}} = \hat{\sigma}_{\overline
  {\mu+1}\, \overline{\mu}} \hat{\sigma}_{\overline{\mu}\, \overline{\nu}} - q 
  \hat{\sigma}_{\overline{\mu}\, \overline{\nu}} \hat{\sigma}_{\overline{\mu+1}
  \, \overline{\mu}}, \qquad \nu < \mu < k.
\end{equation*}

\noindent Using these, we deduce

\begin{alignat*}{2}
\hat{\sigma}_{\overline{\mu+1}\, \overline{\mu-1}} &= \hat{\sigma}_{\overline
  {\mu+1}\, \overline{\mu}} \hat{\sigma}_{\overline{\mu}\, \overline{\mu-1}}
  - q \hat{\sigma}_{\overline{\mu}\, \overline{\mu-1}} \hat{\sigma}_{\overline
  {\mu+1}\, \overline{\mu}} &&\\
&= (E^{\overline{\mu+1}}_{\overline{\mu}} + E^\mu_{\mu+1}) (E^{\overline{\mu}}_
  {\overline{\mu-1}} + E^{\mu-1}_{\mu}) - q (E^{\overline{\mu}}_{\overline
  {\mu-1}} + E^{\mu-1}_{\mu}) (E^{\overline{\mu+1}}_{\overline{\mu}} + 
  &&E^\mu_{\mu+1})\\
&= E^{\overline{\mu+1}}_{\overline{\mu-1}} - q E^{\mu-1}_{\mu +1},&& 1<\mu<k. 
\end{alignat*}

\noindent From this we suspect that $\hat{\sigma}_{\overline{\mu}\, \overline{\nu}} = E^{\overline{\mu}}_{\overline{\nu}} - (-1)^{\mu + \nu} q^{\mu - \nu -1} E^\nu_\mu$.  Assuming this is true for a given $\nu<\mu<k$, from the inductive relations \eqref{indrel} we find

\begin{align*}
\hat{\sigma}_{\overline{\mu+1}\, \overline{\nu}} &= \hat{\sigma}_{\overline
  {\mu+1}\, \overline{\mu}} \hat{\sigma}_{\overline{\mu}\, \overline{\nu}} - q 
  \hat{\sigma}_{\overline{\mu}\, \overline{\nu}} \hat{\sigma}_{\overline{\mu+1}
  \, \overline{\mu}} \\
&=  (E^{\overline{\mu+1}}_{\overline{\mu}} + E^\mu_{\mu+1}) (E^{\overline{\mu}}
  _{\overline{\nu}} - (-1)^{\mu+\nu} q^{\mu-\nu-1} E^{\nu}_{\mu})&&\\
& \qquad - q (E^{\overline{\mu}}_{\overline{\nu}} -(-1)^{\mu+\nu} q^{\mu-\nu-1}
  E^{\nu}_{\mu}) (E^{\overline{\mu+1}}_{\overline{\mu}} + E^\mu_{\mu+1})\\
&= E^{\overline{\mu+1}}_{\overline{\nu}} - (-1)^{(\mu+1) + \nu} q^{(\mu+1) - 
  \nu -1} E^{\nu}_{\mu+1},
\end{align*}

\noindent as expected. Therefore it follows by induction that

\begin{equation}
\boxed{\hat{\sigma}_{\overline{\mu}\, \overline{\nu}} = E^{\overline{\mu}}_
  {\overline{\nu}} - (-1)^{\mu + \nu} q^{(\rho, \delta_{\overline{\nu}} 
  - \delta_{\overline{\mu}}) -1} E^\nu_\mu, \quad 1 \leq \nu < \mu \leq k.}
\end{equation}


\section{Calculating $\hat{\sigma}_{\mu i}, \hat{\sigma}_{\overline{i}\, \overline{\mu}}$} \label{mui}

\noindent  The next step is to find the odd unknowns of the form $\hat{\sigma}_{\mu i}$ and $ \hat{\sigma}_{\overline{i}\, \overline{\mu}}$.  The process is simplified by using the values of $\hat{\sigma}_{ba}$ already calculated.

\

\noindent  In the case of $\hat{\sigma}_{\mu i}$ we first use $\hat{\sigma}_{\nu = k\, i=1}$ and $\hat{\sigma}_{\mu\nu}$ to find the general $\hat{\sigma}_{\mu \, i=1}$ and then apply $\hat{\sigma}_{ji}$ to extend this to all values of $\hat{\sigma}_{\mu i}$.

\

\noindent  Earlier in the chapter we showed that:

\begin{alignat*}{2}
\hat{\sigma}_{\nu = k\, i=1} &= E^{\nu=k}_{i=1} + (-1)^k q 
  E^{\overline{i} = \overline{1}}_{\overline{\nu} = \overline{k}}, &&\\
\hat{\sigma}_{\mu \nu} &= E^{\mu}_{\nu} - (-1)^{\nu+\mu} q^{\nu-\mu-1} 
  E^{\overline{\nu}}_{\overline{\mu}}, \qquad &&\mu <\nu \leq k, \\
\hat{\sigma}_{ji} &= E^j_i - q^{j-i+1} E^{\overline{i}}_{\overline{j}}, \qquad 
 && 1 \leq j < i \leq l, \\
\hat{\sigma}_{j\, l+1} &= E^j_{l+1} - q^{j-l-\frac{1}{2}} E^{\overline{l+1}}
  _{\overline{j}}, && 1 \leq j \leq l,\; m=2l+1.
\end{alignat*}

\noindent  The relations we apply are

\begin{alignat*}{2}
\hat{\sigma}_{\mu\, i=1} &= \hat{\sigma}_{\mu\,\nu=k} \hat{\sigma}_{\nu=k\,i=1}
  - q \hat{\sigma}_{\nu=k\, i=1} \hat{\sigma}_{\mu\, \nu=k}, \qquad&&\mu < k,\\
\hat{\sigma}_{\mu i} &= \hat{\sigma}_{\mu\,j=1} \hat{\sigma}_{j=1\,i} - q^{-1} 
  \hat{\sigma}_{j=1\, i} \hat{\sigma}_{\mu\, j=1}, && i>1.
\end{alignat*}

\noindent Note the second of these is taken from the general form of the inductive relations.  Although one could calculate $\hat{\sigma}_{\mu i}$ using only the relations in Tables \ref{list}, \ref{even} and \ref{odd}, at this point using the general form saves time.

\

\noindent Combining the information, we see

\begin{alignat*}{2}
\hat{\sigma}_{\mu\, i=1}& = \hat{\sigma}_{\mu\,\nu=k} \hat{\sigma}_{\nu=k\,i=1}
  - q \hat{\sigma}_{\nu=k\, i=1} \hat{\sigma}_{\mu\, \nu=k}, && \mu<k \\
&=(E^{\mu}_{\nu=k} - (-1)^{k+\mu} q^{k-\mu-1} E^{\overline{\nu}=\overline{k}}
   _{\overline{\mu}}) (E^{\nu=k}_{i=1} + (-1)^k q E^{\overline{i}=\overline{1}}
   _{\overline{\nu} = \overline{k}}) \\
& \qquad - q  (E^{\nu=k}_{i=1} + (-1)^k q E^{\overline{i} = \overline{1}}_
   {\overline{\nu}=\overline{k}}) (E^{\mu}_{\nu=k} - (-1)^{k+\mu} q^{k-\mu-1}
   E^{\overline{\nu}=\overline{k}}_{\overline {\mu}})  \\
&= E^\mu_{i=1} + (-1)^\mu q^{k-\mu +1} E^{\overline{i} = \overline{1}}_
   {\overline{\mu}}, && \mu \leq k.
\end{alignat*}
\newpage
\noindent Then for $i \leq l$

\begin{alignat*}{2}
\hat{\sigma}_{\mu i}& = \hat{\sigma}_{\mu\,j=1} \hat{\sigma}_{j=1\,i} - q^{-1} 
  \hat{\sigma}_{j=1\, i} \hat{\sigma}_{\mu\, j=1}, && 1<i \leq l \\
&= (E^\mu_{j=1} + (-1)^\mu q^{k-\mu +1} E^{\overline{j} = \overline{1}}_
   {\overline{\mu}})  (E^{j=1}_i - q^{2-i} E^{\overline{i}}_{\overline{j}=
   \overline{1}}) \\
& \qquad - q^{-1} (E^{j=1}_i - q^{2-i} E^{\overline{i}}_{\overline{j}=
   \overline{1}}) (E^\mu_{j=1} + (-1)^\mu q^{k-\mu +1} E^{\overline{j} = 
   \overline{1}}_{\overline{\mu}}) \notag \\
&= E^\mu_i + (-1)^\mu q^{k-\mu+2-i} E^{\overline{i}}_{\overline{\mu}}, &&
   1 \leq i \leq l\\
&= E^\mu_i + (-1)^{\mu} q^{(\rho, \varepsilon_i-\delta_\mu)+1} E^{\overline{i}}
   _{\overline{\mu}},&& 1 \leq i \leq l.
\end{alignat*}

\noindent It only remains to calculate $\hat{\sigma}_{\mu\, l+1}$ when $m=2l+1$.  We find

\begin{align*}
\hat{\sigma}_{\mu\, l+1} &= \hat{\sigma}_{\mu\,i=1} \hat{\sigma}_{i=1\, l+1} - 
  q^{-1} \hat{\sigma}_{i=1\,l+1} \hat{\sigma}_{\mu\,i=1} \notag \\
&=  (E^\mu_{i=1} + (-1)^{\mu} q^{k-\mu+1} E^{\overline{i}=\overline{1}}_
  {\overline{\mu}}) (E^{i=1}_{l+1} - q^{\frac{1}{2}-l}  E^{l+1}_{\overline{i}=
  \overline{1}}) \notag \\
& \qquad - q^{-1} (E^{i=1}_{l+1} - q^{\frac{1}{2}-l}  E^{l+1}_{\overline{i}=
  \overline{1}})  (E^\mu_{i=1} + (-1)^{\mu} q^{k-\mu+1} E^{\overline{i}=
  \overline{1}}_{\overline{\mu}}) \notag \\
&= E^\mu_{l+1} + (-1)^{\mu} q^{k-\mu-l+\frac{1}{2}} E^{l+1}_{\overline{\mu}}
  \notag  \\
&= E^\mu_{l+1} + (-1)^{\mu} q^{(\rho, \varepsilon_{l+1}-\delta_\mu)} 
  E^{l+1}_{\overline{\mu}}.
\end{align*}

\noindent Thus 

\begin{equation}
\boxed{ \hat{\sigma}_{\mu i} = E^\mu_{i} + (-1)^{\mu} q^{(\rho, \varepsilon_{i}
  - \delta_\mu) +(\varepsilon_i, \varepsilon_i)} E^{i}_{\overline{\mu}}, \quad
  1 \leq i \leq \lceil \tfrac{m}{2} \rceil,\; 1 \leq \mu \leq k.}
\end{equation}

\

\noindent To find $\hat{\sigma}_{\overline{i} \overline{\mu}}$ we follow the same procedure.  The known values for $\hat{\sigma}_{ba}$ and recurrence relations we use are:

\begin{alignat*}{2}
&\hat{\sigma}_{\overline{i}=\overline{1}\, \overline{\nu}=\overline{k}} = 
  E^{\overline{i}=\overline{1}}_{\overline{\nu}=\overline{k}} + (-1)^k q^{-1} 
  E^{\nu=k}_{i=1}, && \\
&\hat{\sigma}_{\overline{\nu}\, \overline{\mu}} = E^{\overline{\nu}}_
  {\overline{\mu}} - (-1)^{\nu+\mu} q^{\nu-\mu-1} E^{\mu}_{\nu}, && 1 \leq
  \mu < \nu \leq k, \\
&\hat{\sigma}_{\overline{i}\, \overline{j}} = E^{\overline{i}}_{\overline{j}}
  - q^{j-i+1} E^j_i, && 1 \leq j < i \leq l, \\
&\hat{\sigma}_{l+1\,\overline{j}} = E^{l+1}_{\overline{j}} -q^{j-1+\frac{1}{2}}
  E^j_{l+1}, && 1 \leq j \leq l,\; m=2l+1, \\
&\hat{\sigma}_{\overline{i} = \overline{1}\, \overline{\mu}} = \hat{\sigma}_
  {\overline{i} = \overline{1}\, \overline{\nu} = \overline{k}} \hat{\sigma}_
  {\overline{\nu} = \overline{k}\, \overline{\mu}} - q \hat{\sigma}_{\overline
  {\nu} = \overline{k}\, \overline{\mu}} \hat{\sigma}_{\overline{i} = \overline
  {1}\, \overline{\nu} = \overline{k}}, &\qquad &\mu < k, \\
&\hat{\sigma}_{\overline{i}\, \overline{\mu}} = \hat{\sigma}_{\overline{i}\, 
  \overline{j} = \overline{1}} \hat{\sigma}_{\overline{j} = \overline{1}\, 
  \overline{\mu}} - q \hat{\sigma}_{\overline{j} = \overline{1}\, \overline
  {\mu}} \hat{\sigma}_{\overline{i}\, \overline{j} = \overline{1}}, &&
  i > 1.
\end{alignat*}

\noindent From these we determine

\begin{alignat*}{2}
\hat{\sigma}_{\overline{i} = \overline{1}\, \overline{\mu}} &= \hat{\sigma}_
  {\overline{i} = \overline{1}\, \overline{\nu} = \overline{k}} \hat{\sigma}_
  {\overline{\nu} = \overline{k}\, \overline{\mu}} - q \hat{\sigma}_{\overline
  {\nu} = \overline{k}\, \overline{\mu}} \hat{\sigma}_{\overline{i} = \overline
  {1}\, \overline{\nu} = \overline{k}}, &\qquad &\mu < k \notag \\
&= (E^{\overline{i}=\overline{1}}_{\overline{\nu}=\overline{k}} + (-1)^kq^{-1} 
  E^{\nu=k}_{i=1}) (E^{\overline{\nu} = \overline{k}}_{\overline{\mu}} - 
  (-1)^{k+\mu} q^{k-\mu-1} E^{\mu}_{\nu=k}) && \notag \\
& \qquad - q (E^{\overline{\nu} = \overline{k}}_{\overline{\mu}}-(-1)^{k+\mu}
  q^{k-\mu-1} E^{\mu}_{\nu=k}) (E^{\overline{i}=\overline{1}}_{\overline{\nu}
  =\overline{k}} +(-1)^k q^{-1} E^{\nu=k}_{i=1}) && \notag \\
&= E^{\overline{i} = \overline{1}}_{\overline{\mu}} + (-1)^\mu q^{k-\mu-1}
   E^{\mu}_{i=1}, && \mu \leq k.
\end{alignat*}

\noindent Therefore, for $1 \leq \mu \leq k$, we find

\begin{alignat*}{2}
\hat{\sigma}_{\overline{i}\, \overline{\mu}} &= \hat{\sigma}_{\overline{i}\, 
  \overline{j} = \overline{1}} \hat{\sigma}_{\overline{j} = \overline{1}\, 
  \overline{\mu}} -q^{-1} \hat{\sigma}_{\overline{j} = \overline{1}\, \overline
  {\mu}} \hat{\sigma}_{\overline{i}\, \overline{j} = \overline{1}}, &&
  1 < i \leq l \notag \\
&= (E^{\overline{i}}_{\overline{j}=\overline{1}} - q^{2-i} E^{j=1}_i) 
  (E^{\overline{j} = \overline{1}}_{\overline{\mu}} + (-1)^\mu q^{k-\mu-1}
   E^{\mu}_{j=1}) \notag \\
&\qquad - q^{-1} (E^{\overline{j} = \overline{1}}_{\overline{\mu}} + (-1)^\mu 
  q^{k-\mu-1} E^{\mu}_{j=1}) (E^{\overline{i}}_{\overline{j}=\overline{1}} - 
  q^{2-i} E^{j=1}_i) \notag \\
&= E^{\overline{i}}_{\overline{\mu}} + (-1)^\mu q^{k-\mu-i} E^\mu_i, &\qquad&
  1 \leq i \leq l
\end{alignat*}

\noindent and

\begin{align*}
\hat{\sigma}_{l+1\, \overline{\mu}} &= \hat{\sigma}_{l+1\, \overline{j} = 
  \overline{1}} \hat{\sigma}_{\overline{j} = \overline{1}\, \overline{\mu}} -
  q^{-1} \hat{\sigma}_{\overline{j} = \overline{1}\, \overline{\mu}} 
  \hat{\sigma}_{l+1\, \overline{j} = \overline{1}} \notag \\
&= (E^{l+1}_{\overline{j}=\overline{1}} - q^{\frac{3}{2}-l} E^{j=1}_{l+1}) 
  (E^{\overline{j} = \overline{1}}_{\overline{\mu}} + (-1)^\mu q^{k-\mu-1}
   E^{\mu}_{j=1}) \notag \\
&\qquad - q^{-1} (E^{\overline{j} = \overline{1}}_{\overline{\mu}} + (-1)^\mu 
  q^{k-\mu-1} E^{\mu}_{j=1}) (E^{l+1}_{\overline{j}=\overline{1}} - 
  q^{\frac{3}{2}-l} E^{j=1}_{l+1})  \notag \\
&= E^{l+1}_{\overline{\mu}} + (-1)^{\mu} q^{k-\mu-l-\frac{1}{2}} E^{\mu}_{l+1}.
\end{align*}

\noindent Hence 

\begin{equation}
\boxed{\hat{\sigma}_{\overline{i}\,\overline{\mu}} = E^{\overline{i}}_
  {\overline{\mu}} + (-1)^\mu q^{(\rho, \delta_{\overline{\mu}} - \varepsilon_
  {\overline{i}})-1} E^\mu_i, \quad 1 \leq i \leq \lceil \tfrac{m}{2} \rceil,
  \, 1 \leq \mu \leq k.}
\end{equation}


\section{Calculating $\hat{\sigma}_{i\, \overline{j}}$} \label{ijbar}

\noindent Now we construct the operators of the form $\hat{\sigma}_{i\, \overline{j}}$, starting from the fundamental values associated with the root $\alpha_l$.  As we currently have different operators depending on the parity of $m$, the first step is to unify the two cases.

\

\noindent It has already been shown that when $m$ is even, $\hat{\sigma}_{l\, \overline{l}} = 0$, $\hat{\sigma}_{l-1\, \overline{l}} = E^{l-1}_{\overline{l}} - E^l_{\overline{l-1}}$ and $\hat{\sigma}_{l\, \overline{l-1}}= E^l_{\overline{l-1}} - E^{l-1}_{\overline{l}}$.  We now calculate those operators for the case $m=2l+1$, using the expressions for $\hat{\sigma}_{l\, l+1}$ and $\hat{\sigma}_{l+1\, \overline{l}}$ in Table \ref{basevector} together with those for $\hat{\sigma}_{l-1\, l+1}$ and $\hat{\sigma}_{l+1\, \overline{l-1}}$ as determined in Section \ref{ij}.  From the inductive relations \eqref{indrel} we find

\begin{align*}
\hat{\sigma}_{l\, \overline{l}} &= q \hat{\sigma}_{l\, l+1} \hat{\sigma}_{l+1\,
  \overline{l}} - \hat{\sigma}_{l+1\, \overline{l}} \hat{\sigma}_{l\, l+1} \\
&= q (E^l_{l+1} - q^{-\frac{1}{2}} E^{l+1}_{\overline{l}}) (E^{l+1}_
  {\overline{l}} - q^{\frac{1}{2}} E^l_{l+1}) - (E^{l+1}_{\overline{l}} - 
  q^{\frac{1}{2}} E^l_{l+1}) (E^l_{l+1} - q^{-\frac{1}{2}} 
  E^{l+1}_{\overline{l}}) \\
&= q E^l_{\overline{l}} - E^l_{\overline{l}}.
\end{align*}

\noindent Also,

\begin{align*}
\hat{\sigma}_{l-1\, \overline{l}} &= \hat{\sigma}_{l-1\, l+1} 
  \hat{\sigma}_{l+1\, \overline{l}} - \hat{\sigma}_{l+1\, \overline{l}} 
  \hat{\sigma}_{l-1\, l+1} \\
&= (E^{l-1}_{l+1} - q^{-\frac{3}{2}} E^{l+1}_{\overline{l-1}}) (E^{l+1}_
  {\overline{l}} - q^{\frac{1}{2}} E^l_{l+1}) - (E^{l+1}_{\overline{l}} - 
  q^{\frac{1}{2}} E^l_{l+1}) (E^{l-1}_{l+1} - q^{-\frac{3}{2}} 
  E^{l+1}_{\overline{l-1}}) \notag \\
&= E^{l-1}_{\overline{l}} - q^{-1} E^l_{\overline{l-1}} 
\end{align*}

\noindent and

\begin{align*}
\hat{\sigma}_{l\, \overline{l-1}} &= \hat{\sigma}_{l\,l+1} \hat{\sigma}_{l+1\, 
  \overline{l-1}} - \hat{\sigma}_{l+1\, \overline{l-1}} \hat{\sigma}_{l\,l+1}\\
&= (E^l_{l+1} - q^{-\frac{1}{2}} E^{l+1}_{\overline{l}}) (E^{l+1}_{\overline
  {l-1}} - q^{-\frac{1}{2}} E^{l-1}_{l+1})  - (E^{l+1}_{\overline{l-1}} - 
  q^{-\frac{1}{2}} E^{l-1}_{l+1}) (E^l_{l+1} - q^{-\frac{1}{2}} E^{l+1}_
  {\overline{l}}) \\
&= E^l_{\overline{l-1}} - q^{-1} E^{l-1}_{\overline{l}}.
\end{align*}

\noindent  Hence we can say that for both odd and even $m$,

\begin{align*}
&\hat{\sigma}_{l\, \overline{l}} = q E^l_{\overline{l}} - q^{2l-m+1} 
  E^l_{\overline{l}}, \\
&\hat{\sigma}_{l-1\, \overline{l}} = E^{l-1}_{\overline{l}} - q^{2l-m} 
  E^l_{\overline{l-1}}, \\
&\hat{\sigma}_{l\, \overline{l-1}} = E^l_{\overline{l-1}} - q^{2l-m} E^{l-1}_
  {\overline{l}}.
\end{align*}

\noindent Now that we have a formula for $\hat{\sigma}_{l\, \overline{l}}$ and $\hat{\sigma}_{l-1\, \overline{l}}$ that hold for any $m$, we can use the general form of the inductive relations to calculate the remaining values of $\hat{\sigma}_{i\overline{j}}$.  First we find $\hat{\sigma}_{i\overline{l}}$ for all $i < l$, remembering that

\begin{equation*}
\hat{\sigma}_{ji} = E^j_i - q^{j-i+1} E^{\overline{i}}_{\overline{j}}, \qquad 1 \leq j < i \leq l.
\end{equation*}

\noindent Then

\begin{alignat*}{2}
\hat{\sigma}_{i\overline{l}} &= \hat{\sigma}_{i\, l-1} \hat{\sigma}_{l-1\, 
  \overline{l}} - q^{-1} \hat{\sigma}_{l-1\, \overline{l}} \hat{\sigma}_{i\, 
  l-1}, 
  &\quad& i<l-1 \\
&=  (E^i_{l-1} - q^{i-l+2} E^{\overline{l-1}}_{\overline{i}}) (E^{l-1}_
  {\overline{l}} - q^{2l-m} E^l_{\overline{l-1}}) \\
& \qquad - q^{-1} (E^{l-1}_{\overline{l}} - q^{2l-m} E^l_{\overline{l-1}})
  (E^i_{l-1} - q^{i-l+2} E^{\overline{l-1}}_{\overline{i}}) \notag \\
&= E^i_{\overline{l}} - q^{i+l-m+1} E^l_{\overline{i}}, && i \leq l-1.
\end{alignat*}

\noindent Similarly,

\begin{alignat*}{2}
\hat{\sigma}_{l\overline{j}} &= \hat{\sigma}_{l\,\overline{l-1}} \hat{\sigma}_
  {\overline{l-1}\,\overline{j}} - q^{-1}  \hat{\sigma}_{\overline{l-1}\,
  \overline{j}} \hat{\sigma}_{l\,\overline{l-1}}, & \qquad & j < l-1 \notag \\
&= (E^l_{\overline{l-1}} - q^{2l-m} E^{l-1}_{\overline{l}}) (E^{\overline{l-1}}
  _{\overline{j}} - q^{j-l+2} E^j_{l-1}) \notag \\
& \qquad - q^{-1} (E^{\overline{l-1}}_{\overline{j}} - q^{j-l+2} E^j_{l-1})
  (E^l_{\overline{l-1}} - q^{2l-m} E^{l-1}_{\overline{l}}) \notag \\
&= E^l_{\overline{j}} - q^{j+l-m+1} E^j_{\overline{l}},&& j \leq l-1.
\end{alignat*}

\noindent Now we evaluate the remaining operators of the form $\hat{\sigma}_{i\overline{j}}$, using the formula for $\hat{\sigma}_{\overline{i}\,\overline{j}}$ derived in Section \ref{ij} together with $\hat{\sigma}_{i\overline{l}}$.  We find

\begin{alignat*}{2}
\hat{\sigma}_{i\overline{j}} &=q^{-(\varepsilon_i, \varepsilon_{\overline{j}})}
  \hat{\sigma}_{i\overline{l}} \hat{\sigma}_{\overline{l}\, \overline{j}} - 
  q^{-1} \hat{\sigma}_{\overline{l}\, \overline{j}} \hat{\sigma}_
  {i\overline{l}}, &\qquad& i,j < l \\
&= q^{-(\varepsilon_i, \varepsilon_{\overline{j}})} (E^i_{\overline{l}} - 
  q^{i+l-m+1} E^l_{\overline{i}})(E^{\overline{l}}_{\overline{j}} - q^{j-l+1} 
  E^j_l) \\\
& \qquad - q^{-1} (E^{\overline{l}}_{\overline{j}} - q^{j-l+1} E^j_l) 
  (E^i_{\overline{l}} - q^{i+l-m+1} E^l_{\overline{i}})  \\
&= q^{-(\varepsilon_i, \varepsilon_{\overline{j}})} E^i_{\overline{j}} - 
  q^{i+j-m+1} E^j_{\overline{i}}, && i,j \leq l.
\end{alignat*}

\noindent Hence we have shown

\begin{equation}
\boxed{\hat{\sigma}_{i\, \overline{j}} = q^{-(\varepsilon_i, \varepsilon_{\overline{j}})} E^i_{\overline{j}} - q^{(\rho, \varepsilon_{\overline{j}} - \varepsilon_i)+1} E^j_{\overline{i}}, \qquad 1 \leq i,\;j \leq l.}
\end{equation}


\section{Calculating $\hat{\sigma}_{i\,\overline{\mu}}, \hat{\sigma}_{\mu\, \overline{i}}$} \label{imubar}

\noindent Next we construct the remaining odd operators, $\hat{\sigma}_{i\,\overline{\mu}}$ and $\hat{\sigma}_{\mu\, \overline{i}}$.  These are easily calculated from the operators derived earlier.  We simply combine our previous results for $\hat{\sigma}_{i\, \overline{j}}$ and $\hat{\sigma}_{\overline{j}\, \overline{\mu}}$, using the unified form of the relations (\ref{indrel}).  They tell us that for any $j \neq i,\; 1 \leq i,\;j \leq l,\; 1 \leq \mu \leq k$,  

\begin{equation*}
\hat{\sigma}_{i\, \overline{\mu}} = \hat{\sigma}_{i\, \overline{j}} 
  \hat{\sigma}_{\overline{j}\, \overline{\mu}}- q^{-1} \hat{\sigma}_
  {\overline{j}\, \overline{\mu}} \hat{\sigma}_{i\, \overline{j}}.
\end{equation*}

\noindent Substituting in the values found in Sections \ref{mui} and \ref{ijbar}, we find

\begin{align*}
 \hat{\sigma}_{i\, \overline{\mu}} &= (E^i_{\overline{j}} - q^{i+j-m+1} 
  E^j_{\overline{i}}) (E^{\overline{j}}_{\overline{\mu}} + (-1)^\mu q^{k-\mu-j}
  E^\mu_j) \notag \\
& \qquad - q^{-1} (E^{\overline{j}}_{\overline{\mu}} + (-1)^\mu q^{k-\mu-j} 
  E^\mu_j) (E^i_{\overline{j}} - q^{i+j-m+1} E^j_{\overline{i}}) \notag \\
&= E^i_{\overline{\mu}} + (-1)^\mu q^{i-m+k-\mu} E^\mu_{\overline{i}}.
\end{align*}

\noindent Similarly, for any $j \neq i, \; 1 \leq i \leq l, \; 1 \leq \mu \leq k$ we obtain

\begin{alignat*}{2}
\hat{\sigma}_{\mu\, \overline{i}} &= \hat{\sigma}_{\mu j} \hat{\sigma}_
  {j\, \overline{i}} - q^{-1} \hat{\sigma}_{j\, \overline{i}} \hat{\sigma}_
  {\mu j}, && \qquad i \neq j \\ 
&= (E^\mu_j + (-1)^{\mu} q^{k-\mu + 2-j} E^{\overline{j}}_{\overline{\mu}})
  (E^j_{\overline{i}} - q^{i+j-m+1} E^i_{\overline{j}}) \notag \\
& \qquad - q^{-1} (E^j_{\overline{i}} - q^{i+j-m+1} E^i_{\overline{j}}) 
  (E^\mu_j + (-1)^{\mu} q^{k-\mu + 2-j} E^{\overline{j}}_{\overline{\mu}}) \\
&= E^\mu_{\overline{i}} + (-1)^\mu q^{i-m+k-\mu+2} E^i_{\overline{\mu}}.
\end{alignat*}

\noindent Thus in terms of the graded-half sum of positive roots we have

\begin{equation}
\boxed{\hat{\sigma}_{i\, \overline{\mu}}  = E^i_{\overline{\mu}} + (-1)^\mu 
  q^{(\rho, \delta_{\overline{\mu}} - \varepsilon_i) -1} E^\mu_{\overline{i}}, 
  \quad 1 \leq i \leq l,\; 1 \leq \mu \leq k} 
\end{equation}

\noindent and

\begin{equation}
\boxed{\hat{\sigma}_{\mu\, \overline{i}} = E^\mu_{\overline{i}} + (-1)^\mu
  q^{(\rho, \varepsilon_{\overline{i}} - \delta_{\mu})+1} E^i_{\overline{\mu}},
  \quad 1 \leq i \leq l,\; 1 \leq \mu \leq k.}
\end{equation}


\section{Calculating $\hat{\sigma}_{\mu\, \overline{\nu}}$}

\noindent The only operators yet to be calculated are the $\hat{\sigma}_{\mu\,\overline{\nu}}$.  As with $\hat{\sigma}_{\mu\, \overline{i}}$, it is a straightforward piecing together of values of $\hat{\sigma}_{ba}$ that have already been determined.  In this case we use the results for $\hat{\sigma}_{\mu i}$ from Section \ref{mui} together with those for $\hat{\sigma}_{i\, \overline{\nu}}$ from the previous section.

\

\noindent Using relations \eqref{indrel}, we obtain for all $1 \leq \mu, \nu \leq k$

\begin{alignat*}{2}
\hat{\sigma}_{\mu\, \overline{\nu}} &= q^{(\delta_\mu,\delta_\nu)} \hat{\sigma}
  _{\mu i} \hat{\sigma}_{i\, \overline{\nu}} + q^{-1} \hat{\sigma}_{i\, 
  \overline{\nu}} \hat{\sigma}_{\mu i}, &&\qquad 1 \leq i \leq l \\
&= q^{(\delta_\mu,\delta_\nu)} (E^\mu_i + (-1)^{\mu} q^{k-\mu + 2-i} 
  E^{\overline{i}}_{\overline{\mu}}) (E^i_{\overline{\nu}} + (-1)^\nu 
  q^{i-m+k-\nu} E^\nu_{\overline{i}})&& \\
& \qquad + q^{-1} (E^i_{\overline{\nu}} + (-1)^\nu q^{i-m+k-\nu} 
  E^\nu_{\overline{i}}) (E^\mu_i + (-1)^{\mu} q^{k-\mu + 2-i} E^{\overline{i}}_
  {\overline{\mu}})&& \\ 
&= q^{(\delta_\mu,\delta_\nu)} E^\mu_{\overline{\nu}} + (-1)^{\mu + \nu} 
  q^{n - m - \mu - \nu +1} E^\nu_{\overline{\mu}}.&&
\end{alignat*}

\noindent Thus, completing the generators for the $R$-matrix, we have

\begin{equation} \label{sigmunubar}
\boxed{ \hat{\sigma}_{\mu\, \overline{\nu}} = q^{-(\delta_\mu,\delta_ 
  {\overline{\nu}})} E^\mu_{\overline{\nu}} + (-1)^{\mu+\nu} q^{(\rho, 
  \delta_{\overline{\nu}} - \delta_\mu) -1} E^\nu_{\overline{\mu}}, \quad
  1 \leq \mu, \nu \leq k.}
\end{equation}


\section{Solution for the $R$-matrix in the Vector Representation.} \label{RmatrixVector}

\noindent  In this chapter we have calculated the explicit solution for all the $\hat{\sigma}_{ba},\; \varepsilon_b > \varepsilon_a$, in the vector representation, which form the basis for the $R$-matrix.  Looking over equations \eqref{sigji} through to \eqref{sigmunubar} we can see a general form for $\hat{\sigma}_{ba}$, namely

\begin{equation*}
\hat{\sigma}_{ba} = q^{-(\varepsilon_a, \varepsilon_b)} E^b_a - (-1)^{[b]([a]+[b])} \xi_a \xi_b q^{(\varepsilon_a, \varepsilon_a)} q^{(\rho, \varepsilon_a - \varepsilon_b)} E^{\overline{a}}_{\overline{b}}, \qquad \varepsilon_b > \varepsilon_a.
\end{equation*}

\noindent Thus we have shown the $R$-matrix for the vector representation of $U_q[osp(m|n)]$, $R = (\pi \otimes \pi) \mathcal{R}$, is given by

\begin{equation*}
R = q ^{h_{j} \otimes h^{j}} \Bigl[ I \otimes I + (q - q^{-1}) \sum_{\varepsilon_{b} > \varepsilon_{a}} (-1)^{[b]} E^a_b \otimes \hat{\sigma}_{ba} \Bigr]
\end{equation*}

\noindent where

\begin{equation*}
\hat{\sigma}_{ba} = q^{-(\varepsilon_a, \varepsilon_b)} E^b_a - (-1)^{[b]([a]+[b])} \xi_a \xi_b q^{(\varepsilon_a, \varepsilon_a)} q^{(\rho, \varepsilon_a - \varepsilon_b)} E^{\overline{a}}_{\overline{b}}.
\end{equation*}

\noindent This can be written in a more elegant form.  Recall that the ansatz for $R$ can also be written as

\begin{equation*}
R= \sum_a E^a_a \otimes q^{h_{\varepsilon_a}} \; + \; (q-q^{-1}) \sum_{\varepsilon_{b} > \varepsilon_{a}} (-1)^{[b]} E^a_b \otimes q^{h_{\varepsilon_a}} \hat{\sigma}_{ba}.
\end{equation*}

\noindent  In the vector representation, this is equal to

\begin{equation*}
R = \sum_{a,b} q^{(\varepsilon_{a}, \varepsilon_{b})} E^a_a \otimes E^b_b \; + \; (q-q^{-1}) \sum_{\varepsilon_{b} > \varepsilon_{a}} (-1)^{[b]} E^a_b \otimes \tilde{\sigma}_{ba}
\end{equation*}

\noindent in terms of generators 

\begin{align*}
\tilde{\sigma}_{ba} &= q^{h_{\varepsilon_a}} \hat{\sigma}_{ba} \\
&= E^b_a - (-1)^{[b]([a]+[b])} \xi_a \xi_b q^{(\rho, \varepsilon_a - \varepsilon_b)} E^{\overline{a}}_{\overline{b}}.
\end{align*}

\

\noindent Hence we have the following result:

\begin{theorem}  The $R$-matrix for the vector representation, $R = (\pi \otimes \pi) \mathcal{R}$, is given by 

\begin{equation*}
R = \sum_{a,b} q^{(\varepsilon_{a}, \varepsilon_{b})} E^a_a \otimes E^b_b \; + \; (q-q^{-1}) \sum_{\varepsilon_{b} > \varepsilon_{a}} (-1)^{[b]} E^a_b \otimes \tilde{\sigma}_{ba},
\end{equation*} 

\noindent where

\begin{equation*} 
\tilde{\sigma}_{ba} = E^b_a - (-1)^{[b]([a]+[b])} \xi_a \xi_b q^{(\rho, 
\varepsilon_a - \varepsilon_b)} E^{\overline{a}}_{\overline{b}}, \qquad 
\varepsilon_b > \varepsilon_a.
\end{equation*}
\end{theorem}

\vspace{20mm}

\noindent  We can also explicitly find the opposite $R$-matrix $R^T$, using 

\begin{equation*}
(E^a_b)^\dagger = (-1)^{[a]([a]+[b])} E^b_a.
\end{equation*}

\noindent From equation \eqref{RT} on page \pageref{RT} we have:

\begin{equation*}
 R^T = \sum_a E^a_a \otimes q^{h_{\varepsilon_a}} + (q-q^{-1}) \sum_
  {\varepsilon_b > \varepsilon_a} (-1)^{[a]} E^b_a \otimes \hat{\sigma}_{ab} 
  q^{h_{\varepsilon_a}},
\end{equation*}

\noindent where

\begin{equation*}
\hat{\sigma}_{ab} = (-1)^{[b]([a]+[b])} \hat{\sigma}_{ba}^\dagger, \qquad
  \varepsilon_a < \varepsilon_b.
\end{equation*}

\noindent  Set $\tilde{\sigma}_{ab} = \hat{\sigma}_{ab} q^{h_{\varepsilon_a}}$ for $\varepsilon_a < \varepsilon_b$, so 

\begin{align*}
\tilde{\sigma}_{ab} &= (-1)^{[b]([a]+[b])} \hat{\sigma}_{ba}^\dagger 
  (q^{h_{\varepsilon_a}})^\dagger \\
&=  (-1)^{[b]([a]+[b])} \tilde{\sigma}_{ba}^\dagger \\
&= E^a_b - (-1)^{[a]([a]+[b])} \xi_a \xi_b q^{(\rho, \varepsilon_a - 
  \varepsilon_b)} E^{\overline{b}}_{\overline{a}}.
\end{align*}

\noindent Hence we have the following result for $R^T$:

\begin{theorem}
\noindent The opposite $R$-matrix for the vector representation, \\
$R^T = (\pi \otimes \pi) \mathcal{R}^T$, is given by

\begin{equation*}
R^T = \sum_{a,b} q^{(\varepsilon_a, \varepsilon_b)} E^a_a \otimes E^b_b + (q-q^{-1}) \sum_{\varepsilon_b > \varepsilon_a} (-1)^{[a]} E^b_a \otimes \tilde{\sigma}_{ab},
\end{equation*}

\noindent where 

\begin{equation*}
\tilde{\sigma}_{ab} = E^a_b - (-1)^{[a]([a]+[b])} \xi_a \xi_b q^{(\rho, 
  \varepsilon_a - \varepsilon_b)} E^{\overline{b}}_{\overline{a}}.
\end{equation*}

\end{theorem}

\

\noindent These formulae for $R$ and $R^T$ on the vector representation agree with those given in \cite{Mehta}.  In that thesis the $R$-matrix for the vector representation was calculated using projection operators onto invariant submodules of the tensor product.  The greatest advantage of the current method is it gives a straightforward way of constructing a solution to the Yang-Baxter Equation in an arbitrary representation of $U_q[osp(m|n)]$.

\chapter{Casimir Invariants and their Eigenvalues} \label{casimir}

\noindent The Lax operator can be used not only to construct solutions of the quantum Yang-Baxter Equation, but also to find families of Casimir invariants.  These are an important tool for understanding the representation theory of the superalgebra.  After constructing the Casimir invariants we can use properties of the root system to calculate their eigenvalues when acting on an irreducible highest weight module.

\

\noindent In this chapter we do exactly that, basing our method upon that used in \cite{Bincer} and \cite{Scheunert83} for the classical general and orthosymplectic superalgebras respectively.  This was adapted in \cite{LinksZhang} to cover the quantum superalgebra $U_q[gl(m|n)]$.   Although the concepts are much the same as in those cases, the combination of the $q$-deformation and the more complex root system of $U_q[osp(m|n)]$ makes the calculations in this chapter substantially more technically challenging.

\section{Casimir Invariants of $U_q[osp(m|n)]$}

\noindent Before constructing the Casimir invariants we need to define a new object. Let $h_\rho$ be the unique element of the Cartan subalgebra H satisfying

\begin{equation*}
\alpha_i(h_\rho) = (\rho, \alpha_i), \qquad \forall \alpha_i \in H^*.
\end{equation*}

\

\noindent It is also convenient to define a new operator $\partial$ by

\begin{equation*}
\partial \equiv (\pi \otimes \text{id}) \Delta.
\end{equation*}

\noindent Then from \cite{ZhangGould} we have the following theorem:

\begin{theorem}
Let $V$ be the representation space of $\pi$, an arbitrary finite dimensional representation of $U_q[osp(m|n)]$.  If $\Gamma \in (\text{End }V) \otimes U_q[osp(m|n)]$ satisfies

\begin{equation} \label{partial}
\partial (a) \Gamma = \Gamma \partial (a), \quad \forall a \in U_q[osp(m|n)],
\end{equation}

\noindent then

\begin{equation*}
C = (str \otimes \mathrm{id}) (\pi (q^{2h_\rho}) \otimes I) \Gamma
\end{equation*}

\noindent belongs to the centre of $U_q[osp(m|n)]$.  Above $str$ denotes the supertrace.
\end{theorem}

\noindent Now choose $\pi$ to be the vector representation.  Recalling that an $R$-matrix satisfies

\begin{equation*} 
\mathcal{R} \Delta(a) = \Delta^T(a) \mathcal{R}, \qquad \forall a \in U_q[osp(m|n)],
\end{equation*}

\noindent it is clear that

\begin{equation*}
\partial(a) R^T R = R^T R \,\partial(a), \qquad \forall a \in U_q[osp(m|n)].
\end{equation*}

\noindent Hence if we set $A \in (\text{End}\;V) \otimes U_q[osp(m|n)]$ to be

\begin{equation*}
A = \frac{(R^T R - I \otimes I)}{(q-q^{-1})},\
\end{equation*}

\noindent the operators $A^l$ will satisfy condition \eqref{partial} for all non-negative integers $l$. Thus the operators $C_l$ defined as

\begin{equation*}
C_l = (str \otimes \text{id})(\pi(q^{2h_p}) \otimes I) A^l, \qquad l \in \mathbb{Z}^+,
\end{equation*}

\noindent form a family of Casimir invariants.  Here $A$ coincides with the matrix of Jarvis and Green \cite{Jarvis} in the classical limit $q \rightarrow 1$, as do the invariants $C_l$.
 
\

\noindent Now write the Lax operator $R$ and its opposite $R^T$ in the form

\begin{align*}
R &= I \otimes I + (q-q^{-1}) \sum_{\varepsilon_b \geq \varepsilon_a} E^a_b 
  \otimes X^b_a, \\
R^T &= I \otimes I + (q-q^{-1}) \sum_{\varepsilon_b \leq \varepsilon_a} E^a_b 
  \otimes X^b_a.
\end{align*}

\noindent In terms of the operators $\hat{\sigma}_{ba}$, this implies 

\begin{equation*}
X^b_a = \begin{cases}
\frac{q^{h_{\varepsilon_a}}-I}{q-q^{-1}}, &a=b, \\
(-1)^{[b]} q^{h_{\varepsilon_a}} \hat{\sigma}_{ba}, \quad &\varepsilon_a < 
  \varepsilon_b, \\
(-1)^{[b]} \hat{\sigma}_{ba} q^{h_{\varepsilon_b}}, & \varepsilon_a > 
  \varepsilon_b.
\end{cases}
\end{equation*}

\noindent Writing $A$ as

\begin{equation*}
A = \sum_{a,b} E^a_b \otimes A^b_a,
\end{equation*}

\noindent we obtain

\begin{equation*}
A^b_a = (1+\delta^a_b) X^b_a + (q-q^{-1}) \sum_{\varepsilon_c \leq \varepsilon_a, \varepsilon_b} (-1)^{([a]+[c])([b]+[c])} X^c_a X^b_c.
\end{equation*}

\noindent This produces a family of Casimir invariants 

\begin{equation*}
C_l = \sum_a (-1)^{[a]} q^{(2\rho, \varepsilon_a)} {A^{(l)}}^a_a,
\end{equation*}

\noindent where the operators ${A^{(l)}}^b_a$ are recursively defined as

\begin{equation} \label{Alba}
{A^{(l)}}^b_a = \sum_c (-1)^{([a]+[c])([b]+[c])} {A^{(l-1)}}^c_a A^b_c.
\end{equation}


\section{Setting up the Eigenvalue Calculations}

\noindent Now that we have found a family of Casimir invariants, we wish to calculate their eigenvalues on a general irreducible finite-dimensional module. Let $V(\Lambda)$ be an arbitrary irreducible finite-dimensional module with highest weight $\Lambda$ and highest weight state $| \Lambda \rangle$.  Define $t_a^{(l)}$ to be the eigenvalue of ${A^{(l)}}^a_a$ on this state, so

\begin{equation*}
 {A^{(l)}}^a_a |\Lambda \rangle = t_a^{(l)} |\Lambda \rangle.
\end{equation*}

\noindent Once we have calculated $t_a^{(l)}$ we will use the result to find the eigenvalues of the Casimir invariants $C_l$. 

\

\noindent To evaluate $t_a^{(l)}$, note that if $\varepsilon_b > \varepsilon_a$ then ${A^{(l)}}^b_a$ is a raising operator, implying ${A^{(l)}}^b_a |\Lambda \rangle = 0$. Thus from equation \eqref{Alba} we deduce

\begin{align*}
t_a^{(l)} |\Lambda \rangle &= t_a^{(l-1)} t_a^{(1)} |\Lambda \rangle + 
  \sum_{\varepsilon_a < \varepsilon_b} (-1)^{[a]+[b]} {A^{(l-1)}}^b_a A^a_b 
  |\Lambda \rangle \notag \\
&= t_a^{(l-1)} t_a^{(1)} |\Lambda \rangle + \sum_{\varepsilon_a< \varepsilon_b}
  (-1)^{[a]+[b]} {A^{(l-1)}}^b_a \bigl[ X^a_b + (q-q^{-1}) X^a_b X^a_a \bigr] 
  |\Lambda 
  \rangle \notag \\
&= t_a^{(l-1)} t_a^{(1)} |\Lambda \rangle + \sum_{\varepsilon_a< \varepsilon_b}
  (-1)^{[a]+[b]} q^{(\Lambda, \varepsilon_a)} {A^{(l-1)}}^b_a X^a_b
  |\Lambda \rangle.
\end{align*}

\noindent Now we know that 

\begin{equation} \label{AdelX}
A^l \partial (X^a_b) =  \partial (X^a_b) A^l.
\end{equation}

\noindent This can be used to calculate ${A^{(l)}}^b_a X^a_b|\Lambda \rangle$ for $\varepsilon_a< \varepsilon_b$. First we need an expression for $\Delta(X^a_b)$. The $R$-matrix properties give

\begin{align*}
&(\Delta \otimes I) R = R_{13} R_{23} \\
\Rightarrow \quad &(I \otimes \Delta) R^T = R^T_{12} R^T_{13}.
\end{align*}

\noindent In terms of $X^a_b$, this implies

\begin{align*}
I \otimes I \otimes I + &(q-q^{-1}) \sum_{\varepsilon_a \leq \varepsilon_b} 
  E^b_a \otimes \Delta (X^a_b) \\
&= \bigl( I \otimes I \otimes I + (q-q^{-1}) \sum_{\varepsilon_a \leq 
  \varepsilon_b} E^b_a \otimes X^a_b \otimes I \bigr)  \\
& \hspace{2cm}\times \bigl( I \otimes I \otimes I + (q-q^{-1}) \sum_
  {\varepsilon_a \leq \varepsilon_b} E^b_a \otimes I  \otimes X^a_b \bigr) \\
&= I \otimes I \otimes I + (q-q^{-1}) \sum_{\varepsilon_a \leq \varepsilon_b}
  E^b_a \otimes (X^a_b \otimes I + I \otimes X^a_b) \\
& \hspace{2cm} + (q-q^{-1})^2 \sum_{\varepsilon_a \leq \varepsilon_c \leq 
  \varepsilon_b} (-1)^{([a]+[c])([b]+[c])} E^b_a \otimes X^c_b \otimes X^a_c.
\end{align*}

\noindent  Hence for all $\varepsilon_a < \varepsilon_b$ 

\begin{equation*}
 \Delta (X^a_b) = X^a_b \otimes I + I \otimes X^a_b + (q-q^{-1}) \sum_{\varepsilon_a \leq \varepsilon_c \leq \varepsilon_b} (-1)^{([a]+[c])([b]+[c])} X^c_b \otimes X^a_c.
\end{equation*}

\noindent We also need an expression for $\pi (X^a_b)$ for $\varepsilon_a \leq \varepsilon_b$.  At the end of the previous chapter we found that the generators for $R^T$ in the vector representation are given by

\begin{equation*}
\tilde{\sigma}_{ab} = \hat{\sigma}_{ab} q^{h_{\varepsilon_a}} = E^a_b - (-1)^{[a]([a]+[b])} \xi_a \xi_b q^{(\rho, \varepsilon_a - \varepsilon_b)} E^{\overline{b}}_{\overline{a}}, \qquad \varepsilon_a < \varepsilon_b.
\end{equation*}

\noindent From this we deduce that 

\begin{align*}
\pi (X^a_b) &= (-1)^{[a]} \pi(\hat{\sigma}_{ab} q^{h_{\varepsilon_a}}) \\
&= (-1)^{[a]} E^a_b - (-1)^{[a][b]} \xi_a \xi_b q^{(\rho, \varepsilon_a - 
  \varepsilon_b)} E^{\overline{b}}_{\overline{a}}, \qquad \varepsilon_a < 
  \varepsilon_b.
\end{align*}

\noindent Also, we know

\begin{align*}
\pi(X^a_a) &= (q-q^{-1})^{-1} \pi (q^{h_{\varepsilon_a}} - I)  \\
&= (q-q^{-1})^{-1} (q^{(\varepsilon_a, \varepsilon_a)(E^a_a - 
  E^{\overline{a}}_{\overline{a}})} -I).
\end{align*}

\noindent Applying these, we find that if $\varepsilon_a  < \varepsilon_b$ then

\begin{align*}
\partial (X^a_b) &= (\pi \otimes I) \Delta (X^a_b) \notag \\
&= \pi(X^a_b) \otimes \bigl( I + (q-q^{-1}) X^a_a \bigr) + \bigl( I+ (q-q^{-1})
  \pi (X^b_b) \bigr) \otimes X^a_b \notag \\
& \quad + (q-q^{-1}) \sum_{\varepsilon_a < \varepsilon_c < \varepsilon_b} 
  (-1)^{([a]+[c])([b]+[c])} \pi(X^c_b) \otimes X^a_c \notag \\
&= \bigl( (-1)^{[a]} E^a_b - (-1)^{[a][b]} \xi_a \xi_b q^{(\rho, \varepsilon_a 
  - \varepsilon_b)} E^{\overline{b}}_{\overline{a}} \bigr) \otimes 
  q^{h_{\varepsilon_a}} + q^{(\varepsilon_b, \varepsilon_b)(E^b_b - 
  E^{\overline{b}}_{\overline{b}})} \otimes X^a_b \notag \\
& \quad + (q-q^{-1}) \sum_{\varepsilon_a < \varepsilon_c < \varepsilon_b} 
  (-1)^{([a]+[c])([b]+[c])} \notag \\
& \hspace{45mm}\times \bigl( (-1)^{[c]} E^c_b - (-1)^{[b][c]} \xi_b \xi_c
  q^{(\rho, \varepsilon_c - \varepsilon_b)} E^{\overline{b}}_{\overline{c}} 
  \bigr) \otimes  X^a_c.
\end{align*}

\noindent Substituting this expression into equation \eqref{AdelX} and equating the $(a,b)$ entries, we find

\begin{align*}
(-1)^{[a]} &{A^{(l)}}^a_a q^{h_{\varepsilon_a}} - \delta^a_{\overline{b}}
  (-1)^{[a][b]} \xi_a \xi_b q^{(\rho, \varepsilon_a - \varepsilon_b)} {A^{(l)}}
  ^a_a q^{h_{\varepsilon_a}} + q^{(\varepsilon_b, \varepsilon_b)} {A^{(l)}}^b_a
  X^a_b \notag \\
& \qquad +(q-q^{-1}) \sum_{\varepsilon_a< \varepsilon_c <\varepsilon_b} \bigl( 
  (-1)^{[c]} {A^{(l)}}^c_a X^a_c - \delta^b_{\overline{c}} (-1)^{[b][c]} \xi_b 
  \xi_c q^{(\rho, \varepsilon_c - \varepsilon_b)} {A^{(l)}}^{\overline{b}}_a 
  X^a_c \bigr) \notag \\
=& (-1)^{[a]} q^{h_{\varepsilon_a}} {A^{(l)}}^b_b - \delta^a_{\overline{b}} 
  (-1)^{[a][b]} \xi_a \xi_b q^{(\rho, \varepsilon_a - \varepsilon_b)} 
  q^{h_{\varepsilon_a}} {A^{(l)}}^b_b + (-1)^{[a]+[b]} q^{(\varepsilon_a, 
  \varepsilon_b)} X^a_b {A^{(l)}}^b_a \notag \\
&\qquad -(q-q)^{-1} \delta^a_{\overline{b}} \sum_{\varepsilon_a<\varepsilon_c 
  < \varepsilon_b} (-1)^{[b][c]} \xi_b \xi_c q^{(\rho, \varepsilon_c - 
  \varepsilon_b)} X^a_c {A^{(l)}}^b_{\overline{c}}.
\end{align*}

\noindent Simplifying gives

\begin{align*}
(-1)^{[a]+[b]} q^{(\varepsilon_a,\varepsilon_b)} X^a_b &{A^{(l)}}^b_a -
  q^{(\varepsilon_b, \varepsilon_b)} {A^{(l)}}^b_a X^a_b \notag \\
= &\bigl( (-1)^{[a]} - \delta^a_{\overline{b}} q^{(\rho, \varepsilon_a - 
  \varepsilon_b)} \bigr) q^{h_{\varepsilon_a}} ({A^{(l)}}^a_a - {A^{(l)}}^b_b) 
  \notag \\
& \qquad +(q-q^{-1})\sum_{\varepsilon_a< \varepsilon_c < \varepsilon_b} \bigl( 
  (-1)^{[c]} - \delta^b_{\overline{c}} q^{(\rho, \varepsilon_c -\varepsilon_b)}
  \bigr) {A^{(l)}}^c_a X^a_c \notag \\
& \qquad +(q-q^{-1}) \delta^a_{\overline{b}} \sum_{\varepsilon_a< \varepsilon_c
  < \varepsilon_b} (-1)^{[b][c]} \xi_b \xi_c q^{(\rho, \varepsilon_c - 
  \varepsilon_b)} X^a_c {A^{(l)}}^{\overline{a}}_{\overline{c}}.
\end{align*}

\noindent Remembering that $\varepsilon_a < \varepsilon_b$, we apply this to the highest weight state $|\Lambda \rangle$ to obtain

\begin{multline} \label{AX}
-q^{(\varepsilon_b,\varepsilon_b)} {A^{(l)}}^b_a X^a_b |\Lambda \rangle
=   q^{(\Lambda, \varepsilon_a)} \bigl( (-1)^{[a]} - \delta^a_{\overline{b}}  
   q^{2(\rho, \varepsilon_a)} \bigr) (t_a^{(l)} - t_b^{(l)}) |\Lambda \rangle\\
+ (q-q^{-1})\sum_{\varepsilon_a < \varepsilon_c < \varepsilon_b}
   \bigl( (-1)^{[c]} - \delta^b_{\overline{c}} q^{2(\rho, \varepsilon_c)}\bigr)
  {A^{(l)}}^c_a X^a_c |\Lambda \rangle.
\end{multline}


\noindent The next step is to calculate ${A^{(l)}}^b_a X^a_b |\Lambda \rangle$ for $\varepsilon_a < \varepsilon_b$.  It is first convenient to order the indices according to $b>c \Leftrightarrow \varepsilon_b < \varepsilon_c$.  With this ordering we say an element $a >0$ if $\varepsilon_a <0$, $a=0$ if $\varepsilon_a =0$, and $a<0$ if $\varepsilon_a >0$.  Using this convention, it is apparent the solution to \eqref{AX} will be of the form

\begin{equation} \label{soln}
{A^{(l)}}^b_a X^a_b |\Lambda \rangle =  q^{(\Lambda, \varepsilon_a)} (-1)^{[a]}
   \sum_{a>c \geq b} \alpha^a_{bc} (t_a^{(l)} - t_c^{(l)}) |\Lambda \rangle,
\end{equation}

\noindent where $\alpha^a_{bc}$ is a function of $a,b$ and $c$. Now from equation \eqref{AX} we have

\begin{align*}
(q-q^{-1}) \sum_{a>c>b}&(-1)^{[c]}{A^{(l)}}^c_a X^a_c |\Lambda \rangle \notag\\
&=  -q^{(\varepsilon_b,\varepsilon_b)} {A^{(l)}}^b_a X^a_b |\Lambda \rangle + 
  (q-q^{-1}) \sum_{a>c>b} \delta^b_{\overline{c}} q^{-2(\rho, \varepsilon_b)} 
  {A^{(l)}}^c_a X^a_c  |\Lambda \rangle \notag \\
& \qquad - (-1)^{[a]} q^{(\Lambda, \varepsilon_a)} \bigl( 1 - \delta^a_
  {\overline{b}} (-1)^{[a]} q^{2(\rho, \varepsilon_a)} \bigr) (t_a^{(l)} - 
  t_b^{(l)}) |\Lambda \rangle \notag \\
&= -q^{(\varepsilon_{b+1},\varepsilon_{b+1})} {A^{(l)}}^{b+1}_a X^a_{b+1} 
  |\Lambda \rangle \notag \\
& \qquad + (q-q^{-1}) \sum_{a>c>b+1} \delta^{b+1}_{\overline{c}} 
  q^{-2(\rho, \varepsilon_{b+1})} {A^{(l)}}^c_a X^a_c |\Lambda \rangle \notag\\
& \qquad - (-1)^{[a]} q^{(\Lambda, \varepsilon_a)} \bigl( 1 - \delta^{b+1}_{ 
  \overline{a}} (-1)^{[a]} q^{2(\rho, \varepsilon_a)} \bigr) 
  (t_a^{(l)} - t_{b+1}^{(l)}) |\Lambda \rangle \notag \\
& \qquad + (q-q^{-1}) (-1)^{[b+1]} {A^{(l)}}^{b+1}_a X^a_{b+1}|\Lambda \rangle.
\end{align*} 

\noindent Substituting in the form of the solution given in equation \eqref{soln} produces

\begin{align}
q^{(\varepsilon_b,\varepsilon_b)} &\sum_{a>d\geq b} \alpha^a_{bd} (t_a^{(l)} - 
  t_d^{(l)}) |\Lambda \rangle \notag \\
= \bigl( &q^{(\varepsilon_{b+1},\varepsilon_{b+1})} - (q-q^{-1}) 
  (-1)^{[b+1]} \bigr) \sum_{a>d \geq b+1} \alpha^a_{b+1\, d} (t_a^{(l)} - 
  t_d^{(l)}) |\Lambda \rangle \notag \\
-&\bigl(1 -\delta^a_{\overline{b}} (-1)^{[a]} q^{2(\rho,\varepsilon_a)}
  \bigr) (t_a^{(l)} - t_b^{(l)}) |\Lambda \rangle + \bigl(1 - \delta^a_
  {\overline{b+1}} (-1)^{[a]} q^{2(\rho,\varepsilon_a)} \bigr) (t_a^{(l)} - 
  t_{b+1}^{(l)}) |\Lambda \rangle \notag \\
+& (q-q^{-1}) \sum_{a>c>b} \delta^b_{\overline{c}} q^{-2(\rho, 
  \varepsilon_b)} \sum_{a>d \geq c} \alpha^a_{\overline{b}\,d} (t_a^{(l)} - 
  t_d^{(l)}) |\Lambda \rangle \notag \\
-& (q-q^{-1}) \sum_{a>c>b+1} \delta^{b+1}_{\overline{c}} q^{-2(\rho,
  \varepsilon_{b+1})}  \sum_{a>d \geq c} \alpha^a_{\overline{b+1}\,d} 
  (t_a^{(l)} - t_d^{(l)}) |\Lambda \rangle . \label{qaa}
\end{align}

\newpage

\noindent Set 

\begin{equation*}
\alpha^a_{bd} = \bar{\alpha}_{bd} (1-\delta^a_{\overline{d}} (-1)^{[a]}
  q^{2(\rho, \varepsilon_a)}).
\end{equation*}

\noindent Then from equation \eqref{qaa} we obtain

\begin{equation*}
\bar{\alpha}_{bb} = - q^{-(\varepsilon_b,\varepsilon_b)} 
\end{equation*}

\noindent and

\begin{align*}
\bar{\alpha}_{b\, b+1} &=  q^{-(\varepsilon_b,\varepsilon_b)} \Bigl[ \bigl( 
  q^{(\varepsilon_{b+1},\varepsilon_{b+1})} - (q-q^{-1}) (-1)^{[b+1]} \bigr) 
  \bar{\alpha}_{b+1\,b+1} +1 \\
&\hspace{7cm} + (q-q^{-1}) \delta^b_{\overline{b+1}} 
  q^{-2(\rho, \varepsilon_b)} \bar{\alpha}_{\overline{b}\,b+1} \Bigr] \\
&=  q^{-(\varepsilon_b,\varepsilon_b)} \Bigl[ \bigl( (q-q^{-1}) (-1)^{[b+1]}
  -q^{(\varepsilon_{b+1},\varepsilon_{b+1})} \bigr) q^{-(\varepsilon_{b+1},
  \varepsilon_{b+1})} +1 \\
& \hspace{7cm} -(q-q^{-1}) \delta^b_{\overline{b+1}} q^{-2(\rho,\varepsilon_b)}
   q^{-(\varepsilon_{b+1},\varepsilon_{b+1})}\Bigr] \\
&=  q^{-(\varepsilon_b,\varepsilon_b)-(\varepsilon_{b+1},\varepsilon_{b+1})}
  (q-q^{-1}) \bigl( (-1)^{[b+1]} - \delta^b_{\overline{b+1}} 
  q^{-2(\rho, \varepsilon_b)} \bigr).
\end{align*}

\

\noindent To simplify this expression note that $q^{2(\rho, \varepsilon_{b+1} - \varepsilon_b)} = q^{-(\varepsilon_b,\varepsilon_b) - (\varepsilon_{b+1}, \varepsilon_{b+1})}$ in all cases except for $[b]=0, \,b=l,\, m=2l$, in which case $q^{2(\rho, \varepsilon_{b+1} - \varepsilon_b)} = q^2 q^{-(\varepsilon_b,\varepsilon_b) - (\varepsilon_{b+1}, \varepsilon_{b+1})}$.  However $[b]=0,\,b=l,\, m=2l$ if and only if $\delta^b_{\overline{b+1}}=1$, and in that case we find $\bar{\alpha}_{b\, b+1}=0$.  Hence for all values of $b$ we can write

\begin{equation*}
\bar{\alpha}_{b\, b+1}= (q-q^{-1}) q^{-2(\rho,\varepsilon_b)} \bigl(  
  (-1)^{[b+1]}q^{2(\rho,\varepsilon_{b+1})} - \delta^b_{\overline{b+1}} \bigr).
\end{equation*}

\noindent Now that we have found $\bar{\alpha}_{bb}$ and $\bar{\alpha}_{b\,b+1}$, they can be used to calculate the remaining $\bar{\alpha}_{bd}$.  From equation \eqref{qaa} we observe that if $d>b+1$ then

\begin{align*}
\bar{\alpha}_{bd} &= q^{-(\varepsilon_b,\varepsilon_b)} \bigl 
  (q^{(\varepsilon_{b+1}, \varepsilon_{b+1})} - (q-q^{-1}) (-1)^{[b+1]} \bigr) 
  \bar{\alpha}_{b+1\, d} \\
& \qquad + (q-q^{-1}) q^{-(\varepsilon_b, \varepsilon_b)} \sum_{d \geq c >b} 
  \delta^b_{\overline{c}} q^{-2(\rho, \varepsilon_b)} \bar{\alpha}_{\overline
  {b} d}\\
& \qquad - (q-q^{-1}) q^{-(\varepsilon_b, \varepsilon_{b})} \sum_{d \geq c 
  >b+1} \delta^{b+1}_{\overline{c}} q^{-2(\rho, \varepsilon_{b+1})} 
  \bar{\alpha}_{\overline{b+1} d}.
\end{align*}

\

\noindent  Remembering that $\theta_{xy}$ was defined by

\begin{equation*}
\theta_{xy} = 
\begin{cases} 1 \qquad &x<y,\\
0 & x \geq y,
\end{cases}
\end{equation*}

\noindent this can be rewritten as 

\begin{multline} \label{bad}
\bar{\alpha}_{bd} = q^{-(\varepsilon_b,\varepsilon_b)} \bigl( q^{(\varepsilon
  _{b+1}, \varepsilon_{b+1})} - (q-q^{-1}) (-1)^{[b+1]} \bigr) \bar{\alpha}_
  {b+1\,d} \\
+ (q-q^{-1}) q^{-(\varepsilon_b, \varepsilon_b)}  q^{2(\rho, \varepsilon_c)} 
  \bigl( \theta_{bc} \theta_{c\,d+1} \delta^b_ {\overline{c}} - \theta_
  {b+1\,c} \theta_{c\,d+1} \delta^{b+1}_{\overline{c}} \bigr) \bar{\alpha}_
  {cd}, \quad d > a+1.
\end{multline}

\noindent Now consider $\bar{\alpha}_{bd}$ for any $b >l$.  Both $\theta_{b\overline{b}}$ and $\theta_{b+1\,\overline{b+1}}$ will equal $0$, so

\begin{align*}
\bar{\alpha}_{bd} &= q^{-(\varepsilon_b,\varepsilon_b)} \bigl(q^{(\varepsilon
  _{b+1}, \varepsilon_{b+1})} - (q-q^{-1}) (-1)^{[b+1]} \bigr) \bar{\alpha}
  ^a_{b+1\, d} \\
&= q^{-(\varepsilon_b, \varepsilon_b)} q^{-(\varepsilon_{b+1}, \varepsilon_
  {b+1})} \bar{\alpha}_{b+1\, d}  \\
&= q^{2(\rho, \varepsilon_{b+1} - \varepsilon_b)} \bar{\alpha}_{b+1\, d}.
\end{align*}

\noindent Since 

\begin{equation*}
\bar{\alpha}_{d-1\,d} =(-1)^{[d]} (q-q^{-1}) q^{2(\rho, \varepsilon_d-
  \varepsilon_{d-1})},
\end{equation*}

\noindent we obtain

\vspace{-1cm}

\begin{equation*} 
\bar{\alpha}_{bd} =(-1)^{[d]} (q-q^{-1}) q^{2(\rho, \varepsilon_{d}-
  \varepsilon_b)}, \quad d>b>l.
\end{equation*}

\noindent Substituting this together with our expression for $\bar{\alpha}_{bb}$ into equation \eqref{bad}, we find 

\begin{align}
\bar{\alpha}_{bd} = &q^{-(\varepsilon_b,\varepsilon_b)} \bigl( 
  q^{-(\varepsilon_{b+1}, \varepsilon_{b+1})} - \delta^{b+1}_{\overline{b+1}}
  (q-q^{-1}) \bigr) \bar{\alpha}_{b+1\, d} \notag \\
+& (q-q^{-1})^2 q^{-(\varepsilon_b, \varepsilon_b)} (-1)^{[d]} 
  q^{2(\rho, \varepsilon_d)} \bigl( \theta_{b\overline{b}} \theta_{\overline
  {b}d} - \theta_{b+1\,\overline{b+1}} \theta_{\overline{b+1}d} \bigr)\notag\\
-& (q-q^{-1}) q^{-(\varepsilon_b, \varepsilon_b)} q^{-(\varepsilon_d, 
  \varepsilon_d)} q^{2(\rho, \varepsilon_d)} ( \delta^{\overline{b}}_d - 
  \delta^{\overline{b+1}}_d), \; d > b+1. \label{um}
\end{align}
\newpage
\noindent But for $d>b+1$ 

\begin{align*}
\theta_{b\overline{b}} \theta_{\overline{b}d} - \theta_{b+1\,\overline{b+1}} 
  \theta_{\overline{b+1}\,d} &= \delta^b_l \theta_{\overline{l}d} - 
  \delta^{\overline{b}}_d \theta_{bl} \notag \\
&= \delta^b_l (1-\delta^d_{\overline{l}}) -  \delta^{\overline{b}}_d 
  (1- \delta^b_l) \notag \\
&= \delta^b_l -\delta^{\overline{b}}_d.
\end{align*}

\noindent Also, $-[(-1)^{[d]} (q-q^{-1}) + q^{-(\varepsilon_d,\varepsilon_d)}] \delta^{\overline{b}}_d = - q^{(\varepsilon_d,\varepsilon_d)} \delta^{\overline{b}}_d$, so equation \eqref{um} reduces to

\begin{align*}
\bar{\alpha}_{bd} &= \bigl( q^{2(\rho, \varepsilon_{b+1} - \varepsilon_b)} 
  q^{-2\delta^b_{\overline{b+1}}} -\delta^{b+1}_{\overline{b+1}} q^{-1}
  (q-q^{-1}) \bigr) \bar{\alpha}_{b+1\,d} + \delta^b_l q^{-1}(q-q^{-1})^2 
  (-1)^{[d]} q^{2(\rho, \varepsilon_d)}  \\
& \quad - \delta^{\overline{b}}_d  (q-q^{-1}) q^{2(\rho, \varepsilon_d)} + 
  \delta^{\overline{b+1}}_d (q-q^{-1}) q^{2(\rho, \varepsilon_{b+1}-
  \varepsilon_{b})} q^{-2\delta^b_{\overline{b+1}}} q^{2(\rho, \varepsilon_d)}
  \\
&= \bigl( q^{2(\rho, \varepsilon_{b+1} - \varepsilon_b)} 
  q^{-2\delta^b_{\overline{b+1}}} -\delta^{b+1}_{\overline{b+1}} q^{-1}
  (q-q^{-1}) \bigr) \bar{\alpha}_{b+1\,d} + \delta^b_l q^{-1}(q-q^{-1})^2 
  (-1)^{[d]} q^{2(\rho, \varepsilon_d)} \\
& \quad + (q-q^{-1}) q^{-2(\rho, \varepsilon_b)} (\delta^{\overline{b+1}}_d -
  \delta^{\overline{b}}_d) ,\qquad d>b+1 .
\end{align*}

\noindent Recall that for $b > l$ we have

\begin{equation*}
\bar{\alpha}_{bd} = (-1)^{[d]} (q-q^{-1}) q^{2(\rho,\varepsilon_d -\varepsilon_b)}, \quad d>b.
\end{equation*}

\noindent  Then when $b=l$ we find

\begin{align*}
\bar{\alpha}_{bd} &= \bigl(q^{2(\rho, \varepsilon_{b+1} - \varepsilon_b)} 
  q^{-2\delta^b_{\overline{b+1}}} - \delta^{b+1}_{\overline{b+1}} q^{-1}
  (q-q^{-1}) \bigr) (-1)^{[d]} (q-q^{-1}) q^{2(\rho,\varepsilon_{d}- 
  \varepsilon_{b+1})} \\
& \quad + q^{-1}(q-q^{-1})^2 (-1)^{[d]} q^{2(\rho, \varepsilon_d)} -
  (q-q^{-1}) q^{-2(\rho, \varepsilon_b)} \delta^{\overline{l}}_d \notag \\
&= (-1)^{[d]} (q-q^{-1}) q^{2(\rho,\varepsilon_{d}- \varepsilon_{b})}\notag\\
&\hspace{2cm} \Bigl[ \delta^{b+1}_{\overline{b+1}} \bigl( 1 - (q-q^{-1}) + 
  (q-q^{-1}) \bigr)  + \delta^b_{\overline{b+1}} \bigl( q^{-2} + q^{-1} 
  (q-q^{-1}) \Bigr] \notag \\
& \quad - (q-q^{-1})  q^{-2(\rho, \varepsilon_b)} \delta_d^{\overline{l}} \\
&= (q-q^{-1}) q^{-2(\rho, \varepsilon_{b})} \bigl( (-1)^{[d]} 
  q^{2(\rho,\varepsilon_{d})} - \delta^b_{\overline{d}} \bigr)
\end{align*}

\noindent for all $d > b +1$.  Comparing this with our earlier results for $d = b+1$ and $b>l$, we have

\begin{equation*}
\bar{\alpha}_{bd} = (q-q^{-1}) q^{-2(\rho, \varepsilon_{b})} \bigl( 
  (-1)^{[d]} q^{2(\rho,\varepsilon_{d})} - \delta^b_{\overline{d}} \bigr), 
  \quad \forall b \geq l,\; d > b.
\end{equation*}

\noindent But for $b < l$ we know 

\begin{align*}
\bar{\alpha}_{bd} &= q^{2(\rho, \varepsilon_{b+1}-\varepsilon_b)} 
  \bar{\alpha}_{b+1\,d} + (q-q^{-1}) q^{-2(\rho,\varepsilon_b)}
  (\delta^{\overline{b+1}}_d - \delta^{\overline{b}}_d), \qquad d>b+1 .
\end{align*}

\noindent Hence for all $b$ we obtain

\begin{align*}
\bar{\alpha}_{bd} &= (q-q^{-1}) q^{-2(\rho, \varepsilon_{b})} \bigl( 
  (-1)^{[d]} q^{2(\rho,\varepsilon_{d})} - \sum_{c=b}^{d-1} \delta^{\overline
  {c}}_d + \sum_{c=b}^{d-2} \delta^{\overline{c+1}}_d \bigr)  \\
&= (q-q^{-1}) q^{-2(\rho, \varepsilon_{b})} \bigl( (-1)^{[d]} 
  q^{2(\rho,\varepsilon_{d})} - \delta^{\overline{b}}_d \bigr), \qquad d>b.
\end{align*}

\noindent Thus for all $a>b$

\begin{equation*}
{A^{(l)}}^b_a X^a_b |\Lambda \rangle =  q^{(\Lambda, \varepsilon_a)} (-1)^{[a]}
   \sum_{a>c \geq b} \alpha^a_{bc} (t_a^{(l)} - t_c^{(l)}) |\Lambda \rangle,
\end{equation*}

\noindent  where $\alpha^a_{bc}$ is given by

\begin{equation*}
\alpha^a_{bc} = 
  \begin{cases}
  - q^{-(\varepsilon_b,\varepsilon_b)} (1 -\delta^a_{\overline{b}} (-1)^{[a]} 
    q^{2(\rho, \varepsilon_a)}), & c=b, \\
  (q-q^{-1}) q^{-2(\rho, \varepsilon_{b})} \bigl( (-1)^{[c]} q^{2(\rho,
  \varepsilon_{c})} - \delta^{\overline{b}}_c \bigr)(1-\delta^a_{\overline{c}} 
  (-1)^{[a]} q^{2(\rho, \varepsilon_a)}), \quad & c>b.
  \end{cases}
\end{equation*}

\section{Constructing the Perelomov-Popov \\
         Matrix Equation}

\noindent This expression can now be substituted into the equation

\begin{equation*}
t^{(l)}_a |\Lambda \rangle = t^{(l-1)}_a t^{(1)}_a |\Lambda \rangle + 
  \sum_{\varepsilon_a<\varepsilon_b} (-1)^{[a]+[b]} q^{(\Lambda,\varepsilon_a)}
  {A^{(l-1)}}^b_a X^a_b |\Lambda \rangle
\end{equation*}

\noindent to find a matrix equation for the various $t_a^{(l)}$.  The matrix factor is an analogue of the Perelomov-Popov matrix introduced in \cite{Popov} and \cite{Perelomov}, which has been used to calculate the eigenvalues of the Casimir invariants of various classical Lie algebras.

\

\noindent First recall that

\begin{equation*}
A^b_a = (1+\delta^a_b) X^b_a + (q-q^{-1}) \sum_{c \geq a,b} (-1)^{([a]+[c])([b]+[c])} X^c_a X^b_c,
\end{equation*}

\noindent where

\begin{equation*}
X^b_a =
\begin{cases}
\frac{q^{h_{\varepsilon_a}} -I}{q-q^{-1}}, &a=b, \\
(-1)^{[b]} q^{h_{\varepsilon_a}} \hat{\sigma}_{ba}, \quad &\varepsilon_a <
  \varepsilon_b, \\
(-1)^{[b]} \hat{\sigma}_{ba} q^{h_{\varepsilon_b}}, &\varepsilon_a > 
  \varepsilon_b.
\end{cases}
\end{equation*}

\noindent Then

\begin{alignat*}{2}
&&\;A^a_a |\Lambda \rangle&= 2 X^a_a |\Lambda \rangle + (q-q^{-1}) X^a_a X^a_a 
  |\Lambda \rangle \notag \\
&&&= (q-q^{-1})^{-1} (2(q^{h_{\varepsilon_a}}-1) + (q^{h_{\varepsilon_a}}-1)^2)
   |\Lambda\rangle. \\
&\therefore & t^{(1)}_a & =\frac{q^{2(\Lambda, \varepsilon_a)}-1}{q-q^{-1}}.
\end{alignat*}

\noindent Hence we obtain

\begin{align*}
t^{(l)}_a &= \frac{(q^{2(\Lambda, \varepsilon_a)} -1)}{(q-q^{-1})} t^{(l-1)}_a 
  \notag \\
&\qquad+ \sum_{b<a} (-1)^{[a]+[b]} q^{(\Lambda,\varepsilon_a)} \bigl(q^
  {(\Lambda,\varepsilon_a)} (-1)^{[a]} \sum_{b \leq c < a} \alpha^a_{bc} 
  (t^{(l-1)}_a - t^{(l-1)}_c) \bigr) \notag \\
&= \frac{(q^{2(\Lambda, \varepsilon_a)} -1)}{(q-q^{-1})} t^{(l-1)}_a \notag\\
&\qquad -q^{2(\Lambda,\varepsilon_a)} \sum_{b<a} (-1)^{[b]} q^{-(\varepsilon_b,
  \varepsilon_b)} (1-\delta^a_{\overline{b}}(-1)^{[a]} q^{2(\rho, 
  \varepsilon_a)}) (t^{(l-1)}_a - t^{(l-1)}_b) \notag \\
& \qquad + (q-q^{-1}) q^{2(\Lambda,\varepsilon_a)} \sum_{c < b<a} (-1)^{[c]} 
  q^{-2(\rho,\varepsilon_{c})} (1-\delta^a_{\overline{b}} (-1)^{[a]} q^{2(\rho,
  \varepsilon_a)}) \notag \\
& \hspace{6cm} ((-1)^{[b]} q^{2(\rho,\varepsilon_{b})}-\delta^b_{\overline{c}})
  (t^{(l-1)}_a - t^{(l-1)}_b) .
\end{align*}

\noindent Now consider the function $\gamma_b$ defined by:

\begin{equation*}
\gamma_b = (-1)^{[b]} q^{-(\varepsilon_b,\varepsilon_b)} - (q-q^{-1}) \sum_{c
  < b} (-1)^{[c]} q^{-2(\rho, \varepsilon_{c})} \bigl( (-1)^{[b]} 
  q^{2(\rho,\varepsilon_{b})} - \delta^b_{\overline{c}} )\bigr).
\end{equation*}

\noindent  We evaluate this for all $b$, labelling $C(\Lambda_0) = (\delta_1, \delta_1+2\rho) = m-n-1$ and remembering that $$\rho = \frac{1}{2} \sum_{i=1}^l (m-2i) \varepsilon_i + \frac{1}{2} \sum_{\mu=1}^k (n-m+2-2\mu) \delta_\mu.$$

\

\noindent  First we consider the case $[b]=1$ and $b \leq k$. Here $\gamma_b$ is given by

\begin{align*}
\gamma_b &= -q - (q-q^{-1}) q^{-(n-m+2-2b)}\sum_{c < b} q^{n-m+2-2c} \\
&= -q - (q-q^{-1}) (q^2-1)^{-1} (q^{2b} -  q^{2}) \\
&= (-1)^{[b]} q^{2(\rho, \varepsilon_b)} q^{-C(\Lambda_0)}.
\end{align*} 

\noindent For $[b]=0,\, b \leq l$, we obtain
\begin{align*}
\gamma_b &= q^{-1} - (q-q^{-1}) q^{m-2b} \Bigl[ (q^2-1)^{-1} (q^{2-m}- 
  q^{n-m+2}) + \sum_{\varepsilon_1 \geq \varepsilon_c > \varepsilon_b} 
  q^{2c-m} \Bigr]  \\
&= q^{-1} - q^{-1} (q^{2-2b}- q^{n+2-2b}+1 - q^{2-2b})\\
&= (-1)^{[b]} q^{2(\rho, \varepsilon_b)} q^{-C(\Lambda_0)}.
\end{align*}

\noindent In the case $b=l+1, \, m=2l+1,$ we find

\begin{align*}
\gamma_b &= 1 - q^{-1} (-q^{n-m+2} +q^{2l+2-m}) \\
&=  (-1)^{[b]} q^{2(\rho, \varepsilon_b)} q^{-C(\Lambda_0)}.
\end{align*}

\noindent If $[b]=0$ and $b \geq \overline{l}$ then
\begin{align*}
\gamma_b &= q^{-1} - (q-q^{-1}) q^{2\overline{b}-m}  \Bigl[ (q^2-1)^{-1} 
  ( q^{2l+2-m} -q^{n-m+2}) + (m-2l) + \sum_{b>i \geq \overline{l}} 
  q^{m-2\overline{i}} \Bigr]  \\
& \hspace{2cm} + (q-q^{-1}) q^{2(\rho, \varepsilon_b)} \\
&= q^{-1} - q^{-1} q^{2\overline{b}-m} \bigl[ q^{2l+2-m}-q^{n-m+2} + 
  q^{m-2\overline{b}} - q^{m-2l} + (m-2l)(q^2-1) \bigr] \\ 
& \hspace{2cm} +(q-q^{-1}) q^{2\overline{b}-m} \\
&= q^{-1} - q^{-1}(1 -q^{n-2m+2 +2\overline{b}}) + (q-q^{-1}) (q^{2\overline{b}
  -m} - q^{2\overline{b}-m}) \\
&= (-1)^{[b]} q^{2(\rho, \varepsilon_b)} q^{-C(\Lambda_0)}.
\end{align*}

\noindent Lastly, when $[b]=1$ and $b \geq \overline{k}$ we find

\begin{align*}
\gamma_b &= -q + (q-q^{-1}) q^{n-m+2-2\overline{b}} \Bigl[ (q^2-1)^{-1} (q^m 
  - q^{n-m+2} + q^2-1) - \sum_{b>c \geq \overline{k}} q^{m-n-2+2\overline{c}} 
  \Bigr] \\
& \hspace{2cm} - (q-q^{-1}) q^{n-m+2-2\overline{b}} \\
&= -q + q^{-1} q^{n-m+2-2\overline{b}} ( q^m - q^{n-m+2} -
  q^{m-n+2k} + q^{m-n +2\overline{b}}) \\
&= -q^{2n-2m+3-2\overline{b}} \\
&= (-1)^{[b]} q^{2(\rho, \varepsilon_b)} q^{-C(\Lambda_0)}.
\end{align*}

\noindent Hence for all $b$

\begin{equation*}
\gamma_b = (-1)^{[b]} q^{2(\rho, \varepsilon_b)} q^{-C(\Lambda_0)}.
\end{equation*}

\noindent We also consider the function

\begin{equation*}
\beta_a = 1 - (q-q^{-1}) \sum_{b<a} \gamma_b \bigl( 1 - \delta^a_{\overline{b}}
  (-1)^{[a]} q^{2(\rho, \varepsilon_a)}\bigr),
\end{equation*}

\noindent so that

\begin{equation} \label{tal}
t_a^{(l)} = \frac{(q^{2(\Lambda, \varepsilon_a)} \beta_a -1)}{(q-q^{-1})}
  t_a^{(l-1)} + q^{2(\Lambda, \varepsilon_a)} \sum_{b<a} \gamma_b \bigl( 1-
  \delta^a_{\overline{b}} (-1)^{[a]} q^{2(\rho, \varepsilon_a)}\bigr) 
  t_b^{(l-1)}.
\end{equation}

\noindent As before, we evaluate this by considering the various cases individually.  Firstly, for $[a]=1,\,a \leq k$ we find

\begin{align*}
\beta_a &=1 +(q-q^{-1}) \sum_{b<a} q^{2b-2+m-n} q^{n-m+1} \\
&= 1+q^{-1}(q^{2a-2+m-n}-q^{m-n}) q^{n-m+1} \\
&= q^{(\varepsilon_a, 2\rho +\varepsilon_a) - C(\Lambda_0)}.
\end{align*}

\noindent  Using this, when $[a]=0$ and $a \leq l$ we obtain

\begin{align*}
\beta_a &= q^{2(k+1)-2} - (q-q^{-1}) \sum_{\varepsilon_1 \geq \varepsilon_b > 
  \varepsilon_a} q^{m-2b} q^{n-m+1} \\
&= q^{n} - q^{-1} q^{n-m+1} (q^{m} - q^{m-2a+2}) \\
&= q^{(\varepsilon_a, 2\rho +\varepsilon_a) - C(\Lambda_0)}.
\end{align*}

\noindent In the case $a=l+1,\, m=2l+1,$ we find

\begin{align*}
\beta_a &= q^{m-2(l+1)+1}q^{n-m+1} \\
&= q^{(\varepsilon_a, 2\rho +\varepsilon_a) - C(\Lambda_0)}.
\end{align*}

\noindent  For the remaining values of $a$, we consider the cases $m=2l$ and $m=2l+1$ separately.  Firstly, take $m=2l$.  Then for $[a]=0,\, a \geq \overline{l}$, we obtain

\begin{align*}
\beta_a &= q^{m-2(l+1)+1} q^{- C(\Lambda_0)} + (q-q^{-1}) \gamma_{\overline{a}}
  q^{2(\rho, \varepsilon_a)} - (q-q^{-1}) \sum_{b=\overline{l}}^{a-1} 
  q^{2\overline{b}-m} q^{- C(\Lambda_0)} \\
&= q^{- C(\Lambda_0)} [q^{-1} + (q-q^{-1}) - q^{-1}(q^{2l+2-m} - 
  q^{2\overline{a}+2-m})] \\
&= q^{(\varepsilon_a, 2\rho +\varepsilon_a) - C(\Lambda_0)},
\end{align*}

\noindent whereas for the remaining odd values of $a$, namely $a \geq \overline{k}$, $\beta_a$ is given by

\begin{align*}
\beta_a &=q^{- C(\Lambda_0)} \bigl( q^{-1} - q^{-1}(q^{2} - q^{2-m}) \bigr) - 
  (q-q^{-1}) \gamma_{\overline{a}} q^{2(\rho, \varepsilon_a)} \\
& \qquad + (q-q^{-1}) \sum_{b=\overline{k}}^{a-1} 
  q^{n-m+2-2\overline{b}} q^{-C(\Lambda_0)}\\
&= q^{- C(\Lambda_0)} [q^{-1} - q^{-1}(q^{2}-q^{2-m})+(q-q^{-1}) + 
  q^{-1}(q^{n-m+2-2\overline{a}} -q^{2-m})] \\
&=q^{(\varepsilon_a, 2\rho +\varepsilon_a) - C(\Lambda_0)}.
\end{align*}

\noindent  Now consider $m=2l+1$.  If $[a]=0$ and $a > l+1$ then

\begin{align*}
\beta_a &= q^{- C(\Lambda_0)} -(q-q^{-1})q^{- C(\Lambda_0)} +(q-q^{-1})
  \gamma_{\overline{a}} q^{2(\rho, \varepsilon_a)} -\sum_{b=\overline{l}}^{a-1}
  q^{2\overline{b}-m} q^{- C(\Lambda_0)} \\
&= q^{- C(\Lambda_0)} [1- (q-q^{-1}) + (q-q^{-1}) -q^{-1}(q^{2l+2-m}-
  q^{2\overline{a}+2-m})] \\
&= q^{(\varepsilon_a, 2\rho +\varepsilon_a) - C(\Lambda_0)}.
\end{align*}

\noindent For the remaining case of $[a]=1,\, a \geq \overline{k},$ we find

\begin{align*}
\beta_a &=q^{- C(\Lambda_0)}[1-(q-q^{-1})-q^{-1}(q-q^{2-m})] - (q-q^{-1}) 
  \gamma_{\overline{a}} q^{2(\rho, \varepsilon_a)} \\
& \qquad + (q-q^{-1}) \sum_{b=\overline{k}}^{a-1} q^{n-m+2-2\overline{b}} 
  q^{-C(\Lambda_0)}  \\
&= q^{- C(\Lambda_0)} [1- (q-q^{-1}) -q^{-1}(q-q^{2-m})  +(q-q^{-1}) + q^{-1}
  (q^{n-m+2-2\overline{a}} -q^{2-m})] \\
&= q^{(\varepsilon_a, 2\rho +\varepsilon_a) - C(\Lambda_0)}.
\end{align*}

\noindent Hence $\beta_a$ is given by

\begin{equation*}
\beta_a= q^{(\varepsilon_a, 2\rho +\varepsilon_a) - C(\Lambda_0)}
\end{equation*}

\noindent for any $a$, regardless of the parity of $m$.  Substituting this result together with that for $\gamma_b$ into equation (\ref{tal}) gives

\begin{multline*}
t_a^{(l)} = \frac{(q^{(\varepsilon_a,2\Lambda+ 2\rho +\varepsilon_a) - 
  C(\Lambda_0)}-1)}{(q-q^{-1})} t_a^{(l-1)} \\
  + q^{(2\Lambda, \varepsilon_a)-C(\Lambda_0)} \sum_{b<a} (-1)^{[b]} 
  q^{(2\rho, \varepsilon_b)} \bigl( 1-\delta^a_{\overline{b}} (-1)^{[a]} 
  q^{(2\rho, \varepsilon_a)}\bigr) t_b^{(l-1)}.
\end{multline*}

\noindent  This can be written in the matrix form

\begin{equation*}
\underline{t}^{(l)} = M \underline{t}^{(l-1)},
\end{equation*}

\noindent where $M$ is a lower triangular matrix with entries

\begin{equation*}
M_{ab} = 
\begin{cases}
0, &a<b \\
(q-q^{-1})^{-1}(q^{(\varepsilon_a,2\Lambda+ 2\rho +\varepsilon_a) - 
  C(\Lambda_0)}-1), &a=b, \\
 q^{(2\Lambda, \varepsilon_a)-C(\Lambda_0)} \bigl( (-1)^{[b]} q^{(2\rho, 
  \varepsilon_b)} -\delta^a_{\overline{b}} \bigr),  &a>b.
\end{cases}
\end{equation*}

\noindent Then we have

\begin{equation*}
\underline{t}^{(l)} = M^{l} \underline{t}^{(0)}, \qquad \text{where} \quad 
  t^{(0)}_a =1 \quad \forall a,
\end{equation*}

\noindent where $M$ is an analogue of the Perelomov-Popov matrix.

\section{Finding the Eigenvalues}

\noindent  This matrix equation for $t_a^{(l)}$ can now be used to calculate the eigenvalues of $C_l$.  Loosely speaking, the problem reduces to diagonalising the matrix $M$.  As with the earlier calculations, the $q$-factors from the $q$-deformation and the $\delta$-functions arising from the root system make this somewhat more difficult than in the classical cases studied in \cite{Bincer} and \cite{Scheunert83}, and also than in the case of $U_q[gl(m|n)]$ \cite{LinksZhang}.

\

\noindent Recall

\begin{equation*}
C_l = \sum_a (-1)^{[a]} q^{(2\rho, \varepsilon_a)} {A^{(l)}}^a_a.
\end{equation*}

\noindent Denote the eigenvalue of $C_l$ on $V(\Lambda)$ as $\chi_{\Lambda}(C_l)$. Then we have

\begin{equation*}
\chi_{\Lambda}(C_l) = \sum_a (-1)^{[a]} q^{(2\rho, \varepsilon_a)} t_a^{(l)}
  = \sum_{a,b} (-1)^{[a]} q^{(2\rho, \varepsilon_a)} (M^l)_{ab}.
\end{equation*}

\noindent To calculate this we wish to diagonalise $M$.  We assume the eigenvalues of $M$,

\begin{equation*}
\alpha_a^{\Lambda} = \frac{(q^{(\varepsilon_a,2\Lambda+ 2\rho +\varepsilon_a) -
  C(\Lambda_0)}-1)}{(q-q^{-1})},
\end{equation*}

\noindent are distinct.  Then we need a matrix $N$ satisfying

\begin{equation*}
(N^{-1}MN)_{ab} = \delta^a_b \alpha_a^{\Lambda}, 
\end{equation*}

\noindent  which implies

\begin{equation} \label{chi}
\chi_{\Lambda}(C_l) = \sum_{a,b,c} (-1)^{[a]} q^{(2\rho, \varepsilon_a)} (
  \alpha_b^{\Lambda} )^l N_{ab} (N^{-1})_{bc}.
\end{equation}

\noindent Now

\begin{equation*}
(MN)_{ab} = \alpha_b^{\Lambda} N_{ab}.
\end{equation*}

\noindent Substituting in the values for $M_{ab}$ gives

\begin{equation} \label{alphaN}
\alpha_a^{\Lambda} N_{ab} + q^{(2\Lambda, \varepsilon_a)-C(\Lambda_0)} 
  \sum_{c<a} \bigl( (-1)^{[c]} q^{(2\rho, \varepsilon_c)} - 
  \delta^a_{\overline{c}} \bigr) N_{cb} = \alpha_b^{\Lambda} N_{ab}.
\end{equation}

\noindent Since the eigenvalues $\alpha_a^{\Lambda}$ are distinct, this implies

\begin{equation*}
N_{ab} = 0, \qquad \forall a < b.
\end{equation*}

\noindent Set

\begin{equation} \label{P}
P_{ab} = \sum_{c\leq a} (-1)^{[c]} q^{(2\rho, \varepsilon_c)} N_{cb}.
\end{equation}

\noindent Then equation (\ref{alphaN}) becomes

\begin{alignat*}{2}
&&&(\alpha_b^{\Lambda}- \alpha_a^{\Lambda}) N_{ab} = q^{(2\Lambda, 
  \varepsilon_a) -C(\Lambda_0)} P_{a-1\,b} - \theta_{0a} q^{(2\Lambda, 
  \varepsilon_a)-C(\Lambda_0)} N_{\overline{a}b} \notag \\
&\Rightarrow \qquad&&(\alpha_b^{\Lambda} - \alpha_a^{\Lambda})(-1)^{[a]} 
  q^{(-2\rho,\varepsilon_a)} (P_{ab}-P_{a-1\,b}) \notag \\
&&&\hspace{35mm}= q^{(2\Lambda, \varepsilon_a)-C(\Lambda_0)} P_{a-1\,b} - 
  \theta_{0a} q^{(2\Lambda, \varepsilon_a)-C(\Lambda_0)} N_{\overline{a}b},
\end{alignat*}

\noindent which simplifies to

\begin{equation*} 
P_{ab} = \frac{(\alpha_b^{\Lambda} - \alpha_a^{\Lambda} + (-1)^{[a]} 
  q^{2(\Lambda+\rho, \varepsilon_a) - C(\Lambda_0)})}{(\alpha_b^{\Lambda}- 
  \alpha_a^{\Lambda})} P_{a-1\,b} - \frac{\theta_{0a}(-1)^{[a]} q^{2(\Lambda +
  \rho, \varepsilon_a) -C(\Lambda_0)}}{(\alpha_b^{\Lambda} - 
  \alpha_a^{\Lambda})} N_{\overline{a}b}.
\end{equation*}
\

\noindent Set

\begin{equation*}
\psi^b_{a} = \alpha_b^{\Lambda} - \alpha_a^{\Lambda} + (-1)^{[a]} 
  q^{2(\Lambda+\rho, \varepsilon_a) - C(\Lambda_0)},
\end{equation*}

\noindent so this becomes

\begin{equation} \label{Pab}
P_{ab} = \frac{\psi^b_a}{(\alpha_b^{\Lambda}- \alpha_a^{\Lambda})} P_{a-1\,b} 
  - \frac{\theta_{0a} (-1)^{[a]} q^{2(\Lambda+\rho, \varepsilon_a) -
  C(\Lambda_0)}} {(\alpha_b^{\Lambda} - \alpha_a^{\Lambda})} N_{\overline{a}b}.
\end{equation}

\noindent Without loss of generality we can choose $N_{aa} = 1\; \forall a$, so $P_{bb} = (-1)^{[b]} q^{2(\rho, \varepsilon_b)}$.  Then in the cases $0 \geq a > b$ and $a > b \geq 0$ the last term in equation (\ref{Pab}) vanishes, giving

\begin{equation*}
P_{ab} = (-1)^{[b]} q^{2(\rho, \varepsilon_b)} \prod_{c=b+1}^a \frac{\psi^b_c}
  {(\alpha_b^{\Lambda}- \alpha_c^{\Lambda})}.
\end{equation*}

\noindent Similarly, for $a > \overline{b} > 0$ we obtain

\begin{equation} \label{abbar0}
P_{ab} = P_{\overline{b}b} \prod_{c=\overline{b}+1}^a \frac{\psi^b_c} 
  {(\alpha_b^{\Lambda}- \alpha_c^{\Lambda})} .
\end{equation}

\noindent It remains to find $P_{ab}$ for $\overline{b} \geq a >0$.  In this case, the last term in equation (\ref{Pab}) contributes, giving 

\begin{multline} \label{P_ab}
P_{ab} = (-1)^{[b]} q^{2(\rho, \varepsilon_b)} \prod_{c=b+1}^a \frac{\psi^b_c}
  {(\alpha_b^{\Lambda}- \alpha_c^{\Lambda})}  - \frac{(-1)^{[a]} q^{2(\Lambda+
  \rho, \varepsilon_a) -C(\Lambda_0)}}{(\alpha_b^{\Lambda}-\alpha_a^{\Lambda})}
  N_{\overline{a}b} \\
- \sum_{j=\overline{l}}^{a-1}\frac{(-1)^{[j]} q^{2(\Lambda +\rho,\varepsilon_j)
  -C(\Lambda_0)}}{(\alpha_b^{\Lambda} - \alpha_j^{\Lambda})} N_{\overline{j}b} 
  \prod_{c=j+1}^a \frac{\psi^b_c} {(\alpha_b^{\Lambda}- \alpha_c^{\Lambda})}.
\end{multline}

\noindent Recall that if $b<a<0$, then

\begin{align*}
N_{ab} &= \frac{q^{(2\Lambda,\varepsilon_a)-C(\Lambda_0)}}{(\alpha_b^{\Lambda}-
  \alpha_a^{\Lambda})} P_{a-1\, b} \notag \\
&= \frac{(-1)^{[b]} q^{2(\Lambda,\varepsilon_a)+2(\rho,\varepsilon_b)-C(
  \Lambda_0)}} {(\alpha_b^{\Lambda}- \alpha_a^{\Lambda})} \prod_{c=b+1}^{a-1} 
  \frac{\psi^b_c} {(\alpha_b^{\Lambda}- \alpha_c^{\Lambda})}.
\end{align*} 

\noindent Substituting this into equation \eqref{P_ab}, we find

\begin{align*}
P_{\overline{b}b} &= (-1)^{[b]} q^{2(\rho,\varepsilon_b)} \prod_{c=b+1}^
  {\overline{b}} \frac{\psi^b_c} {(\alpha_b^{\Lambda}- \alpha_c^{\Lambda})} - 
  \frac{(-1)^{[b]} q^{-2(\Lambda+\rho, \varepsilon_b) -C(\Lambda_0)}}
  {(\alpha_b^{\Lambda} - \alpha_{\overline{b}}^{\Lambda})} \notag \\ 
&\qquad -\sum_{j=\overline{l}}^{\overline{b}-1} \frac{(-1)^{[j]+[b]} q^{2(\rho, 
  \varepsilon_j + \varepsilon_b) - 2C(\Lambda_0)}} {(\alpha_b^{\Lambda} - 
  \alpha_j^{\Lambda}) (\alpha_b^{\Lambda} -\alpha_{\overline{j}}^{\Lambda})} 
  \prod_{c=b+1}^{\overline{j}-1} \frac{\psi^b_c} {(\alpha_b^{\Lambda}- 
  \alpha_c^{\Lambda})} \prod_{c=j+1}^{\overline{b}}\frac{\psi^b_c} 
  {(\alpha_b^{\Lambda} - \alpha_c^{\Lambda})},
\end{align*}

\noindent which can also be written as

\begin{align}
P_{\overline{b}b} \prod_{c=b+1}^{\overline{b}} \frac{(\alpha_b^{\Lambda}- 
  \alpha_c^{\Lambda})} {\psi^b_c} &= (-1)^{[b]} q^{2(\rho,\varepsilon_b)} - 
  \frac{(-1)^{[b]} q^{-2(\Lambda+\rho, \varepsilon_b) -C(\Lambda_0)}} 
  {(\alpha_b^{\Lambda} - \alpha_{\overline{b}}^{\Lambda})} \prod_{c=b+1}^
  {\overline{b}} \frac{(\alpha_b^{\Lambda}-\alpha_c^{\Lambda})} {\psi^b_c} 
  \notag \\
&\qquad -\sum_{j=\overline{l}}^{\overline{b}-1} \frac{(-1)^{[j]+[b]} q^{2(\rho,
  \varepsilon_b+\varepsilon_j)-2C(\Lambda_0)}} {(\alpha_b^{\Lambda} - \alpha_j
  ^{\Lambda}) (\alpha_b^{\Lambda} -\alpha_{\overline{j}}^{\Lambda})}\prod_{c=
  \overline{j}}^{j} \frac{(\alpha_b^{\Lambda}-\alpha_c^{\Lambda})} {\psi^b_c}
  \notag \\
&= (-1)^{[b]} q^{2(\rho,\varepsilon_b)} - \frac{q^{-2C(\Lambda_0)}} {\psi^b_
  {\overline{b}} \psi^b_b} \prod_{c=b+1}^{\overline{b}-1} \frac{(\alpha_b^
  {\Lambda}-\alpha_c^{\Lambda})} {\psi^b_c} \notag \\
&\qquad - \sum_{j=\overline{l}}^{\overline{b}-1} \frac{(-1)^{[j]+[b]} 
  q^{2(\rho, \varepsilon_b+\varepsilon_j)-2C(\Lambda_0)}} {\psi^b_j \psi^b_
  {\overline{j}}} \prod_{c=\overline{j}+1}^{j-1} \frac{(\alpha_b^{\Lambda}
  -\alpha_c^{\Lambda})} {\psi^b_c} \notag \\
&= (-1)^{[b]} q^{2(\rho,\varepsilon_b)} \Biggl[ 1-\sum_{j=\overline{l}}
  ^{\overline{b}} \frac{(-1)^{[j]} q^{2(\rho,\varepsilon_j)-2C(\Lambda_0)}}  
  {\psi^b_j  \psi^b_{\overline{j}}} \prod_{c=\overline{j}+1}^{j-1} 
  \frac{(\alpha_b^{\Lambda}-\alpha_c^{\Lambda})}{\psi^b_c}\Biggr]\label{bbbar}.
\end{align}

\noindent From this point we will consider the case $m=2l+1$.  This is marginally more complicated than the case with even $m$.

\

\noindent Define $\Phi^b_j$ to be

\begin{align*}
\Phi^b_j &= \prod_{c=\overline{l}}^{j-1} \frac{(\alpha_b -\alpha_c)(\alpha_b - 
  \alpha_{\overline{c}})}{\psi^b_c \psi^b_{\overline{c}}} \\
&= \frac{(\alpha_b - \alpha_{j-1})(\alpha_b - \alpha_{\overline{j-1}})}
  {\psi^b_{j-1} \psi^b_{\overline{j-1}}} \Phi^b_{j-1}, \quad \Phi^b_
  {\overline{l}}=1.
\end{align*}

\noindent Then $P_{\overline{b}b}$ can be written as

\begin{equation*}
P_{\overline{b}b} = (-1)^{[b]} q^{2(\rho, \varepsilon_b)} \prod_{^{c=b+1}_{c 
  \neq 0}}^{\overline{b}} \frac{\psi^b_c}{(\alpha_b^{\Lambda}-\alpha_c^
  {\Lambda})} \Biggl[ \frac{\psi^b_0}{\alpha_b - \alpha_0} - \sum_{j=\overline
  {l}}^{\overline{b}} \frac{(-1)^{[j]} q^{2(\rho,\varepsilon_j)-2C(\Lambda_0)}}
  {\psi^b_j  \psi^b_{\overline{j}}} \Phi^b_j \Biggr].
\end{equation*}

\noindent Note that for $c \neq 0$,

\begin{align*}
\psi^b_c &= \frac{q^{-C(\Lambda_0)}}{(q-q^{-1})} \bigl( q^{(\varepsilon_b,
  2\rho+2\Lambda +\varepsilon_b)} - q^{(\varepsilon_c, 2\rho+2\Lambda +
  \varepsilon_c)} +(q-q^{-1})(-1)^{[c]} q^{(\varepsilon_c, 2\rho + 2\Lambda)} 
  \bigr) \notag \\
&=\frac{q^{-C(\Lambda_0)}}{(q-q^{-1})} \bigl(q^{(\varepsilon_b, 2\rho +
  2\Lambda +\varepsilon_b)} - q^{(\varepsilon_c, 2\rho+2\Lambda-\varepsilon_c)}
  \bigr) \notag \\
&= \frac{ q^{-C(\Lambda_0)}\,\tilde{\psi}^b_c}{(q-q^{-1})},
\end{align*}

\noindent where

\begin{equation*}
\tilde{\psi}^b_c = q^{(\varepsilon_b, 2\rho +2\Lambda + \varepsilon_b)} - q^{(\varepsilon_c, 2\rho+2\Lambda-\varepsilon_c)}.
\end{equation*}

\noindent So

\begin{align}
\sum_{j=\overline{l}}^{\overline{b}} \frac{(-1)^{[j]}q^{2(\rho,\varepsilon_j)
  -2C(\Lambda_0)}} {\psi^b_j  \psi^b_{\overline{j}}} \Phi^b_j
&=(q-q^{-1}) \sum_{j=\overline{l}}^{\overline{b}} \frac{(-1)^{[j]} (q-q^{-1}) 
  q^{2(\rho,\varepsilon_j)}} {\tilde{\psi}^b_j  \tilde{\psi}^b_{\overline{j}}} 
  \Phi^b_j \notag \\
&=(q-q^{-1})\sum_{j=\overline{l}}^{\overline{b}} \frac{(q^{2(\varepsilon_j, 
  \varepsilon_j)}-1) q^{2(\rho,\varepsilon_j) - (\varepsilon_j,
  \varepsilon_j)}}{\tilde{\psi}^b_j  \tilde{\psi}^b_{\overline{j}}} \Phi^b_j 
\label{sumjlb}
\end{align}

\noindent and

\begin{align*}
\Phi^b_{j+1} &= \frac{(\alpha_b - \alpha_j)(\alpha_b-\alpha_{\overline{j}})}
  {\psi^b_j \psi^b_{\overline{j}}} \Phi^b_{j} \\
&= \frac{(q^{(\varepsilon_b, \varepsilon_b + 2\rho + 2\Lambda)} - 
  q^{(\varepsilon_j, \varepsilon_j + 2\rho + 2\Lambda)})(q^{(\varepsilon_b, 
  \varepsilon_b + 2\rho + 2\Lambda)} - q^{(\varepsilon_j, \varepsilon_j - 2\rho
  - 2\Lambda)})}{\tilde{\psi}^b_j  \tilde{\psi}^b_{\overline{j}}} \Phi^b_{j}
\end{align*}

\noindent for $j \geq \overline{l}$. Now 

\begin{align*}
&(q^{(\varepsilon_b, \varepsilon_b + 2\rho + 2\Lambda)} - q^{(\varepsilon_j, 
  \varepsilon_j + 2\rho + 2\Lambda)})(q^{(\varepsilon_b, \varepsilon_b + 2\rho 
  + 2\Lambda)} - q^{(\varepsilon_j, \varepsilon_j - 2\rho - 2\Lambda)})\notag\\
& \hspace{1cm} = q^{2(\varepsilon_j, \varepsilon_j)} (q^{(\varepsilon_b, 
  \varepsilon_b + 2\rho + 2\Lambda)} - q^{(\varepsilon_j, -\varepsilon_j + 
  2\rho + 2\Lambda)}) (q^{(\varepsilon_b, \varepsilon_b + 2\rho + 2\Lambda)} - 
  q^{-(\varepsilon_j, \varepsilon_j + 2\rho + 2\Lambda)}) \notag\\ 
& \hspace{1cm} \qquad + q^{2(\varepsilon_b,\varepsilon_b + 
  2\rho + 2\Lambda)} (1-q^{2(\varepsilon_j, \varepsilon_j)}) + 
  q^{2(\varepsilon_j, \varepsilon_j)}-1  \notag\\
&\hspace{1cm}=q^{2(\varepsilon_j, \varepsilon_j)} \tilde{\psi}^b_j \tilde{\psi}
  ^b_{\overline{j}} - (q^{2(\varepsilon_b,\varepsilon_b + 2\rho + 2\Lambda)}-1)
  (q^{2(\varepsilon_j, \varepsilon_j)}-1).
\end{align*}

\noindent Then, for $j\geq \overline{l}$,

\begin{equation} \label{cancel}
\frac{\Phi^b_{j+1}}{(q^{2(\varepsilon_b,\varepsilon_b + 2\rho +2\Lambda)}-1)} =
\Bigl[ \frac{q^{2(\varepsilon_j, \varepsilon_j)}}{(q^{2(\varepsilon_b,
  \varepsilon_b + 2\rho + 2\Lambda)}-1)} - \frac{(q^{2(\varepsilon_j, 
  \varepsilon_j)}-1)} {\tilde{\psi}^b_j \tilde{\psi}^b_{\overline{j}}} \Bigr] 
  \Phi^b_{j}.
\end{equation}

\noindent Now for $j= \overline{b}$

\begin{align*}
&\frac{(q^{2(\varepsilon_j, \varepsilon_j)}-1) q^{2(\rho,\varepsilon_j) 
  -(\varepsilon_j, \varepsilon_j)}}{\tilde{\psi}^b_j  \tilde{\psi}^b_{\overline
  {j}}}  \\
& \hspace{25mm} = \frac{(q^{2(\varepsilon_b, \varepsilon_b)}-1) 
  q^{2(\rho,\varepsilon_{\overline{b}}) - (\varepsilon_b, \varepsilon_b)}}
  {(q^{(\varepsilon_b, 2\rho+2\Lambda + \varepsilon_b)} - q^{-(\varepsilon_b, 
  2\rho + 2\Lambda + \varepsilon_b)}) q^{(\varepsilon_b, 2\rho + 2\Lambda)} 
  (q^{(\varepsilon_b, \varepsilon_b)} - q^{-(\varepsilon_b, \varepsilon_b)})}\\
& \hspace{25mm} = \frac{q^{2(\rho,\varepsilon_{\overline{b}}) +(\varepsilon_b,
  \varepsilon_b)}} {(q^{2(\varepsilon_b,\varepsilon_b+2\rho + 2\Lambda)} - 1)},
\end{align*}

\noindent which can be written as

\begin{equation*}
  \frac{(q^{2(\varepsilon_j, \varepsilon_j)}-1) q^{2(\rho,\varepsilon_j) - 
  (\varepsilon_j, \varepsilon_j)}}{\tilde{\psi}^b_j  \tilde{\psi}^b_{\overline
  {j}}} = \frac{q^{2(\rho,\varepsilon_{\overline{b}-1})-(\varepsilon_
  {\overline{b}-1}, \varepsilon_{\overline{b}-1})}} {(q^{2(\varepsilon_b,
  \varepsilon_b+ 2\rho +2\Lambda) } - 1)}
\end{equation*}

\noindent when $b<l$.  Hence equation \eqref{cancel} can be used to pairwise cancel the terms in the sum in equation \eqref{sumjlb}.  Adding the first two terms ($j=\overline{b},\, \overline{b}-1$), we find:

\begin{align*}
q^{2(\rho,\varepsilon_{\overline{b}-1})-(\varepsilon_{\overline{b}-1},
  \varepsilon_{\overline{b}-1})} \Bigl[ \frac{\Phi^b_{\overline{b}}} {(q^{2
  (\varepsilon_b,\varepsilon_b + 2\rho + 2\Lambda)} - 1)} &+ \frac{(q^{2
  (\varepsilon_{\overline{b}-1}, \varepsilon_{\overline{b}-1})}-1)} {\tilde
  {\psi}^b_{\overline{b}-1}  \tilde{\psi}^b_{b+1}} \Phi^b_{\overline{b}-1} 
  \Bigr] \notag \\
&= q^{2(\rho,\varepsilon_{\overline{b}-1})-(\varepsilon_{\overline{b}-1},
  \varepsilon_{\overline{b}-1})} \frac{q^{2(\varepsilon_{\overline{b}-1}, 
  \varepsilon_{\overline{b}-1})}}{(q^{2(\varepsilon_b,\varepsilon_b + 2\rho + 
  2\Lambda)}-1)}\Phi^b_{\overline{b}-1}  \notag \\
&= \frac{q^{2(\rho,\varepsilon_{\overline{b}-2})-(\varepsilon_{\overline{b}-2},
  \varepsilon_{\overline{b}-2})}}{(q^{2(\varepsilon_b,\varepsilon_b + 2\rho + 
  2\Lambda)} - 1)}\Phi^b_{\overline{b}-1}.
\end{align*}

\noindent Continuing to apply equation (\ref{cancel}) in this manner gives

\begin{align}
\sum_{j=\overline{l}}^{\overline{b}} \frac{(q^{2(\varepsilon_j, 
  \varepsilon_j)}-1) q^{2(\rho,\varepsilon_j) - (\varepsilon_j,
  \varepsilon_j)}}{\tilde{\psi}^b_j  \tilde{\psi}^b_{\overline{j}}} \Phi^b_j
&= \frac{q^{2(\rho,\varepsilon_{\overline{l}})+(\varepsilon_l,
  \varepsilon_l)}}{(q^{2 (\varepsilon_b,\varepsilon_b + \rho + \Lambda)} - 1)}
  \Phi^b_{\overline{l}}\notag  \\
&= \frac{q^{2l+1-m}}{(q^{2 (\varepsilon_b,\varepsilon_b +2\rho +2\Lambda)}-1)}.
  \label{sumprod}
\end{align}

\noindent Hence in the case $m=2l+1$ 

\begin{equation*}
P_{\overline{b}b} = (-1)^{[b]} q^{2(\rho,\varepsilon_b)} \Bigl[\frac{\psi^b_0}
  {\alpha_b - \alpha_0} - \frac{(q-q^{-1})} {(q^{2 (\varepsilon_b,\varepsilon
  _b + 2\rho + 2\Lambda)}- 1)} \Big]\prod_{^{c=b+1}_{c \neq 0}}^{\overline{b}} 
  \frac{\psi^b_c}{(\alpha_b^{\Lambda}-\alpha_c^{\Lambda})}.
\end{equation*}

\noindent By substituting in the formulae for $\psi^b_c$ and $\alpha_b$ we obtain

\begin{align*}
P_{\overline{b}b}&=(-1)^{[b]} q^{2(\rho, \varepsilon_b)}  \Bigl[ 1 +(q-q^{-1})
  \Bigl( \frac{1} {q^{(\varepsilon_b,\varepsilon_b +2\rho +2\Lambda)}-1)} - 
  \frac{1} {(q^{2 (\varepsilon_b,\varepsilon_b + 2\rho + 2\Lambda)}- 1)} \Bigr)
  \Bigr] \\
& \hspace{110mm} \prod_{^{c=b+1}_{c \neq 0}}^{\overline{b}} \frac{\psi^b_c} 
  {(\alpha_b^{\Lambda}- \alpha_c^{\Lambda})}  \\
&= (-1)^{[b]} q^{2(\rho, \varepsilon_b)}  \Bigl[1 + (q-q^{-1}) \frac{q^{
  (\varepsilon_b, \varepsilon_b + 2\rho +2\Lambda)}}{(q^{2 (\varepsilon_b,
  \varepsilon_b + 2\rho + 2\Lambda)}- 1)} \Bigr]  \\
& \hspace{6cm} \prod_{c=b+1}^{\overline{b}} \frac{
  (q^{(\varepsilon_b, 2\rho +2\Lambda +\varepsilon_b)} - 
  q^{(\varepsilon_c, 2\rho+2\Lambda-\varepsilon_c)})}{(q^{(\varepsilon_b, 2\rho
  +2\Lambda+\varepsilon_b)}-q^{(\varepsilon_c,2\rho+2\Lambda+\varepsilon_c)})},
\end{align*}

\noindent and thus for $a \geq \overline{b} > 0$

\begin{multline*}
P_{ab} =  (-1)^{[b]} q^{2(\rho, \varepsilon_b)}  \Bigl[1 + (q-q^{-1}) \frac
  {q^{(\varepsilon_b, \varepsilon_b + 2\rho +2\Lambda)}}{(q^{2 
  (\varepsilon_b, \varepsilon_b + 2\rho + 2\Lambda)}- 1)} \Bigr] \\
\prod_{c=b+1}^{a} \frac{(q^{(\varepsilon_b, 2\rho +2\Lambda +
  \varepsilon_b)} - q^{(\varepsilon_c, 2\rho+2\Lambda-\varepsilon_c)})}
  {(q^{(\varepsilon_b, 2\rho + 2\Lambda+\varepsilon_b)}-q^{(\varepsilon_c, 
  2\rho+2\Lambda+\varepsilon_c)})}.
\end{multline*}

\noindent Similarly, we find from equations \eqref{abbar0}, \eqref{bbbar}, \eqref{sumjlb} and \eqref{sumprod} that if $m$ is even then 

\begin{equation*}
P_{ab} = (-1)^{[b]} q^{2(\rho,\varepsilon_b)} \Bigl[1 - \frac{q(q-q^{-1})} 
  {(q^{2 (\varepsilon_b,\varepsilon_b + 2\rho + 2\Lambda)}- 1)} \Big] \prod_
  {c=b+1}^{a} \frac{(q^{(\varepsilon_b, 2\rho +2\Lambda +\varepsilon_b)} - 
  q^{(\varepsilon_c, 2\rho+2\Lambda-\varepsilon_c)})}{(q^{(\varepsilon_b, 2\rho
  +2\Lambda+\varepsilon_b)}-q^{(\varepsilon_c,2\rho+2\Lambda+\varepsilon_c)})}
\end{equation*}

\noindent for $a \geq \overline{b} > 0$. Hence we have found expressions for $P_{ab}$ for all $a,b$ satisfying $a \geq \overline{b} >0$.  At the end of the chapter these, together with the earlier results for $P_{ab}$, will be used to calculate $\chi_\Lambda (C_l)$.

\

\noindent Now we return to our diagonalising matrix $N$.  We know

\begin{equation*}
(N^{-1}M)_{ab} = \alpha_a^\Lambda (N^{-1})_{ab}.
\end{equation*}

\noindent Substituting in the values for $M_{ab}$ gives

\begin{multline} \label{Q} 
\alpha_b^\Lambda (N^{-1})_{ab} + (-1)^{[b]} q^{(2\rho, \varepsilon_b) - 
  C(\Lambda_0)} \sum_{c>b} q^{(2\Lambda, \varepsilon_c)} (1- \delta^c_{
  \overline{b}} (-1)^{[b]} q^{-2(\rho, \varepsilon_b)}) (N^{-1})_{ac} \\
= \alpha_a^{\Lambda} (N^{-1})_{ab}.
\end{multline}

\noindent Set

\begin{equation*}
\hat{Q}_{ab} = \sum_{c \geq b} q^{2(\Lambda, \varepsilon_c)} (N^{-1})_{ac}.
\end{equation*}

\noindent We now solve for $\hat{Q}_{ab}$, using a similar method as for $P_{ab}$.  Once we have an expression for $\hat{Q}_{a\, \nu=1}$ the calculation of $\chi_\Lambda(C_l)$ will be straightforward.  Firstly, from equation \eqref{Q} we have

\begin{alignat*}{2}
&&&(\alpha_a^{\Lambda} - \alpha_b^\Lambda) (N^{-1})_{ab} =(-1)^{[b]} q^{2(\rho,
  \varepsilon_b) - C(\Lambda_0)} \hat{Q}_{a\, b+1} - \theta_{b0}q^{-2(\Lambda, 
  \varepsilon_b) - C(\Lambda_0)} (N^{-1})_{a\overline{b}} \\
&\Rightarrow \quad&&(\alpha_a^{\Lambda} - \alpha_b^\Lambda) q^{-2(\Lambda, 
   \varepsilon_b)} (\hat{Q}_{ab} -\hat{Q}_{a\, b+1})  \\
&&&\hspace{33mm}  = (-1)^{[b]} q^{2(\rho, 
  \varepsilon_b) - C(\Lambda_0)} \hat{Q}_{a\, b+1} - \theta_{b0} q^{-2(\Lambda,
  \varepsilon_b) - C(\Lambda_0)} (N^{-1})_{a\overline{b}}.
\end{alignat*} 

\noindent This gives the recursion relation

\begin{align}
&\hat{Q}_{ab} = \frac{(\alpha_a^{\Lambda} - \alpha_b^\Lambda + (-1)^{[b]} 
  q^{2(\rho+\Lambda, \varepsilon_b) - C(\Lambda_0)})}{(\alpha_a^{\Lambda} - 
  \alpha_b^\Lambda)} \hat{Q}_{a\, b+1} - \frac{\theta_{b0} q^{- C(\Lambda_0)}} 
  {(\alpha_a^{\Lambda} - \alpha_b^\Lambda)} (N^{-1})_{a\overline{b}}, \notag \\
\text{i.e.} \qquad &\hat{Q}_{ab} = \frac{\psi^a_b}{(\alpha_a^{\Lambda} - \alpha_b^\Lambda)} 
  \hat{Q}_{a\,b+1} - \frac{\theta_{b0} q^{- C(\Lambda_0)}} 
  {(\alpha_a^{\Lambda} - \alpha_b^\Lambda)} (N^{-1})_{a\overline{b}}.
  \label{Qab2}
\end{align}

\noindent But $ (N^{-1})_{ab} = 0$ for all $a<b$, and $N_{aa}=1 \Rightarrow \hat{Q}_{aa} = q^{2(\Lambda, \varepsilon_a)}$ for all $a$.  Thus for $0\leq b<a$ and $b < a \leq 0$ we have

\begin{equation*}
\hat{Q}_{ab} = q^{2(\Lambda, \varepsilon_a)} \prod_{c=b}^{a-1} 
  \frac{\psi^a_c}{(\alpha_a^{\Lambda} - \alpha_c^\Lambda)}.
\end{equation*}

\noindent Similarly, for $b < \overline{a} < 0$

\begin{equation*}
\hat{Q}_{ab} = \hat{Q}_{a\overline{a}} \prod_{c=b}^{\overline{a}-1} 
  \frac{\psi^a_c}{(\alpha_a^{\Lambda} - \alpha_c^\Lambda)}.
\end{equation*}

\noindent Now consider the remaining case, namely $\overline{a} \leq b < 0$, in which case the last term of equation (\ref{Qab2}) is significant.  Then

\begin{multline*}
\hat{Q}_{ab} = q^{2(\Lambda, \varepsilon_a)} \prod_{c=b}^{a-1} \frac{\psi^
  a_c} {(\alpha_a^{\Lambda} - \alpha_c^\Lambda)} - \frac{q^{-C(\Lambda_0)}}
  {(\alpha_a^{\Lambda} - \alpha_b^\Lambda)} (N^{-1})_{a\overline{b}} \\
- \sum_{j=b+1}^l \frac{q^{-C(\Lambda_0)}}{(\alpha_a^{\Lambda} - \alpha_j^
  \Lambda)} (N^{-1})_{a\overline{j}} \prod_{c=b}^{j-1} \frac{\psi^a_c}
  {(\alpha_a^{\Lambda} - \alpha_c^\Lambda)}.
\end{multline*}

\noindent However for $a>b>0$

\begin{align*}
(N^{-1})_{ab} &= \frac{(-1)^{[b]} q^{2(\rho, \varepsilon_b) - C(\Lambda_0)}}
  {(\alpha_a^\Lambda - \alpha_b^\Lambda)} \hat{Q}_{a\, b+1} \notag \\
&= \frac{(-1)^{[b]} q^{2(\Lambda, \varepsilon_a) + 2(\rho, \varepsilon_b) - 
  C(\Lambda_0)}} {(\alpha_a^\Lambda - \alpha_b^\Lambda)} \prod_{c=b+1}^{a-1} 
  \frac{\psi^a_c}{(\alpha_a^{\Lambda} - \alpha_c^\Lambda)}.
\end{align*}

\noindent Hence we find

\begin{multline*}
\hat{Q}_{a\overline{a}} = q^{2(\Lambda, \varepsilon_a)} \prod_{c=\overline
  {a}}^{a-1} \frac{\psi^a_c} {(\alpha_a^{\Lambda} - \alpha_c^\Lambda)} - 
  \frac{q^{-C(\Lambda_0)}} {(\alpha_a^{\Lambda} - \alpha_{\overline{a}}^
  \Lambda)}\\
- \sum_{j=\overline{a}+1}^l \frac{(-1)^{[j]} q^{2(\Lambda, \varepsilon_a) - 
  2(\rho, \varepsilon_j) - 2C(\Lambda_0)}} {(\alpha_a^{\Lambda} - \alpha_j^
  \Lambda) (\alpha_a^\Lambda -\alpha_{\overline{j}}^\Lambda)} \prod_{c=
  \overline{j}+1}^{a-1} \frac{\psi^a_c}{(\alpha_a^{\Lambda} -\alpha_c^\Lambda)}
  \prod_{c={\overline{a}}}^{j-1} \frac{\psi^a_c} {(\alpha_a^{\Lambda} - 
  \alpha_c^\Lambda)},
\end{multline*}

\noindent which can also be written as

\begin{align*}
\hat{Q}_{a\overline{a}}&= q^{2(\Lambda, \varepsilon_a)} \prod_{c=\overline
  {a}}^{a-1}\frac{\psi^a_c} {(\alpha_a^{\Lambda} - \alpha_c^\Lambda)} \Biggl[ 1
  - \frac{q^{-2(\Lambda, \varepsilon_a)-C(\Lambda_0)}} {(\alpha_a^{\Lambda} - 
  \alpha_{\overline{a}}^\Lambda)} \prod_{c=\overline{a}}^{a-1} \frac{(\alpha_a
  ^{\Lambda} - \alpha_c^\Lambda)}{\psi^a_c} \notag \\
& \hspace{5cm} - \sum_{j=\overline{a}+1}^l \frac{(-1)^{[j]} q^{- 2(\rho, 
  \varepsilon_j) - 2C(\Lambda_0)}}{(\alpha_a^{\Lambda} - \alpha_j^\Lambda) 
  (\alpha_a^\Lambda -\alpha_{\overline{j}}^\Lambda)} \prod_{c=j}^{\overline{j}}
  \frac{(\alpha_a^{\Lambda} - \alpha_c^\Lambda)} {\psi^a_c} \Biggr] \notag \\
&= q^{2(\Lambda, \varepsilon_a)} \prod_{c=\overline{a}}^{a-1} \frac{\psi^a_c} 
  {(\alpha_a^{\Lambda} - \alpha_c^\Lambda)} \Biggl[ 1-\frac{q^{- 2(\Lambda, 
  \varepsilon_a)-C(\Lambda_0)}}{\psi^a_{\overline{a}}} \prod_{c=\overline{a}+1}
  ^{a-1} \frac{(\alpha_a^{\Lambda} - \alpha_c^\Lambda)}{\psi^a_c} \notag \\
& \hspace{5cm} - \sum_{j=\overline{a}+1}^l \frac{(-1)^{[j]} q^{- 2(\rho, 
  \varepsilon_j) - 2C(\Lambda_0)}}{\psi^a_j \psi^a_{\overline{j}}} 
  \prod_{c=j+1}^{\overline{j}-1} \frac{(\alpha_a^{\Lambda} - \alpha_c^\Lambda)}
  {\psi^a_c} \Biggr] \notag \\
&= q^{2(\Lambda, \varepsilon_a)} \prod_{c=\overline{a}}^{a-1} \frac{\psi^a_c} 
  {(\alpha_a^{\Lambda} - \alpha_c^\Lambda)} \Biggl[ 1-\sum_{j=\overline{a}}^l
  \frac{(-1)^{[j]} q^{- 2(\rho, \varepsilon_j) - 2C(\Lambda_0)}}{\psi^a_j 
  \psi^a_{\overline{j}}}\prod_{c=j+1}^{\overline{j}-1}\frac{(\alpha_a^{\Lambda}
  - \alpha_c^\Lambda)} {\psi^a_c} \Biggr].
\end{align*}

\noindent  For odd $m$ this becomes 

\begin{align*}
\hat{Q}_{a\overline{a}} = q^{2(\Lambda,\varepsilon_a)}\prod_{_{c\neq 0}^
  {c=\overline{a}}}^{a-1} \frac{\psi^a_c} {(\alpha_a^{\Lambda} - \alpha_c^
  \Lambda)} \Biggl[ \frac{\psi^a_0} {(\alpha_a^{\Lambda} -\alpha_0^\Lambda)}&- 
  \sum_{j=\overline{l}}^a \frac{(-1)^{[j]} q^{2(\rho, \varepsilon_j)- 
  2C(\Lambda_0)}}{\psi^a_j \psi^a_{\overline{j}}} \Phi^a_{j} \Biggr] \notag \\
= q^{2(\Lambda,\varepsilon_a)}\prod_{_{c\neq 0}^{c=\overline
  {a}}}^{a-1} \frac{\psi^a_c} {(\alpha_a^{\Lambda} - \alpha_c^\Lambda)}\Biggl[ 
  \frac{\psi^a_0} {(\alpha_a^{\Lambda} - \alpha_0^\Lambda)}& \notag \\
- (q-&q^{-1}) \sum_
  {j=\overline{l}}^a \frac{(q^{2(\varepsilon_j, \varepsilon_j)}-1) q^{2(\rho,
  \varepsilon_j) - (\varepsilon_j,\varepsilon_j)}}{\tilde{\psi}^a_j  \tilde
  {\psi}^a_{\overline{j}}} \Phi^a_j \Biggr]. 
\end{align*}

\noindent Recall that when $j \geq \overline{l}$

\begin{equation*} 
\frac{\Phi^a_{j+1}}{(q^{2(\varepsilon_a,\varepsilon_a + 2\rho +2\Lambda)}-1)} =
\Bigl[ \frac{q^{2(\varepsilon_j, \varepsilon_j)}}{(q^{2(\varepsilon_a,
  \varepsilon_a + 2\rho + 2\Lambda)}-1)} - \frac{(q^{2(\varepsilon_j, 
  \varepsilon_j)}-1)} {\tilde{\psi}^a_j  \tilde{\psi}^a_{\overline{j}}} \Bigr] 
  \Phi^a_{j}.
\end{equation*}

\noindent Also, for $j=a > \overline{l}$ we have

\begin{align*}
&\frac{(q^{2(\varepsilon_j, \varepsilon_j)}-1) q^{2(\rho, \varepsilon_j) - 
  (\varepsilon_j,\varepsilon_j)}}{\tilde{\psi}^a_j  \tilde {\psi}^a_
  {\overline{j}}} \notag \\
& \hspace{25mm}= \frac{(q^{2(\varepsilon_a, \varepsilon_a)}-1) q^{2(\rho, 
  \varepsilon_a) - (\varepsilon_a,\varepsilon_a)}}{(q^{(\varepsilon_a, 
  \varepsilon_a + 2\rho + 2\Lambda)} - q^{-(\varepsilon_a,\varepsilon_a +2\rho 
  + 2\Lambda)}) q^{(\varepsilon_a, 2\rho+2\Lambda)}(q^{(\varepsilon_a, 
  \varepsilon_a)} - q^{-(\varepsilon_a, \varepsilon_a)})} \notag \\
& \hspace{25mm} = \frac{q^{2(\rho, \varepsilon_{a-1}) -(\varepsilon_{a-1},
  \varepsilon_{a-1})}} {(q^{2(\varepsilon_a, \varepsilon_a + 2\rho + 2\Lambda)}
  -1)}.
\end{align*}

\noindent Then, cancelling terms pairwise in the same manner as for $P_{\overline{b}b}$, we find:

\begin{equation*}
\sum_{j=\overline{l}}^a \frac{(q^{2(\varepsilon_j, \varepsilon_j)}-1) 
  q^{2(\rho, \varepsilon_j) - (\varepsilon_j,\varepsilon_j)}}{\tilde{\psi}^a_j 
  \tilde {\psi}^a_{\overline{j}}} \Phi^a_j  = \frac{q^{2l+1-m}}
  {(q^{2(\varepsilon_a, \varepsilon_a + 2\rho + 2\Lambda)} -1)}.
\end{equation*}

\noindent Hence for $m=2l+1$ we have

\begin{equation*}
\hat{Q}_{a\overline{a}} = q^{2(\Lambda,\varepsilon_a)} \Bigl[ 1+(q-q^{-1})
  \frac{q^{(\varepsilon_a, \varepsilon_a + 2\rho +2\Lambda)}}{(q^{2 
  (\varepsilon_a, \varepsilon_a + 2\rho + 2\Lambda)}- 1)} \Bigr] \prod_
  {c=\overline{a}}^{a-1} \frac{(q^{(\varepsilon_a, 2\rho +2\Lambda +
  \varepsilon_a)} - q^{(\varepsilon_c, 2\rho+2\Lambda-\varepsilon_c)})}
  {(q^{(\varepsilon_a, 2\rho +2\Lambda+\varepsilon_a)}-q^{(\varepsilon_c, 2\rho
  +2\Lambda+\varepsilon_c)})},
\end{equation*}

\noindent and for $b \leq \overline{a} < 0$

\begin{equation*}
\hat{Q}_{ab} = q^{2(\Lambda,\varepsilon_a)} \Bigl[ 1 + (q-q^{-1}) \frac
  {q^{(\varepsilon_a, \varepsilon_a + 2\rho +2\Lambda)}}{(q^{2 
  (\varepsilon_a, \varepsilon_a + 2\rho + 2\Lambda)}- 1)}
   \Bigr] 
\prod_{c=b}^{a-1} \frac{(q^{(\varepsilon_a, 2\rho +2\Lambda +
  \varepsilon_a)} - q^{(\varepsilon_c, 2\rho+2\Lambda-\varepsilon_c)})}
  {(q^{(\varepsilon_a, 2\rho +2\Lambda+\varepsilon_a)}-q^{(\varepsilon_c, 2\rho
  +2\Lambda+\varepsilon_c)})}.
\end{equation*}

\noindent Similarly, for even $m$ we find

\begin{equation*}
\hat{Q}_{ab} = q^{2(\Lambda,\varepsilon_a)} \Bigl[ 1 - \frac{q(q-q^{-1})} 
  {(q^{2 (\varepsilon_a,\varepsilon_a + 2\rho + 2\Lambda)}- 1)} \Big]
\prod_{c=b}^{a-1} \frac{(q^{(\varepsilon_a, 2\rho +2\Lambda +
  \varepsilon_a)} - q^{(\varepsilon_c, 2\rho+2\Lambda-\varepsilon_c)})}
  {(q^{(\varepsilon_a, 2\rho +2\Lambda+\varepsilon_a)}-q^{(\varepsilon_c, 2\rho
  +2\Lambda+\varepsilon_c)})}.
\end{equation*}

\noindent for $b \leq \overline{a} < 0$. This completes the difficult calculations. 

\

\noindent To use these results to calculate $\chi_\Lambda (C_l)$ we introduce a new function $Q_{ab}$, defined by:

\begin{equation*}
Q_{ab} = \sum_{c\geq b} (N^{-1})_{ac}.
\end{equation*}

\noindent Then from equations \eqref{chi} and \eqref{P} we deduce 

\begin{equation} \label{PQ}
\chi_{\Lambda}(C_l) = \sum_a {(\alpha_a^{\Lambda})}^l P_{\mu=\overline{1}\, a} Q_{a\, \nu = 1}.
\end{equation}

\noindent However we know

\begin{equation*}
t_a^{(l)} = \frac{q^{2(\Lambda, \varepsilon_a)} - 1}{q-q^{-1}},
\end{equation*}

\noindent and 

\begin{align*}
\sum_b (N^{-1})_{ab} t_b^{(1)} &= \sum_{b,c} (N^{-1})_{ab} M_{bc} 
  t_c^{(0)} \notag \\
\Rightarrow \quad \sum_b (N^{-1})_{ab}\frac{(q^{2(\Lambda,\varepsilon_b)} - 1)}
  {(q-q^{-1})} &= \sum_b (N^{-1}M)_{ab}  \notag \\
&= \sum_b \alpha_a^\Lambda (N^{-1})_{ab}  \notag \\
&= \sum_b (N^{-1})_{ab} \frac{(q^{(\varepsilon_a,2\Lambda+2\rho +\varepsilon_a)
  - C(\Lambda_0)}-1)}{(q-q^{-1})}.
\end{align*}

\noindent Thus, remembering that 

\begin{equation*}
\hat{Q}_{ab} = \sum_{c\geq b} q^{2(\Lambda, \varepsilon_c)} (N^{-1})_{ac},
\end{equation*}

\noindent we have

\begin{equation*}
Q_{a\, \nu=1} = q^{C(\Lambda_0) - (\varepsilon_a, 2\rho + 2\Lambda +
  \varepsilon_a)} \hat{Q}_{a\, \nu=1}.
\end{equation*}

\noindent  Substituting our formulae for $ P_{\mu=\overline{1}\, a}$ and $Q_{a\, \nu=1}$ into equation \eqref{PQ}, noting that for $a \neq 0$ exactly one of $a<0$ and $a>0$ is true, we find the eigenvalues of the Casimir invariants $C_l$ are given by:

\begin{multline*} 
\chi_{\Lambda}(C_l) =  \sum_a (-1)^{[a]} q^{C(\Lambda_0) - (\varepsilon_a, 
  \varepsilon_a)} f(a) \Bigl[ \frac{(q^{(\varepsilon_a, 2 \rho + 2 \Lambda +
  \varepsilon_a)-C(\Lambda_0)}-1)}{(q-q^{-1})} \Bigr]^l  \\ \times
  \prod_{b \neq a} \frac{(q^{(\varepsilon_a, 2\rho +2\Lambda +
  \varepsilon_a)} - q^{(\varepsilon_b, 2\rho+2\Lambda-\varepsilon_b)})}
  {(q^{(\varepsilon_a, 2\rho +2\Lambda+\varepsilon_a)}-q^{(\varepsilon_b, 2\rho
  +2\Lambda+\varepsilon_b)})},
\end{multline*}

\noindent where 

\begin{equation*} 
f(a) = \begin{cases}
1 - (q-q^{-1})\frac{q} {(q^{2 (\varepsilon_a,\varepsilon_a + 2\rho + 2\Lambda)}
  - 1)}, & m=2l, \\
1 + (q-q^{-1}) \frac{q^{(\varepsilon_a,\varepsilon_a + 2\rho +
  2\Lambda)}} {(q^{2(\varepsilon_a, \varepsilon_a + 2\rho + 2\Lambda)}- 1)}, 
  \quad& a \neq 0, \; m=2l+1, \\
1, & a=0, \; m=2l+1.
\end{cases}
\end{equation*}

\noindent  Throughout we assumed the eigenvalues were distinct.  If they are not, the calculations are more complicated but the result is the same.  Thus, summarising the results from this chapter, we have found:

\begin{theorem}

$U_q[osp(m|n)]$, $m > 2,$ has an infinite family of Casimir invariants of the form

\begin{equation*}
C_l = (str \otimes I)(\pi(q^{2h_p}) \otimes I) A^l, \qquad l \in \mathbb{Z}^+,
\end{equation*}

\noindent where 

\begin{equation*}
A = \frac{(R^T R - I \otimes I)}{(q-q^{-1})}.
\end{equation*}

\noindent The eigenvalues of the invariants when acting on an arbitrary irreducible finite-dimensional module with highest weight $\Lambda$ are given by:

\begin{multline*}
\chi_{\Lambda}(C_l) =  \sum_a (-1)^{[a]} q^{C(\Lambda_0) - (\varepsilon_a, 
  \varepsilon_a)} f(a) \Bigl[ \frac{(q^{(\varepsilon_a, \varepsilon_a +2\rho 
  + 2\Lambda)-C(\Lambda_0)}-1)}{(q-q^{-1})} \Bigr]^l  \\ \times
  \prod_{b \neq a} \frac{(q^{(\varepsilon_a, 2\rho +2\Lambda +
  \varepsilon_a)} - q^{(\varepsilon_b, 2\rho+2\Lambda-\varepsilon_b)})}
  {(q^{(\varepsilon_a, 2\rho +2\Lambda+\varepsilon_a)}-q^{(\varepsilon_b, 2\rho
  +2\Lambda+\varepsilon_b)})},
\end{multline*}

\noindent where 

\begin{equation*}
f(a) = \begin{cases}
1 - (q-q^{-1})\frac{q} {(q^{2 (\varepsilon_a,\varepsilon_a + 2\rho + 2\Lambda)}
  - 1)}, & m=2l, \\
1 + (q-q^{-1}) \frac{q^{(\varepsilon_a,\varepsilon_a + 2\rho +
  2\Lambda)}} {(q^{2(\varepsilon_a, \varepsilon_a + 2\rho + 2\Lambda)}- 1)}, 
  \quad& a \neq 0, \; m=2l+1, \\
1, & a = 0,\; m=2l+1.
\end{cases}
\end{equation*}

\end{theorem}

\

\noindent This completes the calculation of the eigenvalues of an infinite family of Casimir invariants of $U_q[osp(m|n)]$ when acting on an arbitrary irreducible highest weight module, provided $m > 2$ and $n=2k \geq 2$. This had already been done for $U_q[osp(2|n)]$ using a different method in \cite{GLZ}. Moreover, as mentioned earlier, every finite dimensional representation of $U_q[osp(1|n)]$ is isomorphic to a finite dimensional representation of $U_{-q}[so(n+1)]$ \cite{Zhang}, whose central elements are well-understood.  Hence the eigenvalues of a family of Casimir invariants when acting on an arbitrary irreducible highest weight module have now been calculated for all non-exceptional quantum superalgebras.

\chapter{Conclusion}

\noindent  One of the major aims of this thesis was to construct a Lax operator for the $B$ and $D$ type superalgebras.  As this provides a solution to the quantum Yang--Baxter equation in an arbitrary representation, this operator is potentially of great use in integrable systems. In Chapters \ref{R-mat} and \ref{close} we found formulae for the fundamental values and developed a set of inductive and commutative relations that could be used to calculate the remaining matrix entries.  A specific example was given in Chapter \ref{vector}, where the $R$-matrix for the vector representation was calculated from the Lax operator, using a method that can be extended to any other finite-dimensional representation.  The only non-exceptional quantum superalgebras for which no Lax operator is known are now the $C$ series, $osp(2|2n)$.  Although they have a different root system, and thus the solution in this thesis may not be valid for them, it should not be difficult to adjust the method developed here to cover that case.

\

\noindent  Another longstanding problem has been to find families of Casimir invariants for quantum superalgebras and to calculate their eigenvalues when acting on a highest weight module. These will be an important tool in understanding the representation theory associated with the integrable models.  In Chapter \ref{casimir} the Lax operator developed earlier in the thesis was used to do exactly that for all the $B$ and $D$ type quantum superalgebras.  As this had already been done for the $A$ and $C$ type, the solution is now complete for all non-exceptional quantum superalgebras.

\newpage
\addcontentsline{toc}{chapter}{\protect\numberline{}{Bibliography}}
\bibliography{thesis}

\begin{thebibliography}{10}

\bibitem{Arnaudon}
D.~Arnaudon, J.~Avan, N.~Crampe, A.~Doikou, L.~Frappat, and E.~Ragoucy.
\newblock Bethe ansatz equations and exact ${S}$ matrices for the $osp({M}|2n)$
  open super-spin chain.
\newblock {\em Nucl. Phys. B}, {\bf 687}:257--278, 2004.

\bibitem{Baxter72}
R.~J. Baxter.
\newblock Partition function of the eight-vertex lattice model.
\newblock {\em Ann. of Phys.}, {\bf 70}:193--228, (1972).

\bibitem{Baxter}
R.~J. Baxter.
\newblock {\em Exactly solved models in statistical mechanics}.
\newblock Academic Press, London, 1982.

\bibitem{Bincer}
A.~M. Bincer.
\newblock Eigenvalues of {C}asimir operators for the general linear and
  orthosymplectic {L}ie superalgebras.
\newblock {\em J. Math. Phys.}, {\bf 24}:2546--2549, (1983).

\bibitem{BGZ}
A.~J. Bracken, M.~D. Gould, and R.~B. Zhang.
\newblock Quantum supergroups and solutions of the {Y}ang--{B}axter equation.
\newblock {\em Mod. Phys. Lett. A}, {\bf 5}:831--840, (1990).

\bibitem{Chaichian}
M.~Chaichian and P.~Kulish.
\newblock Quantum {L}ie superalgebras and $q$-oscillators.
\newblock {\em Phys. Lett. B}, {\bf 234}:72--80, (1990).

\bibitem{Corwin}
L.~Corwin, Y.~Ne'eman, and S.~Sternberg.
\newblock Graded {L}ie algebras in mathematics and physics.
\newblock {\em Reviews of Modern Physics}, {\bf 47}:573 -- 603, (1975).

\bibitem{Deguchi}
T.~Deguchi, A.~Fujii, and K.~Ito.
\newblock Quantum superalgebra ${U}_q{\rm osp}(2,2)$.
\newblock {\em Phys. Lett. B}, {\bf 238}:242--246, 1990.

\bibitem{Delius}
G.~W. Delius, M.~D. Gould, J.~R. Links, and Y.-Z. Zhang.
\newblock On type 1 quantum affine superalgebras.
\newblock {\em Internat. J. Modern Phys. A}, {\bf 10}:3259--3281, (1995).

\bibitem{Dixmier}
J.~Dixmier.
\newblock {\em Enveloping Algebras}.
\newblock North-Holland, Amsterdam, 1977.

\bibitem{Drinfeld}
V.~G. Drinfeld.
\newblock Quantum groups.
\newblock In {\em Proceedings of the International Congress of Mathematicians},
  volume 1,2, pages 798--820. Amer. Math. Soc., (1987).

\bibitem{Foerster}
A.~Foerster and M.~Karowski.
\newblock The supersymmetric $t$-${J}$ model with quantum group invariance.
\newblock {\em Nucl. Phys. B}, {\bf 408}:512--534, (1993).

\bibitem{Frolicher}
A.~Frolicher and Z.~Nijenhuis.
\newblock A theorem on stability of complex structures.
\newblock {\em Proc. Natl. Acad. Sci.}, {\bf 43}:239--241, (1957).

\bibitem{Martins}
W.~Galleas and M.~J. Martins.
\newblock ${R}$-matrices and spectrum of vertex models based on superalgebras.
\newblock arXiv:nlin.SI/04060003.

\bibitem{Ge}
X.~Y. Ge.
\newblock Integrable boundary conditions for the $q$-deformed extended
  {H}ubbard model.
\newblock {\em Mod. Phys. Lett. B}, {\bf 13}:499--507, (1999).

\bibitem{Gerstenhaber}
M.~Gerstenhaber.
\newblock The cohomology structure of an associative ring.
\newblock {\em Ann. Math.}, {\bf 78}:267--288, (1963).

\bibitem{Gonzalez}
A.~Gonzalez-Ruiz.
\newblock Integrable open-boundary conditions for the supersymmetric $t$-${J}$
  model the quantum-group invariant case.
\newblock {\em Nucl. Phys. B}, {\bf 424}:468--486, (1994).

\bibitem{GHLZ}
M.~D. Gould, K.~E. Hibberd, J.~R. Links, and Y-Z. Zhang.
\newblock Integrable electron model with correlated hopping and quantum
  supersymmetry.
\newblock {\em Phys. Lett. A}, 212:156--160, (1996).

\bibitem{GLZ}
M.~D. Gould, J.~R. Links, and Y.-Z. Zhang.
\newblock Eigenvalues of {C}asimir invariants for type 1 quantum superalgebras.
\newblock {\em Lett. Math. Phys.}, {\bf 36}:415--425, (1996).

\bibitem{GLZT}
M.~D. Gould, J.~R. Links, Y.-Z. Zhang, and I.~Tsohantjis.
\newblock Twisted quantum affine superalgebra ${U}\sb q[{\rm sl}(2\vert 2)\sp
  {(2)}],{U}\sb q[{\rm osp}(2\vert 2)]$ invariant ${R}$-matrices and a new
  integrable electronic model.
\newblock {\em J. Phys. A}, {\bf 31}:4313--4325, (1997).

\bibitem{GS}
M.~D. Gould and N.~I. Stoilova.
\newblock Eigenvalues of {C}asimir operators for ${\rm gl}(m/\infty)$.
\newblock {\em J. Phys. A}, {\bf 32}:391--399, (1999).

\bibitem{GouldZhang99}
M.~D. Gould and Y.-Z. Zhang.
\newblock Quasispin graded-fermion formalism and ${\rm gl}(m\vert
  n)\downarrow{\rm osp}(m\vert n)$ branching rules.
\newblock {\em J. Math. Phys.}, {\bf 40}:5371--5386, (1999).

\bibitem{GouldZhang00}
M.~D. Gould and Y.-Z. Zhang.
\newblock Twisted quantum affine superalgebra {$U_q [gl(m|n)^{(2)}]$} and new
  {$U_q[osp(m|n)]$} invariant {$R$}-matrices.
\newblock {\em Nucl. Phys. B}, {\bf 566}:529--546, (2000).

\bibitem{Jarvis}
P.D. Jarvis and H.S. Green.
\newblock Casimir invariants and characteristic identities for generators of
  the general linear, special linear and orthosymplectic graded {L}ie algebras.
\newblock {\em J. Math. Phys.}, {\bf 20}:2115--2122, (1979).

\bibitem{Jimbo}
M.~Jimbo.
\newblock A $q$-difference analogue of {$U(g)$} and the {Y}ang--{B}axter
  equation.
\newblock {\em Lett. Math. Phys.}, {\bf 10}:63--69, (1985).

\bibitem{Jimbo89}
M.~Jimbo.
\newblock Introduction to the {Y}ang--{B}axter equation.
\newblock {\em Int. J. Mod. Phys. A}, {\bf 4}:3759--3777, (1989).

\bibitem{Kac}
V.~Kac.
\newblock {L}ie superalgebras.
\newblock {\em Advances in Math.}, {\bf 26}:8--96, (1977).

\bibitem{Khoroshkin}
S.~M. Khoroshkin and V.~N. Tolstoy.
\newblock {U}niversal ${R}$-matrix for quantized (super)algebras.
\newblock {\em Comm. Math. Phys.}, {\bf 141}:599--617, (1991).

\bibitem{Kulish}
P.~P. Kulish and N.~Yu Reshetikhin.
\newblock Universal {$R$}-matrix of the quantum superalgebra ${\rm osp}(2\vert
  1)$.
\newblock {\em Lett. Math. Phys.}, {\bf 18}:143--149, (1989).

\bibitem{KS}
P.~P. Kulish and E.~K. Sklyanin.
\newblock Quantum spectral transform method. {R}ecent developments.
\newblock {\em {\rm In} Integrable Quantum Field Theories. Lecture Notes in
  Physics}, {\bf 151}:61--119, (1982).
\newblock Springer-Verlag.

\bibitem{LG}
J.~R. Links and M.~D. Gould.
\newblock Classification of unitary and grade star irreps for
  ${U}_q(osp(2|2n))$.
\newblock {\em J. Math. Phys.}, {\bf 36}:531--545, (1995).

\bibitem{LGZ}
J.~R. Links, M.~D. Gould, and R.~B. Zhang.
\newblock Quantum supergroups, link polynomials and representations of the
  braid generator.
\newblock {\em Rev. Math. Phys.}, {\bf 5}:345--361, (1993).

\bibitem{LinksZhang}
J.~R. Links and R.~B. Zhang.
\newblock Eigenvalues of {C}asimir invariants of {$U_q(gl(m|n))$}.
\newblock {\em J. Math. Phys.}, {\bf 34}:6016--6024, (1993).

\bibitem{MR}
M.~J. Martins and P.~B. Ramos.
\newblock Solution of a supersymmetric model of correlated electrons.
\newblock {\em Phys. Rev. B}, {\bf 56}:6376--6379, (1997).

\bibitem{McGuire}
J.~B. McGuire.
\newblock Study of exactly solvable one-dimensional $n$-body problems.
\newblock {\em J. Math. Phys.}, {\bf 5}:622--636, (1964).

\bibitem{Mehta}
M.~Mehta.
\newblock {\em New Solutions of the Yang--Baxter Equation Associated with
  Quantised Orthosymplectic Lie Superalgebras}.
\newblock PhD thesis, The University of Queensland, 2003.

\bibitem{Nijenhuis}
A.~Nijehuis.
\newblock {J}acobi-type identities for bilinear differential concomitants of
  certain tensor fields. i, ii.
\newblock {\em Nederl. Akad. Wetensch. Proc. Ser. A.}, {\bf 58}:390--403,
  (1955).

\bibitem{Nwachuku}
C.~O. Nwachuku and M.~A. Rashid.
\newblock Eigenvalues of the {C}asimir operators of the orthogonal and
  symplectic groups.
\newblock {\em J. Math. Phys.}, {\bf 17}:1611--1616, (1976).

\bibitem{WebLie}
J.~J. O'Connor and E.~F. Robertson.
\newblock {M}arius {S}ophus {L}ie.
\newblock
  {http://www-gap.dcs.st-and.ac.uk/$\sim$history/Mathematicians/Lie.html},
  February 2000.
\newblock Accessed 23 April 2004.

\bibitem{Popov}
A.~M. Perelomov and V.~S. Popov.
\newblock Casimir operators for the orthogonal and symplectic groups.
\newblock {\em Sov. J. Nucl. Phys.}, {\bf 3}:819--824, (1966).

\bibitem{Perelomov}
A.~M. Perelomov and V.~S. Popov.
\newblock Casimir operators for ${\rm u}(n)$ and ${\rm su}(n)$.
\newblock {\em Sov. J. Nucl. Phys.}, {\bf 3}:676--680, (1966).

\bibitem{Scheun}
M.~Scheunert.
\newblock The ${R}$-matrix of the symplecto-orthogonal quantum superalgebra
  ${U}_q(spo(2n|2m))$ in the vector representation.
\newblock arXiv:math.QA/0004032.

\bibitem{Scheunert83}
M.~Scheunert.
\newblock Eigenvalues of {C}asimir operators for the general linear, the
  special linear, and the orthosymplectic {L}ie superalgebras.
\newblock {\em J. Math. Phys.}, {\bf 24}:2681--2688, (1983).

\bibitem{Serre}
J.-P. Serre.
\newblock {\em Alg{\'e}bres de {L}ie semi-simples complexes}.
\newblock W.A. Benjamin, New York, 1966.

\bibitem{Wadati}
M.~Wadati, T.~Deguchi, and Y.~Akutsu.
\newblock Exactly solvable models and knot theory.
\newblock {\em Phys. Rep.}, {\bf 180}:247--332, (1989).

\bibitem{Weyl}
H.~Weyl.
\newblock {\em The Theory of Groups and Quantum Mechanics}.
\newblock Methuen and Company, Ltd., London, 1931.

\bibitem{Witten}
E.~Witten.
\newblock Quantum field theory and the {J}ones polynomial.
\newblock {\em Commun. Math. Phys.}, {\bf 121}:351--399, (1989).

\bibitem{Yamane}
H.~Yamane.
\newblock On defining relations of the affine {L}ie superalgebras and their
  quantized universal enveloping superalgebras.
\newblock {\em Publ. Res. Inst. Math. Sci.}, {\bf 35}:321--390, (1999).

\bibitem{Yang}
C.~N. Yang.
\newblock Some exact results for the many-body problem in one dimension with
  repulsive delta-function interaction.
\newblock {\em Phys. Rev. Lett.}, {\bf 19}:1312--1315, (1967).

\bibitem{Zhang}
R.~B. Zhang.
\newblock Finite dimensional representations of ${U}_q(osp(1/2n))$ and its
  connection with quantum $so(2n+1)$.
\newblock {\em Lett. Math. Phys.}, {\bf 25}:317--325, (1992).

\bibitem{Zhang2}
R.~B. Zhang.
\newblock Universal ${L}$ operator and invariants of the quantum supergroup
  ${U}\sb q({\rm gl}(m/n))$.
\newblock {\em J. Math. Phys.}, {\bf 33}:1970--1979, (1992).

\bibitem{Zhang95}
R.~B. Zhang.
\newblock Quantum supergroups and topological invariants of three-manifolds.
\newblock {\em Rev. Math. Phys.}, {\bf 7}:809--831, (1995).

\bibitem{ZhangGould}
R.~B. Zhang and M.~J. Gould.
\newblock Universal ${R}$-matrices and invariants of quantum supergroups.
\newblock {\em J. Math. Phys.}, {\bf 32}:3261--3267, (1991).

\end{thebibliography}

\appendix

\chapter{Derivation of the relations used to find the Lax operator}
        \noindent Recall equation \eqref{**}, which states: 

\begin{multline} \label{eqR}
q^{-\frac{1}{2} (\alpha_c, \alpha_c - \varepsilon_{a})} \langle a|e_c|a' 
   \rangle \hat{\sigma}_{ba'}-(-1)^{([a]+[b])[c]} q^{\frac{1}{2} 
   (\alpha_c, \varepsilon_b)} \langle b'|e_c|b \rangle \hat{\sigma}_{b'a} \\
= q^{ (\alpha_c,\varepsilon_b)} \hat{\sigma}_{ba} e_c q^{\frac{1}{2} h_c} - 
   (-1)^ {([a]+[b])[c]} q^{-(\alpha_c,\varepsilon_a)} e_c q^{\frac{1}{2} h_c} 
   \hat{\sigma}_{ba}, \quad \varepsilon_b > \varepsilon_a.
\end{multline}

\noindent By examining this many different relations can be obtained, which were summarised as relations \eqref{qcom} and \eqref{indrel}.  This appendix includes the full list and their derivation.


\section{Relations for $\alpha_i = \varepsilon_i - \varepsilon_{i+1}$, $1 \leq i < l$}

\noindent  As shown in Section \ref{alphai}, in the case $\alpha_i = \varepsilon_i - \varepsilon_{i+1}$ equation (\ref{eqR}) implies

\begin{align} \label{Ai}
\delta_{ai} \hat{\sigma}_{b\, i+1}- \delta_{a\, 
  \overline{i+1}} \hat{\sigma}_{b\, \overline{i}} - \delta_{b\, i+1} 
  \hat{\sigma}_{ia} + \delta_{b\,\overline{i}} \hat{\sigma}_{\overline{i+1}\,a}
&=q^{(\alpha_i,\varepsilon_b)} \hat{\sigma}_{ba} \hat{\sigma}_{i\, i+1} - 
  q^{-(\alpha_i,\varepsilon_a)} \hat{\sigma}_{i\, i+1} \hat{\sigma}_{ba} \\
&=q^{-(\alpha_i,\varepsilon_a)} \hat{\sigma}_{\overline{i+1}\, \overline{i}}
  \hat{\sigma}_{ba} - q^{(\alpha_i,\varepsilon_b)} \hat{\sigma}_{ba} 
  \hat{\sigma}_{\overline{i+1}\, \overline{i}}. \label{Aibar}
\end{align}

\noindent To absorb all the information these equations hold, each case must be considered separately.  Throughout recall that $\hat{\sigma}_{i\, i+1} = - \hat{\sigma}_{\overline{i+1}\, \overline{i}} = q^{\frac{1}{2}} e_i q^{\frac{1}{2} h_i}, \; 1 \leq i < l$.

\

\noindent \underline{Case 1: $a = i$}

\noindent As $\varepsilon_b > \varepsilon_a$, we know that $\delta_{b \, i+1} = \delta_{b\, \overline{i}} = 0$.  Hence equation (\ref{Ai}) reduces to

\begin{align*}
\hat{\sigma}_{b\, i+1} &= q^{(\alpha_i, \varepsilon_b)} \hat{\sigma}_{bi} 
  \hat{\sigma}_{i\, i+1} - q^{-1} \hat{\sigma}_{i\, i+1} \hat{\sigma}_{bi} 
  \notag \\
&= \hat{\sigma}_{bi} \hat{\sigma}_{i\, i+1} - q^{-1} \hat{\sigma}_{i\, i+1} 
  \hat{\sigma}_{bi}, \hspace{2cm} \varepsilon_b > \varepsilon_i.
\end{align*}

\noindent \underline{Case 2: $b = \overline{i}$}

\noindent Again, the constraint $\varepsilon_b > \varepsilon_a$ ensures there is only one non-zero term on the left-hand side, giving

\begin{align*}
\hat{\sigma}_{\overline{i+1}\, a} &= q^{-(\alpha_i, \varepsilon_a)} 
  \hat{\sigma}_{\overline{i+1}\, \overline{i}} \hat{\sigma}_{\overline{i}\, a} 
  - q^{-1} \hat{\sigma}_{\overline{i}\, a} \hat{\sigma}_{\overline{i+1}\, 
  \overline{i}} \notag \\
&= \hat{\sigma}_{\overline{i+1}\, \overline{i}} \hat{\sigma}_{\overline{i}\, a}
   - q^{-1} \hat{\sigma}_{\overline{i}\, a} \hat{\sigma}_{\overline{i+1}\, 
   \overline{i}}, \hspace{2cm} \varepsilon_a < -\varepsilon_i.
\end{align*}

\noindent \underline{Case 3: $a = \overline{i+1},\: b \neq i+1$}

\noindent  In this case we are unable to simplify the $q^{(\alpha_i,\varepsilon_b)}$ term of equation \eqref{Aibar}, so we have

\begin{equation*}
\hat{\sigma}_{b\,\overline{i}} = q^{(\alpha_i,\varepsilon_b)} \hat{\sigma}_{b\,
  \overline{i+1}} \hat{\sigma}_{\overline{i+1}\, \overline{i}} - q^{-1} 
  \hat{\sigma}_{\overline{i+1}\,\overline{i}} \hat{\sigma}_{b\,\overline{i+1}},
  \qquad \varepsilon_b > -  \varepsilon_{i+1},\; b \neq i+1.
\end{equation*}

\noindent \underline{Case 4: $b = i+1,\: a \neq \overline{i+1}$}

\noindent Again, the equation cannot be simplified further than

\begin{equation*}
\hat{\sigma}_{ia} = q^{-(\alpha_i, \varepsilon_a)} \hat{\sigma}_{i\, i+1} 
  \hat{\sigma}_{i+1\, a} - q^{-1} \hat{\sigma}_{i+1\,a} \hat{\sigma}_{i\, i+1},
  \quad \varepsilon_a < \varepsilon_{i+1},\; a \neq \overline{i+1}.
\end{equation*}

\noindent \underline{Case 5: $a = \overline{i+1},\, b=i+1$ }

\noindent This is the only case in which two of the terms on the left-hand side contribute.  We obtain

\begin{align*}
-\hat{\sigma}_{i+1\, \overline{i}} - \hat{\sigma}_{i\, \overline{i+1}} &=
  q^{-1} \hat{\sigma}_{i+1\, \overline{i+1}} \hat{\sigma}_{i\, i+1} - q^{-1}
  \hat{\sigma}_{i\, i+1} \hat{\sigma}_{i+1\, \overline{i+1}} \notag \\
\therefore \qquad \hat{\sigma}_{i+1\, \overline{i}} + \hat{\sigma}_{i\, \overline{i+1}} &= q^{-1}
   [\hat{\sigma}_{i\, i+1}, \hat{\sigma}_{i+1\, \overline{i+1}}].
\end{align*}

\noindent \underline{Case 6: $a \neq i, \overline{i+1}$ and $b \neq i+1, \overline{i}$}

\noindent Here the left-hand side of \eqref{Ai} vanishes, giving a commutation-style relation.

\begin{equation*}
q^{(\alpha_i,\varepsilon_b)} \hat{\sigma}_{ba} \hat{\sigma}_{i\, i+1} - q^{-(\alpha_i,\varepsilon_a)} \hat{\sigma}_{i\, i+1} \hat{\sigma}_{ba} = 0, \qquad
\varepsilon_b > \varepsilon_a.
\end{equation*}


\section{Relations for $\alpha_l = \varepsilon_l + \varepsilon_{l-1}$, where $m=2l$}

\noindent In this case $e_l \equiv E^{l-1}_{\overline{l}} - E^l_{\overline{l-1}}$.  Hence

\begin{equation*}
\langle a|e_l = \delta_{a\, l-1} \langle \overline{l}| - \delta_{al} \langle
  \overline{l-1}|, \qquad
e_l|b \rangle = \delta_{b \overline{l}} |l-1 \rangle - \delta_{b\, 
  \overline{l-1}} |l \rangle.
\end{equation*}

\noindent Thus equation (\ref{eqR}) becomes

\begin{multline*} 
q^{-\frac{1}{2} (\alpha_l, \alpha_l)} \bigl( \delta_{a\, l-1} 
  q^{\frac{1}{2} (\alpha_l, \varepsilon_{l-1})} \hat{\sigma}_{b \overline{l}}
  - \delta_{al} q^{\frac{1}{2} (\alpha_l, \varepsilon_l)} \hat{\sigma}_
  {b\, \overline{l-1}} \bigr) \\ 
- \bigl( \delta_{b \overline{l}} q^{-\frac{1}{2} (\alpha_l, \varepsilon_l)} 
  \hat{\sigma}_{l-1\, a} - \delta_{b\, \overline{l-1}} q^{-\frac{1}{2} 
  (\alpha_l, \varepsilon_{l-1})} \hat{\sigma}_{la} \bigr) \\
= q^{(\alpha_l, \varepsilon_b)} \hat{\sigma}_{ba} e_l q^{\frac{1}{2} h_l} - 
  q^{- (\alpha_l, \varepsilon_a)} e_l q^{\frac{1}{2} h_l} \hat{\sigma}_{ba}, 
  \quad \varepsilon_b > \varepsilon_a.
\end{multline*}

\noindent Noting that  $\hat{\sigma}_{l-1\, \overline{l}} = - \hat{\sigma}_{l\, \overline{l-1}} = q^{\frac{1}{2}} e_l q^{\frac{1}{2} h_l}$, this implies

\begin{align*} 
\delta_{a\, l-1} \hat{\sigma}_{b \overline{l}} - \delta_{al} \hat{\sigma}_
  {b\,\overline{l-1}} -  \delta_{b \overline{l}} \hat{\sigma}_{l-1\, a} + 
  &\delta_{b\, \overline{l-1}} \hat{\sigma}_{la}   \notag \\
&=q^{\frac{1}{2} + (\alpha_l, \varepsilon_b)} \hat{\sigma}_{ba} e_l 
  q^{\frac{1}{2} h_l} - q^{\frac{1}{2} - (\alpha_l, \varepsilon_a)} e_l 
  q^{\frac{1}{2} h_l} \hat{\sigma}_{ba} \notag \\
&=q^{(\alpha_l, \varepsilon_b)} \hat{\sigma}_{ba} \hat{\sigma}_{l-1\, 
  \overline{l}} - q^{-(\alpha_l, \varepsilon_a)} \hat{\sigma}_{l-1\,
  \overline{l}} \hat{\sigma}_{ba}, \quad \varepsilon_b >\varepsilon_a\\
&=q^{-(\alpha_l, \varepsilon_a)} \hat{\sigma}_{l\, \overline{l-1}} \hat{\sigma}
  _{ba} - q^{(\alpha_l, \varepsilon_b)} \hat{\sigma}_{ba} \hat{\sigma}_{l\, 
  \overline{l-1}}, \quad \varepsilon_b > \varepsilon_a.
\end{align*}

\noindent To extract all the information contained in these equations we consider the various cases separately, as in the previous section.  The relations we obtain and the conditions under which they hold are:

\begin{alignat*}{2}
&\hat{\sigma}_{b\, \overline{l-1}} = q^{(\alpha_l, \varepsilon_b)} \hat{\sigma}
  _{bl} \hat{\sigma}_{l\, \overline{l-1}} - q^{-1} \hat{\sigma}_{l\, 
  \overline{l-1}} \hat{\sigma}_{bl}, && \varepsilon_b > \varepsilon_l, \\
&\hat{\sigma}_{b \overline{l}} = \hat{\sigma}_{b\,l-1} \hat{\sigma}_{l-1 \, 
  \overline{l}} - q^{-1} \hat{\sigma}_{l-1\,\overline{l}} \hat{\sigma}_{b\,l-1}
  , && \varepsilon_b > \varepsilon_{l-1}, \\
&\hat{\sigma}_{la} =\hat{\sigma}_{l\, \overline{l-1}} \hat{\sigma}_{\overline
  {l-1}\,a} - q^{-1} \hat{\sigma}_{\overline{l-1}\, a} \hat{\sigma}_{l\, 
  \overline{l-1}}, && \varepsilon_a < -\varepsilon_{l-1}, \\
&\hat{\sigma}_{l-1\, a} = q^{-(\alpha_l, \varepsilon_a)} \hat{\sigma}_{l-1\, 
  \overline{l}} \hat{\sigma}_{\overline{l}a} - q^{-1} \hat{\sigma}_{\overline
  {l}a}\hat{\sigma}_{l-1\,\overline{l}},&\qquad&\varepsilon_a<-\varepsilon_l,\\
&q^{(\alpha_l, \varepsilon_b)} \hat{\sigma}_{ba} \hat{\sigma}_{l-1\, 
  \overline{l}} - q^{- (\alpha_l, \varepsilon_a)} \hat{\sigma}_{l-1\, 
  \overline{l}} \hat{\sigma}_{ba} = 0,&& \varepsilon_b > \varepsilon_a,\; a 
  \neq l,\, l-1 \text{ and } b \neq \overline{l-1},\, \overline{l}.
\end{alignat*}


\section{Relations for $\alpha_l = \varepsilon_l$, where $m=2l+1$}

\noindent Here $e_l \equiv E^l_{l+1} - E^{l+1}_{\overline{l}}$, and thus

\begin{equation*}
\langle a|e_l = \delta_{al} \langle l+1| - \delta_{a\,l+1} \langle 
  \overline{l}|, \qquad
e_l|b \rangle = \delta_{b\, l+1} |l \rangle - \delta_{b \overline{l}} |l+1
  \rangle.
\end{equation*}

\noindent Applying these to (\ref{eqR}) gives

\begin{multline*}
q^{-\frac{1}{2}(\alpha_l,\alpha_l)} \bigl( \delta_{al} q^{\frac{1}{2} 
  (\alpha_l, \varepsilon_l)} \hat{\sigma}_{b\, l+1} - \delta_{a\,l+1} 
  q^{\frac{1}{2} (\alpha_l, \varepsilon_{l+1})} \hat{\sigma}_{b \overline{l}} 
  \bigr) \\
- \bigl( \delta_{b\,l+1} q^{\frac{1}{2}(\alpha_l,\varepsilon_{l+1})}
  \hat{\sigma}_{la} - \delta_{b \overline{l}} q^{-\frac{1}{2} (\alpha_l, 
  \varepsilon_l)} \hat{\sigma}_{l+1\, a} \bigr) \\
= q^{(\alpha_l, \varepsilon_b)} \hat{\sigma}_{ba} e_l q^{\frac{1}{2}h_l} - 
  q^{-(\alpha_l, \varepsilon_a)} e_l q^{\frac{1}{2} h_l} \hat{\sigma}_{ba}, 
  \quad \varepsilon_b > \varepsilon_a.
\end{multline*}

\noindent Recalling that $\hat{\sigma}_{l\, l+1} = - q^{-\frac{1}{2}} \hat{\sigma}_{l+1\, \overline{l}} = e_l q^{\frac{1}{2} h_l}$, we find

\begin{alignat*}{2} 
q^{\frac{1}{2}} \delta_{al} \hat{\sigma}_{b\, l+1} - \delta_{a\,l+1} 
  \hat{\sigma}_{b \overline{l}} - q^{\frac{1}{2}} &\delta_{b\, l+1} 
  \hat{\sigma}_{la} + \delta_{b\overline{l}} \hat{\sigma}_{l+1\, a}&& \notag \\
&=q^{\frac{1}{2} + (\alpha_l, \varepsilon_b)} \hat{\sigma}_{ba} e_l 
  q^{\frac{1}{2}h_l} - q^{\frac{1}{2} - (\alpha_l, \varepsilon_a)} e_l 
  q^{\frac{1}{2} h_l} \hat{\sigma}_{ba}&& \notag  \\
&=q^{\frac{1}{2} + (\alpha_l, \varepsilon_b)} \hat{\sigma}_{ba} \hat{\sigma}_
  {l\, l+1} - q^{\frac{1}{2} - (\alpha_l,\varepsilon_a)} \hat{\sigma}_{l\, l+1}
  \hat{\sigma}_{ba}, &\quad &\varepsilon_b > \varepsilon_a \\
&=q^{- (\alpha_l,\varepsilon_a)} \hat{\sigma}_{l+1\, \overline{l}} \hat{\sigma}
  _{ba} - q^{(\alpha_l, \varepsilon_b)} \hat{\sigma}_{ba} \hat{\sigma}_{l+1\, 
  \overline{l}}, && \varepsilon_b > \varepsilon_a.
\end{alignat*}

\noindent By examining the various cases we deduce the following relations:

\begin{alignat*}{2}
&\hat{\sigma}_{b\,l+1} =  \hat{\sigma}_{bl} \hat{\sigma}_{l\,l+1} - q^{-1}
  \hat{\sigma}_{l\,l+1} \hat{\sigma}_{bl},&& \varepsilon_b > \varepsilon_l, \\
&\hat{\sigma}_{b\, \overline{l}} = q^{(\alpha_l, \varepsilon_b)} \hat{\sigma}
  _{b\, l+1} \hat{\sigma}_{l+1\,\overline{l}} -\hat{\sigma}_{l+1\,\overline{l}}
  \hat{\sigma}_{b\, l+1}, && \varepsilon_b > 0, \\
&\hat{\sigma}_{la} = q^{-(\alpha_l, \varepsilon_a)} \hat{\sigma}_{l\, l+1} 
  \hat{\sigma}_{l+1\, a} - \hat{\sigma}_{l+1\, a} \hat{\sigma}_{l\, l+1}, 
  && \varepsilon_a < 0, \\
&\hat{\sigma}_{l+1\, a} = \hat{\sigma}_{l+1\, \overline{l}} \hat{\sigma}_
  {\overline{l}a} - q^{-1} \hat{\sigma}_{\overline{l}a} \hat{\sigma}_{l+1\, 
  \overline{l}}, && \varepsilon_a < -\varepsilon_l, \\
&q^{(\alpha_l, \varepsilon_b)} \hat{\sigma}_{ba} \hat{\sigma}_{l\, l+1} - 
  q^{-(\alpha_l, \varepsilon_a)} \hat{\sigma}_{l\, l+1} \hat{\sigma}_{ba} = 0, 
  &\quad& \varepsilon_b > \varepsilon_a,\: a \neq l,\, l+1 \text{ and } b \neq 
  l+1,\, \overline{l}.
\end{alignat*}


\section{Relations for $\alpha_\mu = \delta_\mu - \delta_{\mu+1}$, $1 \leq \mu < k$}

In this case $e_\mu \equiv E^\mu_{\mu+1} + E^{\overline{\mu+1}}_{\overline{\mu}}$.  Therefore

\begin{equation*}
\langle a|e_\mu = \delta_{a\mu} \langle \mu+1| + \delta_{a\,\overline{\mu+1}} 
  \langle \overline{\mu}|, \qquad
e_\mu|b \rangle = \delta_{b\, \mu+1} |\mu \rangle + \delta_{b \overline{\mu}} 
|\overline{\mu+1} \rangle.
\end{equation*}

\noindent Applying this to (\ref{eqR}) produces

\begin{multline*}
q^{-\frac{1}{2} (\alpha_\mu, \alpha_\mu)} \bigl( \delta_{a\mu} q^{\frac{1}{2} 
  (\alpha_\mu, \delta_\mu)} \hat{\sigma}_{b\, \mu+1} + \delta_{a\,\overline
  {\mu+1}} q^{-\frac{1}{2} (\alpha_\mu, \delta_{\mu+1})} \hat{\sigma}_{b 
  \overline{\mu}} \bigr) \\
- \bigl( \delta_{b\, \mu+1} q^{\frac{1}{2} (\alpha_\mu, \delta_{\mu+1})} 
  \hat{\sigma}_{\mu a} + \delta_{b \overline{\mu}} q^{-\frac{1}{2} (\alpha_\mu,
  \delta_{\mu})} \hat{\sigma}_{\overline{\mu+1}\, a} \bigr) \\
= q^{ (\alpha_\mu,\varepsilon_b)} \hat{\sigma}_{ba} e_\mu q^{\frac{1}{2} h_\mu}
  - q^{-(\alpha_\mu,\varepsilon_a)} e_\mu q^{\frac{1}{2} h_\mu} 
  \hat{\sigma}_{ba}, \quad \varepsilon_b > \varepsilon_a.
\end{multline*}

\noindent However $\hat{\sigma}_{\mu \, \mu +1} = \hat{\sigma}_{\overline{\mu+1} \, \overline{\mu}} = q^{-\frac{1}{2}} e_\mu q^{\frac{1}{2} h_\mu}$.  Thus this equation reduces to

\begin{alignat*}{2}
\delta_{a\mu} \hat{\sigma}_{b\, \mu+1} + \delta_{a\,\overline{\mu+1}} 
  \hat{\sigma}_{b \overline{\mu}} - &\delta_{b\, \mu+1} \hat{\sigma}_{\mu a} - 
  \delta_{b \overline{\mu}} \hat{\sigma}_{\overline{\mu+1}\, a}&& \notag \\
&=q^{-\frac{1}{2}} \bigl( q^{(\alpha_\mu,\varepsilon_b)} \hat{\sigma}_{ba} 
  e_\mu q^{\frac{1}{2} h_\mu} - q^{-(\alpha_\mu,\varepsilon_a)} e_\mu 
  q^{\frac{1}{2} h_\mu} \hat{\sigma}_{ba} \bigr) &&\notag \\
&=q^{(\alpha_\mu,\varepsilon_b)} \hat{\sigma}_{ba} \hat{\sigma}_{\mu \, \mu +1}
  - q^{-(\alpha_\mu,\varepsilon_a)} \hat{\sigma}_{\mu \, \mu +1} \hat{\sigma}_
  {ba},&\quad& \varepsilon_b > \varepsilon_a \\
&=q^{(\alpha_\mu,\varepsilon_b)} \hat{\sigma}_{ba} \hat{\sigma}_{\overline
  {\mu+1} \, \overline{\mu}} - q^{-(\alpha_\mu,\varepsilon_a)} \hat{\sigma}_
  {\overline{\mu+1} \, \overline{\mu}} \hat{\sigma}_{ba},&& \varepsilon_b > 
  \varepsilon_a .
\end{alignat*}

\noindent From these we obtain:

\begin{alignat*}{2}
&\hat{\sigma}_{\nu\, \mu+1} = \hat{\sigma}_{\nu \mu} \hat{\sigma}_{\mu\, \mu+1}
  - q \hat{\sigma}_{\mu\, \mu+1}  \hat{\sigma}_{\nu \mu}, && \nu < \mu<k, \\
&\hat{\sigma}_{\overline{\mu+1}\, \overline{\nu}} = \hat{\sigma}_{\overline
  {\mu+1} \, \overline{\mu}} \hat{\sigma}_{\overline{\mu} \overline{\nu}} - q 
  \hat{\sigma}_{\overline {\mu} \overline{\nu}} \hat{\sigma}_{\overline{\mu+1} 
  \,\overline{\mu}}, & \quad & \nu < \mu < k ,\\
&\hat{\sigma}_{b \overline{\mu}} = q^{(\alpha_\mu, \varepsilon_b)} \hat{\sigma}
  _{b\, \overline{\mu+1}} \hat{\sigma}_{\overline{\mu+1}\, \overline{\mu}} - q
  \hat{\sigma}_{\overline{\mu+1}\, \overline{\mu}} \hat{\sigma}_{b\, \overline
  {\mu+1}}, & \quad & \varepsilon_b > - \delta_{\mu + 1}, b \neq \mu + 1, \\
&\hat{\sigma}_{\mu a} = q^{-(\alpha_\mu, \varepsilon_a)} \hat{\sigma}_{\mu\, 
  \mu+1} \hat{\sigma}_{\mu + 1\, a} - q \hat{\sigma}_{\mu + 1\, a} 
  \hat{\sigma}_{\mu\, \mu+1},&& \varepsilon_a < \delta_{\mu + 1}, a \neq 
  \overline{\mu + 1}, \\
&\hat{\sigma}_{\mu+1\, \overline{\mu}} - \hat{\sigma}_{\mu\,\overline{\mu +1}} 
  = q [ \hat{\sigma}_{\mu + 1\,\overline{\mu +1}},\hat{\sigma}_{\mu\, \mu+1}],
  \\
&q^{(\alpha_\mu,\varepsilon_b)} \hat{\sigma}_{ba} \hat{\sigma}_{\mu \, \mu +1}
  - q^{-(\alpha_\mu,\varepsilon_a)} \hat{\sigma}_{\mu \, \mu +1} \hat{\sigma}_
  {ba} = 0, && \varepsilon_b > \varepsilon_a,\,a \neq \mu, \overline{\mu+1} 
  \text{ and } b \neq \mu+1, \overline{\mu}.
\end{alignat*}


\section{Relations for $\alpha_s = \delta_k - \varepsilon_1$}

\noindent The raising generator corresponding to $\alpha_s$ is $e_s \equiv E^{\mu =k}_{i=1} + (-1)^k E^{i = \overline{1}}_{\mu = \overline{k}}$.  Hence

\begin{align*}
\langle a|e_s &= \delta_{a\,\mu=k} \langle i=1| + (-1)^k \delta_{a\, 
  i= \overline{1}} \langle \mu = \overline{k}|, \\
e_s|b \rangle &= \delta_{b\, i=1} |\mu =k \rangle + (-1)^k \delta_{b\, \mu = 
  \overline{k}} |i = \overline{1} \rangle.
\end{align*}

\noindent Substituting these into  equation (\ref{eqR}) produces

\begin{multline*}
q^{-\frac{1}{2} (\alpha_s, \alpha_s)} \bigl( \delta_{a\,\mu=k} q^{\frac{1}{2} 
  (\alpha_s, \delta_k)} \hat{\sigma}_{b\, i=1}+(-1)^k\delta_{a\,i=\overline{1}}
  q^{-\frac{1}{2} (\alpha_s,\varepsilon_1)} \hat{\sigma}_{b\,\mu=\overline{k}}
  \bigr) \\
- (-1)^{[a]} \bigl( \delta_{b\, i=1} q^{\frac{1}{2} (\alpha_s, \varepsilon_1)}
  \hat{\sigma}_{\mu=k\, a} - (-1)^k \delta_{b\, \mu=\overline{k}} q^{-\frac
  {1}{2} (\alpha_s, \delta_k)} \hat{\sigma}_{i=\overline{1}\, a} \bigr) \\
= q^{ (\alpha_s,\varepsilon_b)} \hat{\sigma}_{ba} e_s q^{\frac{1}{2} h_s}
  - (-1)^{[a]+[b]} q^{-(\alpha_s,\varepsilon_a)} e_s q^{\frac{1}{2} h_s} 
  \hat{\sigma}_{ba}, \quad \varepsilon_b > \varepsilon_a.
\end{multline*}

\noindent  Using $\hat{\sigma}_{\mu=k \, i=1} = (-1)^k q\, \hat{\sigma}_{i = \overline{1} \, \mu = \overline{k}} = q^{\frac{1}{2}} e_s q^{\frac{1}{2} h_s}$ and simplifying, we obtain

\begin{alignat*}{2}
\bigl( \delta_{a\,\mu=k} \hat{\sigma}_{b\, i=1} + &(-1)^k q \delta_{a\,
  i=\overline{1}} \hat{\sigma}_{b\,\mu=\overline{k}} \bigr) - (-1)^{[a]} 
  \bigl( \delta_{b\, i=1}  \hat{\sigma}_{\mu=k\, a}- (-1)^k q \delta_{b\, \mu=
  \overline{k}} \hat{\sigma}_{i=\overline{1}\, a} \bigr) && \notag\\ 
&= q^{\frac{1}{2}} \bigl( q^{ (\alpha_s,\varepsilon_b)} \hat{\sigma}_{ba} e_s 
  q^{\frac{1}{2} h_s} - (-1)^{[a]+[b]} q^{-(\alpha_s,\varepsilon_a)} e_s 
  q^{\frac{1}{2} h_s} \hat{\sigma}_{ba} \bigr) && \notag \\
&= q^{(\alpha_s, \varepsilon_b)} \hat{\sigma}_{ba} \hat{\sigma}_{\mu=k \, i=1}
  - (-1)^{[a]+[b]} q^{-(\alpha_s,\varepsilon_a)} \hat{\sigma}_{\mu=k \, i=1}
  \hat{\sigma}_{ba} && \varepsilon_b > \varepsilon_a \\
&= (-1)^k q \bigl( q^{(\alpha_s, \varepsilon_b)} \hat{\sigma}_{ba} 
  \hat{\sigma}_{i = \overline{1} \, \mu = \overline{k}} - (-1)^{[a]+[b]} 
  q^{-(\alpha_s, \varepsilon_a)} \hat{\sigma}_{i = \overline{1} \, \mu = 
  \overline{k}} \hat{\sigma}_{ba} \bigr), && \varepsilon_b > \varepsilon_a.
\end{alignat*}

\noindent From these equations the following relations can be deduced:

\begin{alignat*}{2}
&\hat{\sigma}_{\nu\, i=1} = \hat{\sigma}_{\nu\, \mu=k} \hat{\sigma}_{\mu=k \, 
  i=1} - q \hat{\sigma}_{\mu=k \, i=1} \hat{\sigma}_{\nu\, \mu=k}, & \quad
  &\nu < k, \\
&\hat{\sigma}_{i=\overline{1}\, \overline{\nu}} = \hat{\sigma}_{i=\overline{1}
  \, \mu = \overline{k}} \hat{\sigma}_{\mu = \overline{k}\, \overline{\nu}}
  - q \hat{\sigma}_{\mu = \overline{k}\, \overline{\nu}} \hat{\sigma}_{i = 
  \overline{1}\, \mu = \overline{k}}, && \nu < k, \\
&\hat{\sigma}_{\mu=k\, a} = q^{-(\alpha_s, \varepsilon_a)} \hat{\sigma}_{\mu=k
  \,i=1} \hat{\sigma}_{i=1\, a} - (-1)^{[a]} q^{-1} \hat{\sigma}_{i=1\, a} \hat
  {\sigma}_{\mu=k \,i=1},&& \varepsilon_a < \varepsilon_1, \varepsilon_a \neq  
  -\varepsilon_1, \\
&\hat{\sigma}_{b\, \mu = \overline{k}} = q^{(\alpha_s, \varepsilon_b)} 
  \hat{\sigma}_{b\, i= \overline{1}} \hat{\sigma}_{i = \overline{1}\, \mu = 
  \overline{k}} - (-1)^{[b]} q^{-1}  \hat{\sigma}_{i = \overline{1}\, \mu = 
  \overline{k}} \hat{\sigma}_{b\, i= \overline{1}}, && \varepsilon_b > - 
  \varepsilon_1, \varepsilon_b \neq \varepsilon_1, \\
&\hat{\sigma}_{\mu=k\, \overline{i} = \overline{1}} - (-1)^k q \hat{\sigma}_
  {i=1\, \overline{\mu} = \overline{k}}  = q^{-1} [ \hat{\sigma}_{\mu=k\, i=1},
  \hat{\sigma}_{i=1\, \overline{i} = \overline{1}}], \\
&q^{(\alpha_s, \varepsilon_b)} \hat{\sigma}_{ba} \hat{\sigma}_{\mu=k \, i=1}
  - (-1)^{[a]+[b]} q^{-(\alpha_s,\varepsilon_a)} \hat{\sigma}_{\mu=k \, i=1}
  \hat{\sigma}_{ba} = 0, && \varepsilon_b > \varepsilon_a, \varepsilon_a \neq 
  \delta_k, -\varepsilon_1 \\
&&&\text{ and } \varepsilon_b \neq \varepsilon_1, -\delta_k.
\end{alignat*}


\section{Summary of Relations}

The relations arising from \eqref{**} common to the $m=2l$ and $m=2l+1$ cases are listed in Table \ref{list} on the following page.  Tables \ref{even} and \ref{odd}, on the page after that, contain the extra relations applying only in the cases of even and odd values of $m$ respectively.  The $q$-commutation relations \eqref{qcom} and inductive relations \eqref{indrel} were deduced from these lists.

\begin{table} 
\caption{The relations for the operators $\hat{\sigma}_{ba}$ common to all values of $m$} \label{list}
\centering
\begin{tabular}{|ll|} \hline
\hspace{25mm} Relation & \hspace{10mm} Conditions \\ \hline 

$\hat{\sigma}_{b\, i+1} = \hat{\sigma}_{bi} \hat{\sigma}_{i\, i+1} - q^{-1} 
  \hat{\sigma}_{i\, i+1} \hat{\sigma}_{bi},$&$i<l,\: \varepsilon_b > 
  \varepsilon_i$ \\
$\hat{\sigma}_{\overline{i+1}\,a} = \hat{\sigma}_{\overline{i+1}\,\overline{i}}
  \hat{\sigma}_{\overline{i}\, a} - q^{-1} \hat{\sigma}_{\overline{i}\, a} 
  \hat{\sigma}_{\overline{i+1}\,\overline{i}}, $&$ i<l,\: \varepsilon_a < - 
  \varepsilon_i$\\
$\hat{\sigma}_{b\,\overline{i}} = q^{(\alpha_i, \varepsilon_b)} \hat{\sigma}_
  {b\, \overline{i+1}} \hat{\sigma}_{\overline{i+1}\, \overline{i}} - q^{-1} 
  \hat{\sigma}_{\overline{i+1}\,\overline{i}} \hat{\sigma}_{b\,\overline{i+1}},
  $&$ i<l,\: \varepsilon_b > - \varepsilon_{i+1},$ \\
&$ b \neq i+1$ \\
$\hat{\sigma}_{ia} = q^{-(\alpha_i, \varepsilon_a)} \hat{\sigma}_{i\, i+1} 
  \hat{\sigma}_{i+1\, a} - q^{-1} \hat{\sigma}_{i+1\, a} \hat{\sigma}_{i\,i+1},
  $&$ i<l,\: \varepsilon_a < \varepsilon_{i+1},$ \\&$ a \neq \overline{i+1}$ \\
$\hat{\sigma}_{i+1\,\overline{i}} + \hat{\sigma}_{i\, \overline{i+1}} = q^{-1}
   [\hat{\sigma}_{i\, i+1}, \hat{\sigma}_{i+1\, \overline{i+1}}],$& $i<l$\\
$q^{(\alpha_i,\varepsilon_b)} \hat{\sigma}_{ba} \hat{\sigma}_{i\, i+1} - 
  q^{-(\alpha_i,\varepsilon_a)} \hat{\sigma}_{i\, i+1} \hat{\sigma}_{ba} = 0,$
  & $i<l;\: \varepsilon_b > \varepsilon_a;$\\&$ a \neq i, \overline{i+1}$ \\
& and $b \neq i+1, \overline{i}$ \\
$\hat{\sigma}_{\nu\, \mu+1} = \hat{\sigma}_{\nu \mu} \hat{\sigma}_{\mu\, \mu+1}
  - q \hat{\sigma}_{\mu\, \mu+1}  \hat{\sigma}_{\nu \mu}, $&$ \nu < \mu<k$ \\
$\hat{\sigma}_{\overline{\mu+1}\, \overline{\nu}} = \hat{\sigma}_{\overline
  {\mu+1} \, \overline{\mu}} \hat{\sigma}_{\overline{\mu} \overline{\nu}} - q 
  \hat{\sigma}_{\overline {\mu} \overline{\nu}} \hat{\sigma}_{\overline{\mu+1} 
  \,\overline{\mu}}, $&$ \nu < \mu < k$ \\
$\hat{\sigma}_{b \overline{\mu}} = q^{(\alpha_\mu, \varepsilon_b)} \hat{\sigma}
  _{b\, \overline{\mu+1}} \hat{\sigma}_{\overline{\mu+1}\, \overline{\mu}} - q
  \hat{\sigma}_{\overline{\mu+1}\, \overline{\mu}} \hat{\sigma}_{b\, \overline
  {\mu+1}}, $&$ \mu<k,\: \varepsilon_b > -\delta_{\mu + 1},$\\&$b\neq \mu +1$\\
$\hat{\sigma}_{\mu a} = q^{-(\alpha_\mu, \varepsilon_a)} \hat{\sigma}_{\mu\, 
  \mu+1} \hat{\sigma}_{\mu + 1\, a} - q \hat{\sigma}_{\mu + 1\, a} 
  \hat{\sigma}_{\mu\, \mu+1},$&$ \mu<k,\: \varepsilon_a < \delta_{\mu + 1},$\\
&$a \neq \overline{\mu + 1}$ \\
$\hat{\sigma}_{\mu+1\, \overline{\mu}} - \hat{\sigma}_{\mu\,\overline{\mu +1}} 
  = q [ \hat{\sigma}_{\mu +1 \, \overline{\mu +1}},\hat{\sigma}_{\mu\,\mu+1}],$
  & $\mu<k$ \\
$q^{(\alpha_\mu,\varepsilon_b)} \hat{\sigma}_{ba} \hat{\sigma}_{\mu \, \mu +1}
  - q^{-(\alpha_\mu,\varepsilon_a)} \hat{\sigma}_{\mu \, \mu +1} \hat{\sigma}_
  {ba} = 0, $&$\mu<k;\: \varepsilon_b > \varepsilon_a;$ \\
&$a\neq \mu, \overline{\mu+1}$ \\ &and $b \neq \mu+1, \overline{\mu}$ \\
$\hat{\sigma}_{\nu\, i=1} = \hat{\sigma}_{\nu\, \mu=k} \hat{\sigma}_{\mu=k \, 
  i=1} - q \hat{\sigma}_{\mu=k \, i=1} \hat{\sigma}_{\nu\, \mu=k},$&$\nu <k$\\
$\hat{\sigma}_{i=\overline{1}\, \overline{\nu}} = \hat{\sigma}_{i=\overline{1}
  \, \mu = \overline{k}} \hat{\sigma}_{\mu = \overline{k}\, \overline{\nu}}
  - q \hat{\sigma}_{\mu = \overline{k}\, \overline{\nu}} \hat{\sigma}_{i = 
  \overline{1}\, \mu = \overline{k}}, $&$ \nu < k $\\
$\hat{\sigma}_{\mu=k\, a} = q^{-(\alpha_s, \varepsilon_a)} \hat{\sigma}_{\mu=k
  \,i=1} \hat{\sigma}_{i=1\, a} - (-1)^{[a]} q^{-1} \hat{\sigma}_{i=1\, a} \hat
  {\sigma}_{\mu=k \,i=1},$&$ \varepsilon_a < \varepsilon_1,\: a\neq i=\overline
  {1} $\\
$\hat{\sigma}_{b\, \mu = \overline{k}} = q^{(\alpha_s, \varepsilon_b)} 
  \hat{\sigma}_{b\, i= \overline{1}} \hat{\sigma}_{i = \overline{1}\, \mu = 
  \overline{k}} - (-1)^{[b]} q^{-1}  \hat{\sigma}_{i = \overline{1}\, \mu = 
  \overline{k}} \hat{\sigma}_{b\, i= \overline{1}}, $&$ \varepsilon_b > - 
  \varepsilon_1,\: b \neq i = 1$ \\
$\hat{\sigma}_{\mu=k\, \overline{i} = \overline{1}} - (-1)^k q \hat{\sigma}_
  {i=1\, \overline{\mu} = \overline{k}}  = q^{-1} [ \hat{\sigma}_{\mu=k\, i=1},
  \hat{\sigma}_{i=1\, \overline{i} = \overline{1}}],$& \\
$q^{(\alpha_s, \varepsilon_b)} \hat{\sigma}_{ba} \hat{\sigma}_{\mu=k \, i=1}
  - (-1)^{[a]+[b]} q^{-(\alpha_s,\varepsilon_a)} \hat{\sigma}_{\mu=k \, i=1}
  \hat{\sigma}_{ba} = 0 $&$ \varepsilon_b > \varepsilon_a;\: \varepsilon_a \neq
  \delta_k, -\varepsilon_1$ \\
&and $\varepsilon_b \neq \varepsilon_1, -\delta_k$ \\ \hline
\end{tabular}
\end{table}

\begin{table}
\caption{The relations for the operators $\hat{\sigma}_{ba}$ that hold only for even $m$} \label{even}
\centering
\begin{tabular}{|ll|} \hline
\hspace{20mm} Relation & Conditions \\ \hline

$\hat{\sigma}_{b\, \overline{l-1}} = q^{(\alpha_l,\varepsilon_b)} \hat{\sigma}
  _{bl} \hat{\sigma}_{l\, \overline{l-1}} - q^{-1} \hat{\sigma}_{l\, \overline
  {l-1}} \hat{\sigma}_{bl},$&$ \varepsilon_b > \varepsilon_l $\\
$\hat{\sigma}_{b \overline{l}} = \hat{\sigma}_{b\,l-1} \hat{\sigma}_{l-1 \, 
  \overline{l}} - q^{-1} \hat{\sigma}_{l-1\,\overline{l}} \hat{\sigma}_{b\,l-1}
  ,$&$ \varepsilon_b > \varepsilon_{l-1} $\\
$\hat{\sigma}_{la} = \hat{\sigma}_{l\, \overline{l-1}} \hat{\sigma}_{\overline
  {l-1}\, a} - q^{-1} \hat{\sigma}_{\overline{l-1}\, a} \hat{\sigma}_{l\, 
  \overline{l-1}},$&$ \varepsilon_a < -\varepsilon_{l-1} $\\
$\hat{\sigma}_{l-1\, a} = q^{-(\alpha_l, \varepsilon_a)} \hat{\sigma}_{l-1\, 
  \overline{l}}\hat{\sigma}_{\overline{l}a} -q^{-1}\hat{\sigma}_{\overline{l}a}
  \hat{\sigma}_{l-1\,\overline{l}},$&$ \varepsilon_a < -\varepsilon_l$\\ 
$q^{(\alpha_l, \varepsilon_b)} \hat{\sigma}_{ba} \hat{\sigma}_{l-1\, 
  \overline{l}}- q^{-(\alpha_l,\varepsilon_a)} \hat{\sigma}_{l-1\,\overline{l}}
  \hat{\sigma}_{ba} = 0,$&$ \varepsilon_b > \varepsilon_a;\: a \neq l, l-1;\: 
  b \neq \overline{l-1}, \overline{l} $ \\ \hline
\end{tabular}
\end{table}

\begin{table}
\caption{The relations for the operators $\hat{\sigma}_{ba}$ that hold only for odd $m$} \label{odd}
\centering
\begin{tabular}{|ll|} \hline
\hspace{20mm} Relation & Conditions \\ \hline
$\hat{\sigma}_{b\,l+1} =  \hat{\sigma}_{bl} \hat{\sigma}_{l\,l+1} - q^{-1}
  \hat{\sigma}_{l\,l+1} \hat{\sigma}_{bl},$&$ \varepsilon_b > \varepsilon_l$\\
$\hat{\sigma}_{b\, \overline{l}} = q^{(\alpha_l, \varepsilon_b)} \hat{\sigma}
  _{b\,l+1} \hat{\sigma}_{l+1\,\overline{l}} -\hat{\sigma}_{l+1\, \overline{l}}
  \hat{\sigma}_{b\, l+1},$&$ \varepsilon_b > 0$\\
$\hat{\sigma}_{la} = q^{-(\alpha_l, \varepsilon_a)} \hat{\sigma}_{l\, l+1} 
  \hat{\sigma}_{l+1\, a} - \hat{\sigma}_{l+1\, a} \hat{\sigma}_{l\, l+1}, 
  $&$ \varepsilon_a < 0 $\\
$\hat{\sigma}_{l+1\, a} = \hat{\sigma}_{l+1\, \overline{l}} \hat{\sigma}_
  {\overline{l}a} - q^{-1} \hat{\sigma}_{\overline{l}a} \hat{\sigma}_{l+1\, 
  \overline{l}}, $&$\varepsilon_a < -\varepsilon_l $\\ 
$q^{(\alpha_l, \varepsilon_b)} \hat{\sigma}_{ba} \hat{\sigma}_{l\, l+1} - 
  q^{-(\alpha_l, \varepsilon_a)} \hat{\sigma}_{l\,l+1} \hat{\sigma}_{ba} = 0,$
  &$\varepsilon_b >\varepsilon_a;\: a \neq l, l+1;\: b\neq l+1,\overline{l}$\\ 
\hline
\end{tabular}
\end{table}


\end{document}